\documentclass[10pt]{amsart}

\usepackage{eucal}
\usepackage{amscd}
\usepackage[all,cmtip]{xy}
\usepackage{rotating}
\usepackage{hyperref,mathrsfs}
\usepackage[listings]{tcolorbox}
\usepackage{tikz}
\usetikzlibrary{positioning,decorations.markings,arrows}
\usepackage{amsmath,amssymb,soul,minibox}
\usepackage{multicol}
\usepackage{tcolorbox}

\setul{}{1.3pt}
 \renewcommand\fbox{\fcolorbox{black}{white}}

\newtheorem{theorem}{Theorem }[section]

\newtheorem{lemma}[theorem]{Lemma}
\newtheorem{conjecture}[theorem]{Conjecture}

\newtheorem{corollary}[theorem]{Corollary}
\newtheorem{proposition}[theorem]{Proposition}
\newtheorem{question}[theorem]{Question}
\theoremstyle{definition}

\newtheorem{remark}[theorem]{Remark}
\newtheorem{observation}[theorem]{Observation}

\newcommand{\vb}{\subsubsection*{Example.}}
\newcommand{\vbn}{\subsubsection*{Examples}}
\newcommand{\defi}{\subsubsection*{Definition}}
\newcommand{\defs}{\subsubsection*{Definitions}}
\newcommand{\notation}{\subsubsection*{Notation}}

\newcommand{\proj}{\mathrm{proj}}
\def\I{\mathbf{I}}
\def\nI{\mathop{\not\mathrm{I}}}
\def\para{\mathop{\parallel}}
\def\notpara{\mathop{\not\|}}
\def\ultrapara{\mathop{\widehat{\parallel}}}
\def\spg{\mathop{\mathrm{spg}}\nolimits}
\def\pg{\mathop{\mathrm{pg}}\nolimits}
\def\srg{\mathop{\mathrm {srg}}\nolimits}
\def\min{\mathop{\mathrm{min}}\nolimits}
\def\PG{\mathbf{PG}}
\def\AG{\mathrm{AG}}
\def\GF{\mathrm{GF}}
\def\net{\mathcal{N}}
\def\eop{\hspace*{\fill}{\footnotesize$\blacksquare$}}
\newcommand{\Aut}{\mathrm{Aut}}
\newcommand{\id}{\mathrm{id}}

\newcommand{\wis}[1]{{\text{\em \usefont{OT1}{cmtt}{m}{n} #1}}}

\newcommand{\mO}{\mathcal{O}}
\newcommand{\A}{\mathbb{A}}
\newcommand{\hP}{\mathbb{P}}
\newcommand{\B}{\overline{\mathbf{B}}}
\newcommand{\cB}{\mathcal{B}}
\newcommand{\mD}{\mathcal{D}}
\newcommand{\Fun}{\mathbb{F}_1}
\newcommand{\Spec}{\wis{Spec}}
\newcommand{\USpec}{\underline{\wis{Spec}}}
\newcommand{\Z}{\mathbb{Z}}
\newcommand{\Proj}{\wis{Proj}}
\newcommand{\mE}{\mathcal{E}}
\newcommand{\bG}{\mathbb{G}}
\newcommand{\Hom}{\wis{Hom}}
\newcommand{\Adj}{\wis{Adj}}
\newcommand{\Inc}{\wis{Inc}}
\newcommand{\Gal}{\wis{Gal}}

\def\doubleprod#1#2{\ooalign{$#1\prod$\cr$#1\coprod$\cr}}
\DeclareMathOperator*{\Rprod}{\mathpalette\doubleprod\relax}

\newcommand{\fp}{\frak{p}}
\newcommand{\rd}{\mathrm{d}}

\newcommand{\mC}{\mathcal{C}}
\newcommand{\bbL}{\mathbb{L}}
\newcommand{\mF}{\mathcal{F}}
\newcommand{\cJ}{\mathcal{J}}
\newcommand{\mK}{\mathcal{K}}
\newcommand{\mL}{\mathcal{L}}
\newcommand{\cM}{\mathcal{M}}
\newcommand{\cN}{\mathcal{N}}
\newcommand{\cO}{\mathcal{O}}
\newcommand{\mP}{\mathcal{P}}
\newcommand{\fP}{\mathbf{P}}
\newcommand{\cQ}{\mathcal{Q}}
\newcommand{\mM}{\mathcal{M}}
\newcommand{\mS}{\mathcal{S}}
\newcommand{\cT}{\mathcal{T}}
\newcommand{\cU}{\mathcal{U}}
\newcommand{\cP}{\mathcal{P}}
\newcommand{\hW}{\mathbf{W}}
\newcommand{\mX}{\mathcal{X}}
\newcommand{\cY}{\mathcal{Y}}
\newcommand{\cZ}{\mathcal{Z}}
\newcommand{\wt}{\widetilde}
\newcommand{\ol}{\overline}
\newcommand{\hT}{\mathbf{T}}
\newcommand{\mB}{\mathcal{B}}
\newcommand{\F}{\mathbb{F}}
\newcommand{\mA}{\mathcal{A}}
\newcommand{\bA}{\mathbb{A}}
\newcommand{\bP}{\mathbb{P}}
\newcommand{\C}{\mathbb{C}}
\newcommand{\mY}{\mathcal{Y}}
\newcommand{\mV}{\mathcal{V}}
\newcommand{\mQ}{\mathcal{Q}}
\newcommand{\mR}{\mathcal{R}}
\newcommand{\mG}{\mathcal{G}}
\newcommand{\mH}{\mathcal{H}}
\newcommand{\mU}{\mathcal{U}}
\newcommand{\mT}{\mathbf{T}}
\newcommand{\hTr}{\texttt{Tr}}

\newcommand{\PGL}{\mathbf{PGL}}
\newcommand{\PGaL}{\mathbf{P\Gamma L}}
\newcommand{\cl}{\mathrm{cl}}
\newcommand{\Top}{\texttt{top}}

\newcommand{\prf}{\textit{Proof. }}
\newcommand{\Sch}{\wis{Sch}}
\newcommand{\bL}{\mathbb{L}}
\usetikzlibrary{matrix}
\usetikzlibrary{arrows}
\usepackage{pgf}
\usepackage{lscape}
\usepackage{listings}

\usepackage{enumitem}

\title{Virtual motives for synthetic geometries, A.\\ Definition and properties of $K_0(\mQ_\ell)$}
\author{Koen Thas}
\address{Ghent University, Department of Mathematics, Krijgslaan 281, S25, B-9000 Ghent, Belgium}
\email{koen.thas@gmail.com}

\begin{document}
\maketitle

\begin{abstract}
In this note, we introduce the first basics on Grothendieck rings for incidence geometries as a new motivic way and tool to study synthetic geometry. In this first instance, we concentrate on generalized quadrangles and related geometries. Many questions, new properties and insights arise along the way. 
\end{abstract}

\bigskip
\setcounter{tocdepth}{1}

\medskip
\begin{tcolorbox}
\tableofcontents
\end{tcolorbox}

\medskip

\section{Introduction}
\label{intro}

\medskip
\subsection{PPC}

We start this section with an infamous conjecture in combinatorics, which has resisted many attacks, but is still standing strong: 

\begin{conjecture}[Prime Power Conjecture (ppc)] 
The order of a finite projective plane is a prime power.
\end{conjecture}

\medskip
To our understanding, there is no clear essential reason {\em why} this conjecture should be true; all the known examples of finite projective planes|which come in many infinite classes|satisfy the conjecture, but to our best of knowledge there has never been discovered a common property for all projective planes of the same finite order that 
points to the direction of the ppc.

The ppc has proved to be an extremely important motivation|sometimes in an implicit guise|to develop a very deep theory underlying axiomatic projective planes.  The first breakthrough in the quest for understanding the order of finite planes was famously published by Bruck and Ryser in 1949 \cite{BrRy}: {\em  if there is a finite projective plane of order $N$, and 
$N \equiv 1\mod{4}$ or $N \equiv 2\mod{4}$, then $N$ is the sum of two squares.}  (The proof is short but ingenious:  the authors express the assumption  in terms of incidence matrices, and then pass to Hasse-Minkowski theory of rational equivalence of quadratic forms.)  Despite this marvelous result, this is the \ul{only known result on the order of a finite plane} whose statement does not impose further geometric or group theoretical restrictions on  the plane. This was 70 years ago!  

\medskip
\subsection{Coordinatization}

One promising general approach to analyzing planes, with possible applications to the prime power conjecture, certainly came with the coordinatization method devised by Hall in 
\cite{Hall}.  Defining an algebraic structure in which the binary operations are governed by incidence relations in a given plane, opened up a startling new vista  to algebraic methods. 
The importance of this approach is very clearly outlined in chapters V--VI--VII--VIII--IX of Hughes and Piper's seminal book \cite{HP}. 
But the prime power conjecture did not come out of the many papers on the algebraic side of coordinatization, although many sub-results were obtained which tried to characterize planes by using the underlying ternary rings. Wedderburn's famous result that states that finite division rings are always finite fields \cite{Wed} (and hence have a prime power number of elements), and later Artin-Zorn's generalization to finite alternative division rings \cite{ArtZor}, are the model results which inspired much of the algebraic theory in the finite case.\footnote{Besides the nonexistence of projective planes of order $N$ which violate Bruck-Ryser's condition \cite{BrRy}, one very famous nonexistence result is the nonexistence of projective planes of order 10, which was obtained in \cite{Lam} by Lam, Thiel and Swiercz, after a 3000 hour computer search on a CRAY-1A super computer (and using theoretical work of many other authors, such as John Thompson's work related to combinatorial projective curves).} 

Still, coordinatization is perhaps the only known method which handles all projective planes ``at the same time.'' The only other results, conjectures and applications known today which give a prime power result, all rely on hypothesized {\em local or global group actions on  projective planes}. Also, the fact that the Andr\'{e}-Bruck-Bose approach to translation planes \cite{And,BrBo} gave|besides a clear prime power result|a representation method for planes in projective spaces (which also works in the infinite case), gave the theory an unexpected boost beyond the initial goals. Note the author's recent breakthrough on the Andr\'{e}-Bruck-Bose theorem for other generalized polygons \cite{TGQkernel}. In Andr\'{e}-Bruck-Bose theory one associates a vector space over a division ring $\ell$ to the said automorphism group, and the characteristic $\ell$ is the prime (including $0$) one is looking for. 

Other hypothesized group actions often led to classical planes $\mathbb{P}^2(k)$, and hence lost any grip on general planes. Still, { under the assumption that a finite projective plane admits a free transitive automorphism group, it is a famous conjecture that the ppc is true}, even when the group is assumed to be cyclic! 

We mention that for generalized quadrangles|the protagonists of the present paper| there is also a well-developed theory of coordinatization introduced by Hanssens and Van Maldeghem \cite{HaHVM1,HaHVM2}, and in Van Maldeghem's book \cite[section 3.2]{POL}, coordinatization for ``general generalized $2n$-gons'' with an order is described.  

Much more about prime power conjectures of point-line incidence geometries can be found in the author's paper \cite{Order}.

\medskip
\subsection{The present paper: ``cohomology without vector spaces''}

One of our main goals of the present paper is to add a new chapter to the algebraic theory of incidence geometries, and in this paper (which is the first instance in a series) we will almost solely concentrate on generalizations of generalized quadrangles. (The projective planes are next.) We will associate a commutative ring to a category $\mQ_\ell$ which is generated by all (known and unkown) generalized quadrangles with $\ell + 1$ points incident with a line. We hope that this formalism introduces a new way to study (in this case) generalized quadrangles from the viewpoint of cohomology (more precisely: virtual motives) without having the appropriate commutative rings at our disposal to build a Zariski topology. 
Our construction is directly inspired by the {\em Grothendieck ring of $k$-varieties} with $k$ a field, and the latter contains cohomological/motivic information of the category of $k$-varieties. Let us first define this ring.  

\subsubsection{Grothendieck rings and Zariski topology}

Let $k$ be a given field. 
The {\em Grothendieck ring of varieties over $k$}, denoted as $K_0(\mV_k)$ (where $\mV_k$ is the category of $k$-varieties), is formally freely generated as an additive abelian group by the isomorphism classes $[X]$ of varieties ${X}$ over $k$ (call this group $\mA(\mV_k)$), and endowed with two types of relations:
\begin{equation}
\Big[A\Big]\ =\  \Big[A\setminus B \Big]\ +\  \Big[B\Big] 
\end{equation}
for any closed subset $B$ of $A$ in the Zariski topology, and $A$ any $k$-variety, and with a product structure given by
\begin{equation}
\Big[U \Big]\cdot\Big[V\Big]\ := \ \Big[U \times_{\Spec(k)} V\Big].
\end{equation}
 
The ring $K_0(\mV_k)$ captures certain aspects of motives, and hence of various (Weil-)cohomological theories associated to $\mV_k$|we refer to section \ref{defK0V} for a short discussion. For being able to define the ring $K_0(\mV_k)$, we need several ingredients which come with $\mV_k$:
\begin{itemize}
\item
the Zariski topology of the elements of $\mV_k$;
\item
a good notion of ``isomorphism'' for the elements in $\mV_k$; 
\item
a product 
\[ \otimes:\ \Big(A, B \Big) \ \mapsto\ A \otimes B \]
which maps couples of $k$-varieties to $k$-varieties.
\end{itemize}

Furthermore, as we will see, we also need to dispose of a good theory of lines. In order to have a Zariski topology at hand, the starting point is the category of commutative rings, and the main ingredients are the notions of maximal ideal, prime ideal and ideal to begin with (for a given commutative ring $A$, say). Once the affine scheme $\Spec(A)$ is defined, 
the maximal ideals correspond to closed points of $\Spec(A)$, the prime ideals|which define the points of the topology|correspond to irreducible subvarieties in $\Spec(A)$, and the general ideals give rise to the ``general subvarieties.'' The Zariski topology is a natural topology which involves all varieties which lie on $\Spec(A)$. In our case, our main geometrical objects are generalizations of generalized quadrangles, and in general, we do not have commutative rings associated to these geometries in a natural way.




\subsubsection{Geometrical approach}

Let $k$ be a finite field $\F_q$.
Irreducible $k$-quadrics in low dimensions ($n = 3, 4, 5$) verify the axioms for generalized quadrangles; quadrics, together with low-dimensional Hermitian $k$-varieties (in which case $k = \F_{q^2}$), represent|up to point-line duality|the so-called finite ``classical examples'' of these geometries. They are all fully embedded in the ambient projective space, have highly symmetrical automorphism groups (of Lie type), and because of that, enjoy highly symmetrical combinatorial properties. Even more so, the points and lines of these classical generalized quadrangles are the $k$-rational points and lines of an associated projective $k$-variety, and so all these examples come with a well-defined Zariski toplogy (and so define classes in the ring  $K_0(\mV_k)$).  (Quadrics are birational to projective spaces, and that is the key fact for their classes to be inside $\Z[\bL]$; they come with virtual mixed Tate motives. We will consider them in some detail in Appendix \ref{quadapp}.)

On the other hand, there is a deeply 
developed theory of {\em synthetic generalized quadrangles}, and many infinite classes of non-classical generalized quadrangles are known; those are not embedded in projective space at all (by a theorem of Buekenhout and Lef\`{e}vre \cite{BueLef}), usually have rather rigid automorphism groups, etc. Still, in these examples, one can see obvious subsets which should be closed sets in any reasonably defined Zariski topology. 
In fact, what is very clear throughout the synthetic theory of generalized quadrangles is the fundamental role those subgeometries  play in characterizations and classification results, which exactly are combinatorial versions of what would be irreducible sub $k$-varieties in the classical quadrangles.

{\em Simple example.}\quad 
In any of the generalized quadrangles which are embedded in a finite projective space, the set of points $a^\perp$ collinear with an arbitrary point $a$ is the intersection of the quadrangle with a hyperplane. So this point set is the set of rational points of some projective variety. It makes much sense to see sets $a^\perp$ in abstract generalized quadrangles in this way, and hence to see them as synthetic irreducible closed subvarieties (prime ideals) in a Zariski topology.

This ``translation method'' is our starting point for founding the first definitions.


Obviously the topology we obtain will lack a lot of information in comparison to the scheme-theoretic picture (where for example a $k$-scheme with $k$ a field  covers the geometry over field extensions of $k$ as well via  Galois descent). In some classes of examples (where field extensions make sense), a much richer Zariski topology can be established, and we describe such a class in detail. In such cases, we do have a Galois descent at our disposal.

\subsubsection{Theory of lines}

In the theory of Grothendieck rings of $k$-varieties the class $[\mathbb{L}]$ of the affine $k$-line plays an exceptional role. Many foundational questions involve this class, and the subring $\mathbb{Z}[\mathbb{L}]$ describes virtual mixed Tate motives. Arguably the most striking identity in $K_0(\mV_k)$ is
\[  [\bP^m(k)] \ =\ \mathbb{L}^m \ + \ \mathbb{L}^{m - 1}\ + \ \ldots\ +\ [\mathrm{point}]                          \]
which arises inductively from the identity 
\[   [\bP^m(k)] \ = \ [\bP^{m - 1}(k)]\ +\  [\mathbb{A}^m(k)].               \]

In our theory, its natural avatar will come from the theory of geometrical hyperplanes in generalized quadrangles and more general geometries. On the other hand, we won't dispose of one line anymore in $\mQ_\ell$ (the category we are using in this paper): as we will point out, due to the lack of natural structure of lines in the incidence-geometrical context (in stark contrast with algebraic geometry), we will refine the concept of line (and other degenerate geometries) in order to make the theory more capable. It would be very interesting to understand the subring of $K_0(\mQ_\ell)$ generated by the isomorphism classes of all lines in $\mQ_\ell$. We refer to \cite{SPL} for a detailed discussion about synthetic projective lines, inspired by this paper.

\medskip
\subsection{Organization of the paper}

In section \ref{defK0V} we give a short introduction to Grothendieck rings of varieties (and motivate some of the underlying mechanisms). In sections \ref{subq} and \ref{span}, we introduce some elementary notions in the theory of generalized quadrangles which will be needed throughout the paper.  We define the version of Zariski topology we will use in section \ref{Zariski}, and for that matter we have a detailed look at how primal geometries should look like in the category $\mQ_\ell$. 
Then, in section \ref{isom} we study isomorphisms of degenerate subgeometries of generalized quadrangles. In the subsequent section \ref{tracegeom}, we introduce and study trace geometries|local geometries which naturally arise in generalized quadrangles. Some interesting questions arise. In section \ref{theline} we investigate some consequences of the concepts of lines and birationality in the theory of translation generalized quadrangles. In the next section, section \ref{prod}, we introduce the notion of product we will use in the present paper (which is based on a product defined on underlying graphs). In section \ref{grogro} we define $K_0(\mQ_\ell)$. In section \ref{Krull} we define a natural notion of dimension in generalized quadrangles (and related geometries), and formulate some natural properties and questions. In section \ref{Lef} we discuss phenomena which arise in $K_0(\mQ_\ell)$ when $\ell$ is not finite, especially in the context of lines, and we complete the discussion started in section \ref{Zariski} on affine geometries. In section \ref{BASE} we formulate some ideas regarding base extension and Galois descent, and concentrate on a class of quadrangles coming from a specified class of algebraic curves. Finally, in the appendix (section \ref{quadapp}), which does not contains original research,  we have a look at decompositions of quadrics in $K_0(\mV_k)$ for the benefit of the incidence-geometrical reader.

\medskip

\section{The Grothendieck ring of varieties over a field}
\label{defK0V}

\begin{quote}{\footnotesize Parmi toutes les choses math\'{e}matiques que j'avais eu le privil\`{e}ge de d\'{e}couvrir et d'amener au jour, cette r\'{e}alit\'{e} des motifs m'appara\^{i}t encore comme la plus fascinante, la plus charg\'{e}e de myst\`{e}re|au coeur m\^{e}me de l'identit\'{e} profonde entre la ``g\'{e}om\'{e}trie" et l'``arithm\'{e}tique.'' Et le ``yoga des motifs'' auquel m'a conduit cette r\'{e}alit\'{e} longtemps ignor\'{e}e est peut-\^{e}tre le plus puissant instrument de d\'{e}couverte que j'aie d\'{e}gag\'{e} dans cette premi\`{e}re p\'{e}riode de ma vie de math\'{e}maticien.
\begin{flushright}A. Grothendieck, \emph{R\'{e}coltes et Semailles}\end{flushright}}
\end{quote}

\subsection{The ring $K_0(\mV_k)$}

Let $k$ be a given field. 
The {\em Grothendieck ring of varieties over $k$}, denoted as $K_0(\mV_k)$ (where $\mV_k$ is the category of $k$-varieties), is formally freely generated as an additive abelian group by the isomorphism classes $[X]$ of varieties ${X}$ over $k$ (call this group $\mA(\mV_k)$), and endowed with two types of relations:
\begin{equation}
\Big[A\Big]\ =\  \Big[A\setminus B \Big]\ +\  \Big[B\Big] \ \ \mbox{{\bf (``scissor relations'')}}
\end{equation}
for any closed subset $B$ of $A$, and $A$ any $k$-variety, and with a product structure given by
\begin{equation}
\Big[U \Big]\cdot\Big[V\Big]\ := \ \Big[U \times_{\Spec(k)} V\Big].
\end{equation}

Note that the set of all elements $[A] - [A \setminus B] - [B]$ (A a $k$-variety, $B$ closed in $A$), forms an ``ideal'' in $\mA(\mV_k)$.\\
 
We denote by $\bL = [\mathbb{A}^1(k)]$ the class of the affine line over $k$. The neutral element for addition is given by $[\emptyset] =: 0$, and the neutral element for multiplication is $[\Spec(k)] =: 1$. \\

The specific form of the scissor relations stems from motivic aspects, as we will indicate below.
 
Suppose $k$ is contained in $\mathbb{C}$, and let $X$ be a $k$-variety. Define
\[ \chi^c(X) := \sum_{k \geq 0}(-1)^kb_k^c(X),  \]
with $b_k^c(X) := \mathrm{dim}_{\mathbb{C}}H^k_c\Big(X(\mathbb{C}),\mathbb{C}\Big)$ (where we consider singular cohomology with compact support). Then for $Y$ closed in $X$, we have that 
\[ \chi^c(X) = \chi^c(X \setminus Y) + \chi^c(Y).   \]

Also, Gillet and Soul\'{e} proved in \cite{GiSo} that there is a homomorphism from $K_0(\mV_k)$ to the $K_0$ of the 
category of ``pure motives'' if $k$ has characteristic $0$. \\

\medskip
\subsection{The motivic formalism}

Grothendieck's theory of \emph{motives}\index{motive} predicts a universal cohomology 
theory $h$ for ``good'' cohomology theories of (say) varieties over fields.

The functor $h$ must satisfy:
\begin{enumerate}
\item[(M$_1$)]
the K\"{u}nneth formula\index{K\"{u}nneth formula}
\begin{equation}
h(V \times W) = h(V) \otimes h(W);
\end{equation}
\item[(M$_2$)]
the property that disjoint unions are translated into direct sums;
\item[(M$_3$)]
certain additional axioms to obtain the Lefschetz formula.\\
\end{enumerate}

Let $k$ be a field.
The universal cohomology theory $h$ should take values in a category of motives $\wis{Mot}(k)$\index{$\wis{Mot}(k)$} which should look like the category of finite-dimensional $\mathbb{Q}$-vector spaces|more precisely:
\begin{enumerate}
\item[(M$_1'$)]
homomorphism groups should ``be'' $\mathbb{Q}$-vector spaces;
\item[(M$_2'$)]
$\wis{Mot}(k)$ should be an abelian category (or even better a ``tannakian'' category, cf.~\cite{DelMil}).\\
\end{enumerate}

Every Weil (= ``good'') cohomology theory $\mathbf{H}$ with coefficients in some field $\wp$ (such as the \'{e}tale $\ell$-adic one) should fit (uniquely) into the following diagram
\begin{equation}
\mathbf{H}\colon \widetilde{\mV_k} \overset{h}{\longrightarrow} \wis{Mot}(k) \xrightarrow{\omega_{\mathbf{H}}\,} \{ \text{graded }\wp\text{-vector spaces}\},
\end{equation}
where we denote the category of nonsingular projective varieties over $k$ by $\widetilde{\mV_k}$\index{$\wis{Var}(\F_q)$}.
Here, $\omega_{\mathbf{H}}$\index{$\omega_{\mathbf{H}}$} is a functor which comes with the cohomology theory $\mathbf{H}$, such that
\begin{equation}
\omega_{\mathbf{H}}\bigl(h(X)\bigr) = \mathbf{H}^*(X) = \bigoplus_{i = 0}^{2\mathrm{dim}(X)}H^i(X).
\end{equation}

The functor $\omega_{\mathbf{H}}$ gives rise to decompositions of the $h(X)$ in ``pure pieces'' which yield a lot of extra information about $X$. In case of a projective space $\bP^n(k)$ over any base field $k$, it can be proven that 
\[  h(\bP^n(k)) \ = \ \bL^n \ +\ \bL^{n - 1}\ + \ \cdots\ +\ 1,       \]
where $\bL := h(\bA^1(k))$ is the so-called {\em Lefschetz motive}. 
This is the star example of what Serre calls ``m{e}ccano des motifs'' \cite{Serre}. \\

If we now turn back to the ring $K_0(\mV_k)$ ($k$ a field) for a moment, we see that  
due to the scissor relations, we also have identities such as
\[  \Big[\mathbb{P}^n(k) \Big] \ = \Big[\bP^{n - 1}(k)\Big] \ + \ \Big[ \bA^n(k) \Big],        \]
given us the typical motivic composition of the kind discussed above: 

\begin{align} 
\Big[\bP^n(k) \Big] \ &=\ \Big[\bA^{n}(k)\Big] \  + \Big[\bA^{n - 1}(k)\Big]\ +\ \cdots\ +\ 1 \nonumber \\
&= \ \bL^n \ +\ \bL^{n - 1}\ + \ \cdots\ + \ 1.
\end{align}

Many identities of this type in $K_0(\mV_k)$ correspond to analogous motivic decompositions in $\wis{Mot}(k)$. It is one of its many charms.

\subsubsection{Additive invariants}

Let $R$ be a commutative ring.
An \emph{additive invariant}\index{additive invariant} $\chi\colon {\mV_k} \to R$ is a map with the following properties (where $X, Y$ are objects in $\mV_k$):\\

\medskip
$\left\lgroup
\begin{tabular}{p{0.9\textwidth}}
\begin{description}
\item[\rm\textsc{Isomorphism Invariance:}]
$\chi(X) = \chi(Y)$ if $X \cong Y$;
\medskip
\item[\rm\textsc{Multiplicativity:}]$\chi(X \times Y) = \chi(X)\chi(Y)$;
\medskip
\item[\rm\textsc{Inclusion-Exclusion:}]$\chi(X) = \chi(Y) + \chi(X \setminus Y)$ for $Y$ closed in $X$.
\end{description}
\end{tabular}
\right.$
\medskip

A standard example is the \emph{topological Euler characteristic}\index{Euler characteristic}. It is clear that giving an additive invariant $\chi$ is the same as giving a ring morphism
\begin{equation}
\chi\colon K_0(\mV_k) \longrightarrow R.
\end{equation}

Let $\mathbb{L}$\index{$\mathbb{L}$} be the class $[\A^1(k)]$ in $K_0(\mV_k)$ | 
it is called the virtual \emph{Lefschetz motive}\index{virtual!Lefschetz motive}.
The subring $\mathbb{Z}[\mathbb{L}] \subset K_0(\mV_k)$ is the subring of virtual \emph{mixed Tate motives}\index{virtual!mixed Tate motive}. Let us say that a $k$-variety $X$ has a \ul{mixed Tate motive} if $[X] \in \mathbb{Z}[\mathbb{L}]$.

\subsection{Mixed Tate motives and counting polynomials}

If $X$ is a $\Z$-variety, and $X$ corresponds to a mixed Tate motive (over $k$), it has class
\begin{equation}
[X \times_{\Spec(\Z)} \Spec(k)] \in \Z[\mathbb{L}]\subset K_0(\mV_k).
\end{equation}
Up to a finite number of primes (where bad reduction phenomena could occur), it follows that 
\begin{equation}
N_{p^m}(X) := \#\bigl(X \times_{\Spec(\Z)}\Spec(\F_{p^m})\bigr)
\end{equation}
(with $p$ a prime and $m$ a nonzero integer) is a polynomial in $p^m$, since the counting function $N_{p^m}(\cdot)$ is also an additive
invariant and since $N_{p^m}\bigl(\A_1(\F_{p^m})\bigr) = p^m$. We refer to \cite{Andre} for more details.

\medskip
\subsection{Naive motives versus ``real'' motives}

Although $K_0(\mV_k)$ is often considered as a naive model for working with ``real'' motives, it is still the case that $K_0(\mV_k)$ can see certain properties which real motives cannot. A very interesting discussion can be found in the thread \cite{MOvirtual}. Under the cut-and-paste conjecture, two $k$-varieties are equivalent in the Grothendieck ring if and only if they can both be decomposed into the same set of locally closed pieces. There are actually many pairs of $k$-varieties with the same motive but which can not be cut-and-pasted into each other. {\em Example:} $\bP^2(k)$ and a fake projective plane over $k$. So the Grothendieck ring {\em does} remember quite some information about the varieties. In particular, after work of Larsen and Lunts \cite{LarLun}, one can show that two smooth projective $\C$-varieties with the same class in the Grothendieck ring $K_0(\mV_\C)$ have isomorphic fundamental groups, which is a highly nonabelian invariant! So the category of motives loses some information which is contained in the Grothendieck ring | in some sense the ``nonabelian information'' \cite{MOvirtual}. \\

Poonen was the first to show (in \cite{Poonen}) that if $k$ is a field of characteristic $0$, then $K_0(\mV_k)$ is not a domain (by finding elliptic curves $E$ and $E'$ over $k$ for which $[E] \ne [E']$ but nevertheless $([E] + [E'])([E] - [E']) = 0$). Many other related results have been found since \cite{Poonen}|see e.g. \cite{Kollar}. \\

On the other hand, $K_0({\mV_k})$ does not remember nontrivial extensions, but Voevodsky's category of motives does. For example, let $\mC$ be a smooth projective curve of positive genus, and let $x$ and $y$ be two rational points of $\mC$. Then in Voevodsky's version, $h^1(\mC \setminus \{ x, y\})$ depends nontrivially on the class of $[x] - [y]$ in the Jacobian of $\mC$. But the expression $[\mC] - [x] - [y]$ in $K_0({\mV_k})$ does not depend on the choice of $x$ and $y$. 
In our approach to ``synthetic Grothendieck rings'' the latter property is in general \ul{not} true (and such classes will depend on the choice).

\medskip
\subsection{An open problem}

One problem we are interested in, is the following: 
\medskip
\begin{tcolorbox}
\begin{question}
Suppose given a commutative unital ring $\Big(A,+,\cdot\Big)$ with the additional information that $A \cong K_0(\mV_k)$ for some unknown field $k$.  
\begin{itemize}
\item[{\rm (A)}]
Can $k$, or at least the characteristic of $k$, be read from $A$?  
\item[{\rm (B)}]
{\em Variation}: if $k$ and $\ell$ are fields of different characteristic, then can  $K_0(\mV_k) \cong K_0(\mV_\ell)$?  
\end{itemize}
\end{question}
\end{tcolorbox}

This set of problems appears to be wide open, and difficult to access. One idea to single out the characteristic of $K_0(\mV_k)$ would be to look at the subring $K_0(\mV_k^0)$ generated by the $0$-dimensional $k$-varieties \cite{Deligne_lett}. Note first that for any field $k$, a $0$-dimensional scheme is just a finite union of points. 
In what follows, we suppose that $k$ is a finite field. Then $K_0(\mV_k^0)$ can be identified with the free $\Z$-module with basis the isomorphism classes of finite transitive permutation representations of $\Gal(\overline{k}/k)$ \cite{Deligne_lett}. \\

$\Big($To any finite field extension $k \mapsto \ell$, one associates the set $\mathrm{hom}_k(\ell,\overline{k})$, which admits a natural transitive Galois action. And conversely, as a finite transitive $\Gal(\overline{k}/k)$-action can be identified with $(G,G/H)$ for some subgroup $H \leq \Gal(\overline{k}/k) = G$ of finite index, just define $\ell$ as the fixed subfield in $\overline{k}$ of $H$.$\Big)$\\

As 
\[ \widehat{\Z}\ \cong\ \Gal(\overline{k}/k)   \]
is independent of the choice of finite field $k$ (refer to section \ref{absgal} for more details on $\widehat{\Z}$), we cannot distinguish one finite field from another through this approach \cite{Deligne_lett}. 

Especially when we will develop the $K_0(\cdot)$-theory of synthetic projective planes (in the next part of the series), a positive answer to an appropriate variation of question A 
would be highly desirable!

\medskip
\section{Quadrangles and subquadrangles}
\label{subq}

A {\em generalized quadrangle} of {\em order $(s,t)$} is a point-line incidence structure $\Gamma = (\mP,\mL,\I)$, where $\mP$ and $\mL$ are respectively the point and line set, and 
$\I \subseteq (\mP \times \mL) \cup (\mL \times \mP)$ is a symmetric incidence relation which satisfies the following conditions. 

\begin{itemize}
\item[(1)]
Each line is incident with $s + 1$ points.
\item[(2)]
Each point is incident with $t + 1$ lines.
\item[(3)]
If $(x,Y)$ is a non-incident point-line pair, there is precisely one point $z$ such that $z \I Y$ and $z$ and $x$ are collinear. 
\end{itemize}

If $s \geq 2$ and $t \geq 2$, we say that $\Gamma$ is {\em thick}; otherwise it is {\em thin}. In this paper, we allow that $s = 0$ or $t = 0$, and 
even that $(s,t) = (0,0)$; in the latter case, the quadrangle consists of an incident point-line pair. If $\Gamma$ has order $(s,0)$, we 
have one line with $s + 1$ incident points, and if the order is $(0,t)$, we have one point with $t + 1$ incident lines.

\begin{center}
\item
\begin{tikzpicture}[scale=1]
\draw [line width=0.45mm, black ] (-4,0) -- (0,0) node[below] {$Y$};

\draw [line width=0.45mm, black ] (-2,2) -- (-2,0) node[below] {$z$};

\draw[fill=blue] (-2,2) circle (.5ex) node[right] {$x$};
\draw[fill=blue] (-2,0) circle (.5ex) node[above] {};


\end{tikzpicture}
\item
\item
Axiom (3).
\end{center}

Generalized quadrangles, or ``GQs'' for short, were introduced in 1959 by Tits in a seminal paper \cite{Tits1959}, as natural geometric modules on which semisimple Lie groups of relative rank $2$ act. For that reason, but also for their uncanny ability to be used and  applied in Incidence Geometry and far beyond, for the deep construction, classification and characterization theory of classical and  nonclassical examples (think for instance of Tits and Weis' famous classification of Moufang polygons \cite{TiWe}, or the monograph Payne and Thas \cite{PTsec}), they are considered to be the central incidence geometries of rank $2$. 

\subsection{Duality}

If we replace the notions ``point" and ``line" in a point-line geometry which satisfies (3), then we obtain again a point-line geometry 
which satisfies (3).

\subsection{Subquadrangles}

Let $\Gamma = (\mP,\mL,\I)$ be a generalized quadrangle. Then the point-line incidence geometry $\Gamma' = (\mP',\mL',\I')$ 
is called a {\em subquadrangle} (subGQ) of $\Gamma$ if $\mP' \subseteq \mP$, $\mL' \subseteq \mL$, and if 
moreover $\I'$ is the incidence relation inherited from $\I$. If $\Gamma'$ is a subquadrangle of $\Gamma$, we will sometimes express this by writing $\Gamma' \leq \Gamma$. \\

{\em Examples}.\quad
In section \ref{finex} below (we refer also to that section for the notation used hereafter) we will see that 
\[ \mQ(3,q) \ \leq\ \mQ(4,q)\ \leq\ \mQ(5,q)       \]
for each finite field $\F_q$. The action of the automorphism groups of $\mQ(4,q)$ resp. $\mQ(5,q)$ yields the fact that $\mQ(3,q)$-quadrangles resp. $\mQ(4,q)$-quadrangles live in $\mQ(4,q)$ resp. $\mQ(5,q)$ in abundance. \\

We refer to \cite[chapter 1]{PTsec} for more theory on subquadrangles. For now, we will content ourself with its mere definition.

\subsection{Geometries satisfying axiom (3)}

We also consider other thin generalized quadrangles than the ones already mentioned: 
those which do not have an order, but that still satisfy axiom (3).

They are of the following types:

\begin{itemize}
\item
the empty geometry $\varnothing$.
\item
a set of points but no lines.
\item
a set of lines but no points.
\item
$u \times v$-grids: two sets $\mR_1$ and $\mR_2$ of mutually nonconcurrent lines, called {\em reguli}, such that each line of one regulus is 
concurrent with each line of the other, with $\vert \mR_1 \vert = u$ and $\vert \mR_2 \vert = v$.
\item
dual $u \times v$-grids: the duals of the previous examples.
\item
perp geometries: one point $x$ incident with $t + 1$ lines, each line being incident with an arbitrary number ($\ne 0$) of points (including $x$). We again allow $t = 0$ (a case which was already mentioned).
\item
dual perp geometries.
\end{itemize}

\begin{proposition}
If a point-line incidence structure $\mS = (\mP,\mL,\I)$ satisfies axiom (3), then it is isomorphic to one of the types described in this section.
\end{proposition}

{\em Proof}.\quad
First note that if $\mS$ does not contain a generalized quadrangle $(1,1)$ as subquadrangle, it is either isomorphic to $\varnothing$, a point set, a line set, a perp geometry, or a dual perp geometry. If it does contain a subquadrangle of order $(1,1)$, and it is not thick, show that it is either a grid or a dual grid. \eop \\

In this paper we only consider perp geometries with a constant number of points incident with the lines, and we suppose the dual property for dual 
perp geometries.

\subsection{Axiom $\widetilde{(3)}$}

A point-line geometry $\Gamma$ is said to satisfy {\em axiom $\widetilde{(3)}$} if it contains no triangles as ordinary subgeometries. Equivalently, we have:
\begin{itemize}
\item[$\widetilde{(3)}$]
If $(x,Y)$ is a non-incident point-line pair, there is at most one point $z$ such that $z \I Y$ and $z$ and $x$ are collinear. 
\end{itemize}

\medskip
\section{Symmetry, SPGQ theory and TGQs}
\label{span}

The theory of {\em span-symmetric generalized quadrangles} and {\em translation generalized quadrangles} will be used in some of our guiding examples in the set-up 
of the isomorphism theory we will formulate for local geometries. In this section, we will give a short overview with detailed references to 
the most useful properties and results.

\subsection{Symmetry}

Let $\Gamma$ be a thick generalized quadrangle of order $(s,t)$. We define {\em automorphisms} of $\Gamma$ in the usual 
manner, and denote by $\Aut(\Gamma)$ its automorphism group. 

Let $U$ be a line in $\Gamma$. A {\em symmetry} with axis $U$ is an automorphism which fixes each line in $U^{\perp}$. Let $V \sim U$ be 
any line in $U^{\perp} \setminus \{U\}$. If the group $S(U)$ of all symmetries with axis $U$ acts transitively on the points of $V \setminus  \{ U \cap V \}$, we say that $U$ is an {\em axis of symmetry}. The following two properties can then be shown: 
\begin{itemize}
\item[(1)]
the definition is independent of the choice of $V$;
\item[(2)]
the action of $S(U)$ on $V \setminus \{ U \cap V\}$ is sharply transitively. 
\end{itemize}

It is also straightforward to prove that $s \leq t$ if $\Gamma$ has an axis of symmetry.

\subsection{Grids and regularity}

Abstractly, we will sometimes call a thin quadrangle with parameters $(m,1)$ a {\em grid}. We say that a pair $(U,V)$ of nonconcurrent lines in a thick generalized quadrangle $\Gamma$ of order $(s,t)$ is {\em regular} if they are contained in a (necessarily unique) grid $\Gamma'$ of order $(s,1)$. (If the parameter $s$ would be infinite, we furthermore demand that the lines of the grid are ``full,'' meaning that each point in $\Gamma$ which is incident with $U$ and/or $V$ is also a point of $\Gamma'$.) 
If $(U,V)$ is regular for all noncurrent line-pairs $(U,V)$ with $U$ fixed, then $U$ is by definition a {\em regular line}. It is easy to show that an axis of symmetry is a regular line.

\subsection{SPGQs}

If a thick generalized quadrangle $\Gamma$ has nonconcurrent axes of symmetry $A$ and $B$, we call it {\em span-symmetric}, and 
one can show that each line in $\{ A, B\}^{\perp\perp}$ is an axis of symmetry. The set $\{A,B\}^{\perp\perp}$ is called the {\em base-span}. 
A span-symmetric generalized quadrangle is also abbreviated as ``SPGQ.''

If $\Gamma$ is thick, finite and with parameters $(s,s)$, then it has been shown independently by Kantor \cite{kantor} and the author \cite{spgq}, that 
$s$ is a prime power, and $\Gamma \cong \mQ(4,s)$ (cf. section \ref{orth} below for the definition of the latter quadrangle). Interestingly enough, this implies that $\Aut(\Gamma)$ acts transitively on all lines, starting from a local group-theoretical condition. \\

The theory of SPGQs of order $(s,t)$ with $s \ne t$ is more challenging, and we refer to the papers \cite{two,tgqduals} and the monograph \cite{SFGQ} for the details. Let us just mention | and only for finite quadrangles | that:
\begin{itemize}
\item[(a)]
one can prove that $t = s^2$ with $s$ a prime power;
\item[(b)]
(a) follows from the fact that an SPGQ of order $(s,t)$ with $s \ne t$ has $s + 1$ different subGQs, all isomorphic to $\mQ(4,s)$ (cf. section \ref{orth} below) and mutually intersecting in the grid with parameters $(s,1)$ that contains the base-span;
\item[(c)]
contrary to the case $(s,s)$, there are nonclassical examples of SPGQs, even with many different base-spans. We will meet one such class of examples in section \ref{KanKnu}. 
\end{itemize}

\subsection{TGQs}

Let $\Gamma$ be a thick generalized quadrangle. If every line incident with some point $u$ is an axis of symmetry, then 
the automorphism group $T$ generated by the symmetries about all axes of symmetry on $u$, is abelian and acts sharply transitively on the points opposite $u$ \cite[\S 8.3]{PTsec}. Vice versa, if $\Gamma$ would have an abelian automorphism group which fixes all lines incident with some point $u$ and which happens to act sharply transitively on the points opposite $u$, then each line on $u$ is an axis of symmetry. (This is an easy exercise). 

This invariantly carries over to thin quadrangles of order $(s,1)$, and also to thick infinite quadrangles. \\

If $\Gamma$, $u$ and $T$ are as above, we call $(\Gamma,u,T)$ (or $\Gamma$, if $u$ and $T$ are clearly defined), a {\em translation generalized quadrangle} (TGQ) with {\em translation point $u$} and {\em translation group $T$}. There is a pretty intricate literature on TGQs, and we welcome the reader to the monograph \cite{TGQ} for much more.

\subsection{Some finite examples}
\label{finex}

\subsubsection{Orthogonal quadrangles}
\label{orth}

Let $\F_q$ be a finite field. Below, we work with homogeneous coordinates $(x_0 : x_1 : \cdots : x_m)$ in a projective 
space $\mathbb{P}^m(q)$. \\

Let $m = 3$. The $\F_q$-rational points and lines lying on a (hyperbolic) quadric with defining equation
\[ X_0X_1\ +\ X_2X_3\ =\ 0  \] 
define a thin generalized quadrangle $\mQ(3,q)$ with parameters $(q,1)$. \\

Let $m = 4$. The $\F_q$-rational points and lines lying on a (parabolic) quadric with defining equation
\[ X_0X_1\ +\ X_2X_3\  + X_4^2\ =\ 0  \] 
define a thick generalized quadrangle $\mQ(4,q)$ with parameters $(q,q)$. \\

Let $m = 5$. The $\F_q$-rational points and lines lying on a(n) (elliptic) quadric with defining equation
\[ X_0X_1\ +\ X_2X_3\  +\ f(X_4,X_5)\ =\ 0,  \]
where $f(X_4,X_5)$ is an irreducible quadratic form in the parameters $X_4, X_5$,  
define a thick generalized quadrangle $\mQ(5,q)$ with parameters $(q,q^2)$. \\

Each generalized quadrangle isomorphic to one of the quadrangles defined above is called {\em orthogonal}, and they share the following properties:
\begin{itemize}
\item[(a)]
each line is an axis of symmetry, so each line is regular;
\item[(b)]
by (a), each point is a translation point, so a finite orthogonal quadrangle is a TGQ relative to any given point.  
\end{itemize}

\subsubsection{Kantor-Knuth quadrangles}
\label{KanKnu}

Although we will not give their precise definition here | for that, we refer to \cite[\S 4.5 and \S 5.7]{TGQ} | the {\em Kantor-Knuth quadrangles} will play some role in the construction of examples in this paper. We will only consider finite Kantor-Knuth quadrangles. Let $q$ be any odd prime power. Let $n \in \F_q$ be nonsquare, and let $\gamma$ be a field automorphism of $\F_q$. Then a generalized quadrangle $\mK(q,n,\gamma)$ of order $(q,q^2)$ can be constructed. 

We have the following properties/notions:
\begin{itemize}
\item[(a)]
$\mK(q,n,\gamma) = \mK(q',n',\gamma')$ if and only if $(q,\gamma) = (q',\gamma')$ (and so the construction is independent of the choice of the nonsquare $n$); 
\item[(a)$'$]
we call $(q,\gamma)$ the {\em type} of $\mK(q,n,\gamma)$; 
\item[(b)]
$\mK(q,n,\gamma) \cong \mQ(5,q)$ if and only if $\gamma = \id$; 
\item[(c)]
If $\gamma \ne \id$, $\mK(q,n,\gamma)$ has a line $U$ for which each element of $U^{\perp}$ is an axis of symmetry, 
and $\Aut(\mK(q,n,\gamma))$ fixes this line. (So there are no other axes of symmetry outside $U^{\perp}$.) Each line on $U$ hence is a translation point.
\item[(d)]
If $\gamma \ne \id$, $\mK(q,n,\gamma)$ has two orbits $\Omega_1$ and $\Omega_2$ of full $\mQ(4,q)$-subGQs which all contain the special line $U$. We use the notation $\Omega_1$ for the orbit of the smallest size. 
\end{itemize} 

Each generalized quadrangle isomorphic to a $\mK(q,n,\gamma)$ is also called {\em Kantor-Knuth quadrangle}.

\medskip
\section{The Zariski topology}
\label{Zariski}

In order to define the Grothendieck ring, we need to assign a Zariski topology to each generalized quadrangle (since they will be representants of 
the generators). The obvious obstruction is that we do not have defining equations, nor coordinate rings at our disposal, so we want to circumvent this lack of underlying 
algebra. Comparing to the classical Zariski topology of a variety over an algebraically closed field $k$, there would be a straightforward way to define the points of the topology: 
in the classical Zariski topology, they are the $k$-rational points, so if $\Gamma$ is a generalized quadrangle, its points could be defined as the points of the topology. If $k$ is 
a finite field, the classical Zariski topology is a discrete topology, since all rational points are closed points, and since finite sets of closed points are closed. The same thing happens if $\Gamma$ is finite: all points of $\Gamma$ should be compared to rational, so closed points, and we obtain a discrete topology. In the infinite case, the situation is different (as the rest of this section will show). 

We will introduce a different, richer, topology, which we model on the modern notion of Zariski topology. In the latter, points correspond to prime ideals, and in their turn those correspond to irreducible subvarieties. So the first thing we do is define irreducible subgeometries of generalized quadrangles. We do not see how to introduce analogons to closed points which are not rational, although we will indicate some situations in which one perfectly can: in some sense the topology will behave {\em as if the quadrangle will be defined over an algebraically closed field.} Also, as we will see, for classical quadrangles such as $\mQ(4,q)$, which consist of the rational points and lines of 
{\em varieties}, the synthetically defined Zariski topology will be coarser than the real Zariski topology. 
 
After that, the closed sets will be introduced. 

Abstractly | that is, in case of generalized quadrangles which are not defined through classes | the finite case has much less structure than the Zariski topology of projective varieties  over finite fields: the finite set of rational points of such a variety does not tell much about the variety, and nor does it about isomorphisms between varieties. Therefore, we need to be cautious when introducing isomorphisms, and on that level we will greatly deviate from the viewpoint of algebraic geometry.

\subsection{Closed subsets (1)}
\label{closed1}

In this section we list a number of types of subsets which we will regard as (corresponding to) closed sets. So they will form a base, or at least
part of a base, for the closed set ``abstract Zariski topology.'' We first list the base elements which should be trivially included. The first candidates are thin, and need closer attention on the level of isomorphisms.

\subsection*{``Point and line sets''}

Here, we consider finite sets. We need to decide what the role of collinearities and concurrencies is in this context. 

\subsection*{``Perps''}

We consider full lines in ``perp geometries,'' but will allow non-ideal perps.

\subsection*{``Subquadrangles''}

An obvious candidate. Here, we will again only consider full examples.

\subsection*{``Geometrical hyperplanes''}

A very important class of candidates, which contains some of the extremal examples of previous items. (We refer to section \ref{affine} for a precise definition.) 
In $K_0(\mV_k)$, we have the typical motivic identity 
\begin{equation}
\Big[\bP^n(k) \Big]\ = \ \Big[\bA^n(k)\Big]\ + \ \Big[\bP^{n - 1}(k)\Big]
\end{equation}
at our disposal. So if $\mQ$ would be a geometrical hyperplane of a generalized quadrangle $\Gamma$, then a similar identity would be 
\begin{equation}
[\Gamma]\ = \ [\bA]\ + \ [\mQ],
\end{equation}
where $\bA$ is the affine quadrangle which arises by deleting $\mQ$. So we definitely want to see each geometric hyperplane of a generalized quadrangle as a closed set. All such identities will remain true in $K_0(\mV_k)$. (See also Appendix \ref{quadapp}.)

\subsection{Points}
\label{points}

We first obtain a result about decompositions of generalized quadrangles with order. The theorem applies to thick quadrangles and thin quadrangles with order, both 
in the finite and infinite case. 
We remark that a highly relevant paper is Van Maldeghem's \cite{hvmunion}, in which he studies finite generalized quadrangles which are unions of a number of full subquadrangles. 

We start with a lemma.

\begin{lemma}
\label{lemsub}
Suppose $\Gamma$ is a generalized quadrangle of order $(s,t)$ with $s \geq 2$, and let $\mS$ be a subgeometry of $\Gamma$.
Let axiom (3) be satisfied. If $\mS$ has one full line $U$, and at least one other line $V$ not meeting $U$, then either $\mS$ is a full thick subquadrangle, 
or $\mS$ is a grid with at least one full regulus.
\end{lemma}

{\em Proof}.\quad
By axiom (3), $V$ is also a full line, and also by (3), all lines in $\{ U,V\}^{\perp}$ are lines of $\mS$. 
Suppose $W$ is a line which is not contained in $\{ U, V\}^{\perp} \cup \{ U, V\}^{\perp\perp}$. Then $W$ cannot be concurrent with 
both $U$ and $V$. It follows easily that $\mS$ is a full thick subGQ of $\Gamma$. \eop \\

\begin{theorem}[Decomposition Theorem]
\label{decthm}
Let $\Gamma$ be a GQ of order $(s,t)$. Then $\Gamma$ cannot be the union of two different proper subquadrangles (where the definition of subquadrangle is taken as in
\S \ref{subq}), unless we are in the following cases.\\

$\left\lgroup
\begin{tabular}{p{0.9\textwidth}}
\begin{itemize}
\item[{\rm (i)}]
$\Gamma$ is a GQ of order $(s,1)$, and $\Gamma_1$ and $\Gamma_2$ are $(s + 1) \times r_i$-grids ($i = 1,2$) which share one regulus. 
\item[{\rm (ii)}]
$\Gamma$ is a GQ of order $(s,1)$, and $\Gamma_2$ is a GQ of order $(s,0)$. 
Then $\Gamma_1$ is an $(s + 1) \times s$-grid. 
\item[{\rm (iii)}]
$\Gamma$ is a GQ of order $(s,1)$, and $\Gamma_2$ is the perp geometry of one point. 
Then $\Gamma_1$ is an $s \times s$-grid. 
\item[{\rm (iv)}]
$\Gamma$ is a GQ of order $(1,t)$, and $\Gamma_1$ and $\Gamma_2$ are dual $r_i \times (t + 1)$-grids ($i = 1,2$) which share one dual regulus. 
\item[{\rm (v)}]
$\Gamma$ is a GQ of order $(1,t)$, and $\Gamma_2$ is a GQ of order $(0,t)$. Then $\Gamma_1$ is a dual $(t + 1) \times t$-grid.
\item[{\rm (vi)}]
$\Gamma$ is a GQ of order $(1,t)$, and $\Gamma_2$ is a dual grid with parameters $(1,t - 1)$, and $\Gamma_1$ is the perp geometry of one line. 
\item[{\rm (vii)}]
$s = 2$ and $t = 2$. In this case, $\Gamma \cong \mQ(4,2)$, and $\Gamma_2$ is a grid with parameters $(2,1)$. A simple combinatorial exercise yields that 
$\Gamma_1$ is a dual grid with parameters $(1,2)$. 
\end{itemize}
\end{tabular}
\right.$
\end{theorem}

{\em Proof}.\quad
Let $\Gamma$ be a generalized quadrangle of order $(s,t)$, and suppose by way of contradiction that $\Gamma = \Gamma_1 \cup \Gamma_2$. By this, we mean that each point and each line of $\Gamma$ is either in $\Gamma_1$ or $\Gamma_2$. \\

First note the following. If $\Gamma_2$ would not have lines, then all lines of $\Gamma$ are lines of $\Gamma_1$. Then any point of $\Gamma_2$ is also a point of $\Gamma_1$, since it is incident with two intersecting lines of $\Gamma_1$. So $\Gamma_2 \leq \Gamma_1$, contradiction. So both $\Gamma_1$ and $\Gamma_2$ contain 
points and lines; in particular, neither one can be an ovoid, or merely a set of points. \\

Below, we will distinguish between the case where at least one of $\Gamma_1, \Gamma_2$ has only full lines (or only ideal points), and the case where none of 
$\Gamma_1, \Gamma_2$ have. From the viewpoint of this paper, this is a natural distinction. \\

\ul{First let $\Gamma_2$ have only full lines.} {\em For now, we suppose that $\Gamma$ is thick}.

Suppose first that $\Gamma_2$ does not have an order. Then since $\Gamma_2$ only has full lines, it must be the geometry of a partial perp. Say it is a subgeometry
of $x^{\perp}$ with $x$ a point of $\Gamma_2$, and note that we allow $\Gamma_2$ to have order $(s,0)$. Then for any $y \ne x$, $y$ a point of $\Gamma_2$, we have that $y$ is incident with at least two lines of $\Gamma_1$. So $y$ is a point of $\Gamma_1$. Whence each line of $\Gamma_2$ is also a line of $\Gamma_1$, and it easily follows that $\Gamma_2 \leq \Gamma_1$, contradiction. \\

Next, let $\Gamma_2$ have an order; then each point of $\Gamma_2$ is incident with at least one line of $\Gamma_1 \setminus \Gamma_2$. 
If each such point is incident with at least two such lines, then each point of $\Gamma_2$ is also a point of $\Gamma_1$, so that $\Gamma_2 \leq \Gamma_1$, contradiction. 
So suppose now each point of $\Gamma_2$ is incident with precisely one line of $\Gamma_1 \setminus \Gamma_2$. Since $\Gamma_1$ does not contain $\Gamma_2$, 
there is at least one line of $\Gamma_1 \setminus \Gamma_2$ which meets $\Gamma_2$ in one point which is not a point of $\Gamma_1$. Now let $y$ be a point in $\Gamma_1 \setminus \Gamma_2$.  Project $y$ onto all lines of $\Gamma_2$; since each line of $\Gamma_2$ is full, we obtain a set of points in $\Gamma_2$, which we denote by $\mO_y$, and 
which meets line of $\Gamma_2$ in precisely one point. Now let $x$ be a point in $\Gamma_2$ which is not in $\mO_y$; Then $\{ x,y\}^{\perp}$ is, up to one point (which is the one point of $\{ x,y \}^{\perp}$ which is not a point of $\Gamma_2$), contained in $\mO_y$. So $x^{\perp}$ contains an ovoid of $\Gamma_2$ minus one point.  

Now suppose $\Gamma_2$ is a thick subGQ. In that case, let $y'$ be the point of $\mO_y$ which is not in $x^{\perp}$. Let $U$ be a line of $\Gamma_2$ incident with $x$. As each line of $\Gamma_2$ must meet $\mO_y$, it follows that if we project all the points of $\mO_y$ on $U$, $U$ must be covered. It follows easily that $s \leq 2$. 
Let $\Gamma_2$ be thin; then $\Gamma_2$ is a either one full line (together with the incident points), or it is a dual grid with parameters $(1,t - 1)$.  

By inspection of the possibilities (and checking in each case whether an ovoid can occur of the above form), the following cases can occur:
\begin{itemize}
\item[(a)]
$s = 1$ and $t = 1$. In this case, $\Gamma_2$ is a line with two points, and $\Gamma_1$ is a root.
\item[(b)]
$s = 1$ and $t$ is arbitrary. In that case, $\Gamma_2$ is a dual grid with parameters $(1,t - 1)$, and $\Gamma_1$ is a panel. 
\item[(c)]
$s = 2$ and $t = 2$. In this case, $\Gamma \cong \mQ(4,2)$, and $\Gamma_2$ is a grid with parameters $(2,1)$. A simple combinatorial exercise yields that 
$\Gamma_1$ is a dual grid with parameters $(1,2)$. 
\end{itemize}

In all these cases, $\Gamma_1$ and $\Gamma_2$ are disjoint. \\

{\em Now suppose that $\Gamma$ is thin}. We handle the case where its order is $(s,1)$.  
Let $U$ be a line of $\Gamma_2 \setminus \Gamma_1$. Then by assumption $U$ is full. If there is no other line
in $\Gamma_2 \setminus \Gamma_1$, then $\Gamma_2$ is a GQ of order $(s,0)$. $\Gamma_1 \ne \Gamma_2$ consists of the other lines and points of $\Gamma$: we 
obtain an $(s + 1)\times s$-grid. 
 
 Let $\mR_1$ be the regulus containing $U$.
 If $\Gamma_2 \setminus \Gamma_1$ contains other lines besides $U$, then $\Gamma_2$ contains all lines of the other regulus $\mR_2$. It easily follows that $\Gamma_2$ is an $(s + 1)\times r$-grid for some $r \geq 2$. Since $\Gamma_1$ contains lines $\widetilde{U} \in \mR_1$ which are not in $\Gamma_2$ (as otherwise $\Gamma = \Gamma_2$), we see that 
 $\widetilde{U}$ is full. (If $u \I \widetilde{U}$ would not be in $\Gamma_1$, axiom (3) would force $\widetilde{U}$ to be in $\Gamma_2$.) It follows that 
 $\Gamma_1$ either is a GQ of order $(s,0)$, or an $(s + 1) \times \widetilde{r}$-grid with $\widetilde{r} \ne 0$. On the other hand, $\Gamma_2$ does not only have full lines, so this case cannot occur here. \\
 
   If $\Gamma$ is a dual grid, one dualizes the possibilities for $(\Gamma_1,\Gamma_2)$.\\

\ul{Now let $\Gamma_2$ have only ideal points.}  This means that we should have the duals of (a)--(b)--(c) included in the list of possibilities:

\begin{itemize}
\item[(a$'$)]
$s = 1$ and $t = 1$. In this case, $\Gamma_2$ is a line with two points, and $\Gamma_1$ is a root. This case is self-dual.
\item[(b$'$)]
$t = 1$ and $s$ is arbitrary. In that case, $\Gamma_2$ is a grid with parameters $(s - 1,1)$, and $\Gamma_1$ is a dual panel. 
\item[(c$'$)]
$s = 2$ and $t = 2$. In this case, $\Gamma \cong \mQ(4,2)$, and $\Gamma_2$ is a grid with parameters $(2,1)$. Also,  
$\Gamma_1$ is a dual grid with parameters $(1,2)$. This is also a self-dual case.
\end{itemize}

\ul{Next, suppose neither $\Gamma_1$ nor $\Gamma$ have only full lines.} 
By duality, we suppose the same about the points. \\

{\em First suppose $\Gamma$ is thick.}  Let $u$ be a point of $\Gamma_2$ which is not ideal, and let $U \I u$ be a line which is not a line of $\Gamma_2$ (it only has $u$ in common with $\Gamma_2$). 
Let $v$ be a point of $\Gamma_1 \setminus \Gamma_2$ on $U$. If there would be a line $W \I v$ which is a line of $\Gamma_2$, then obviously axiom (3) would be violated for $\Gamma_2$. So $v$, and any other point of $\Gamma_1$ on $U$ different from $v$, is an ideal point of $\Gamma_1$. In the same way one proves that $\Gamma_1$ has full lines. 
So any line incident with the points of $U$ different from $u$, is full by axiom (3) for $\Gamma_1$. Applying Lemma \ref{lemsub}, we conclude that $\Gamma_1$ is a thick
subquadrangle of $\Gamma$, and since it contains ideal points and full lines, it coincides with $\Gamma$. \\

{\em Suppose $\Gamma$ is not thick}. We suppose that $\Gamma$ is a quadrangle of order $(s,1)$ | the case with order $(1,t)$ is handled in a dual fashion. Call the reguli of $\Gamma$ respectively
$\mR_1$ and $\mR_2$. Suppose $U$ is a line of $\Gamma_1$ which is not a line of $\Gamma_2$. 

Suppose $u \I U$ is not  in $\Gamma_1$, and let $U'$ be the other line of $\Gamma$ that is incident with $u$.
Then axiom (3) in $\Gamma_2$ implies that no line intersecting $U$ except $U'$ is a line of 
$\Gamma_2$, as otherwise $U \in \Gamma_2$. The line $V$ cannot contain points of $\Gamma_1$, as otherwise axiom (3) forces $V$ to be in $\Gamma_1$. So $V$ is 
full in $\Gamma_2$, and we proceed as below. 

Now suppose that all points of $U$ are points of $\Gamma_1$.
Then if $\widetilde{U}$ is another line of $\Gamma_1$ in the same regulus,  it is also full. It follows easily that $\Gamma_1$ is an $(s + 1)\times r$-grid for some $1 \leq r \leq  s + 1$. 
If there is no such other line, $\Gamma_1$ would have order $(s,0)$, which is not possible, since then $\Gamma_1$ only has ideal points and lines. So we are in the former case. 
Let $U \in \mR_1$. If all lines of $\mR_1$ are in $\Gamma_1$, then $\Gamma_1 = \Gamma$, contradiction. So some 
lines of $\mR_1$ are in $\Gamma_2 \setminus \Gamma_1$, and hence $\Gamma_2$ is an $(s + 1)\times r'$-grid which shares all the lines of $\mR_2$ with $\Gamma_1$, and such that their union is $\Gamma$. 
\eop 

\medskip
\subsection{Irreducible varieties}

If $A$ is a commutative ring, then the points of the Zariski topology $\Spec(A)$ consist, besides $\emptyset$ and $A$, of the prime ideals of $A$. 
The prime ideals correspond, through the bijection
\begin{equation}
\frak{p}\ \mapsto\ \overline{\{\frak{p}\}}
\end{equation} 
to the closed irreducible subsets of $A$, that is, in the context of varieties, to the irreducible subvarieties. A closed subset in a topology $T$ is {\em irreducible} if it is not the union 
of two nonempty closed proper subsets $T_1$ and $T_2$. (Equivalently, all nonempty open sets are dense, that is, any two nonempty open sets have nonempty intersection.) 
In our setting, our  model base ``varieties'' are the generalized quadrangles, extended with the geometries satisfying axiom (3) (as described in this section). So we need to specify those 
generalized quadrangles \ul{which are not union of two proper subquadrangles}. This is where Theorem \ref{decthm} comes into play. 


\subsection{Affine examples}
\label{affine}

The above points are analogons of projective $k$-varieties in our theory, but we do not have analogons of affine $k$-varieties yet, such as affine $k$-spaces. If $k = \F_q$ for some prime power $q$, affine spaces have one point less per line than projective spaces, but introducing the affine analogons through this counting principle is a bad way for obvious reasons. Since affine varieties are complements of (large) closed varieties (and in particular affine spaces are complements of projective hyperplanes), there is a better way to do it. A {\em geometric hyperplane} $\mG$ in a generalized quadrangle $\Gamma$ (be it thick or thin), is a subquadrangle with the property that for each line $U \in \Gamma$, we have that either $U$ is a full 
line of $\mG$, or $U$ is incident with precisely one point of $\mG$. 
When $\Gamma$ is thick, and of order $(s,t)$,  three classes exist:

\begin{itemize}
\item
ovoids: these are point sets $\mO$ with the property that each line of $\Gamma$  contains one point of $\mO$; 
\item
subGQs $\mS$ of order $(s,t/s)$:  in fact, this is only literally true when $\Gamma$ is finite; when $\Gamma$ is not, we need to additionally ask that each line of $\Gamma$ meets $\mS$;
when $\Gamma$ is a GQ of order $(s,s)$, $\mS$ is a full $(s + 1)\times(s + 1)$-grid (with the additional property in the infinite case);
\item
maximal perp geometries (i. e., the special point is an ideal point). 
\end{itemize}

When one removes a geometrical hyperplane from a thick GQ, the resulting geometry is called {\em affine quandrangle} | see Pralle \cite{Pralle} for an axiomatic theory. We will also use the same nomenclature for thin quadrangles. In section \ref{varieties}, we will define general affine elements.

Below are the primal points of a point-line geometry with $\widetilde{(3)}$, and where every line has $s + 1$ lines. 
(with details if necessary about the geometrical hyperplanes); although we use the outcome of Theorem \ref{decthm}, we will explain our motivation   when we deviate from that result.  We also describe the affine geometries which arise by taking away one geometric hyperplane out of an ``old'' primal point.

\ul{List}:
\medskip

\begin{itemize}
\item[(0)]
the empty quadrangle $\varnothing$.
\item[(a)]
``ordinary points.''
\item[(b)]
full lines.
\item[(c)]
ideal perp geometries with full lines in thin or thick full (not necessarily proper) subquadrangles (or in other words, intersections of global ideal full perp geometries with thin or thick full (not necessarily proper) subquadrangles at a point of the latters). 
\item[(d)]
full grids (that is, grids for which all the lines are full). $\Big($Note that we deviate from Theorem \ref{decthm}; such grids can be decomposed according to (iii) of that theorem, but we have to allow such a decomposition because one has to see full $(s + 1)\times(s + 1)$-grids and the sub $(s \times s)$-grids which arise as in (iii), as analogons of varieties and open subvarieties. The relative viewpoint is important here: the lines of such $(s \times s)$-grids are not full, all miss a point, and only then can we identify them as open sets
in the full grids. Having agreed on that, full $(s + 1)\times(s + 1)$-grids are not the union anymore of disjoint nonempty closed sets if $s > 1$.$\Big)$
\item[(e)]
thick subGQs of order $(s,t)$, where $(s,t) \ne (2,2)$. Here, again,  ``order $(s,t)$'' means a.o. that each line of these thick subGQs are full. 
\end{itemize}

The affine counterparts are:  
\begin{itemize}
\item[(0-aff)]
the empty affine quadrangle $\varnothing$.
\item[(a-aff)]
``ordinary points.''
\item[(b-aff)]
affine quadrangles coming from full lines: full lines minus one point.
\item[(c-aff)]
affine quadrangles coming from ideal perp geometries: they arise from ideal perp geometries by removing one full line (and then one obtains, in general, disconnected examples consisting of full lines minus one point, which meet ``at infinity''), or an ovoid (example: a point set of the form $\{x, y\}^{\perp}$, where $x$ is the special point of the ideal perp geometry, $y$ 
is not collinear with $x$, and $y$ can be projected onto each line of the perp geometry (this is a condition!); other special example: the ovoid $\{ \mbox{special point} \}$), or a full perp subgeometry as described in (c).
\item[(d-aff)]
affine quadrangles coming from grids with only full lines: examples are described in Theorem \ref{decthm}[(iii)]; others arise by taking away an ovoid of a full grid.
\item[(e-aff)]
affine quadrangles constructed from examples in (e). \\
\end{itemize}

For a further discussion on how to distinguish between primal geometries and their affine versions, we refer to section \ref{laws}.


\subsection{Closed subsets (2)}

In the modern Zariski topology, for instance of the ring 
\begin{equation}
R := k[x_1,\ldots,x_m] 
\end{equation}
with  $k$ any field and $m \geq 1$, closed sets are defined in function of the ideals of $R$; 
if $I$ is an ideal, the corresponding closed set $C(I)$ consists of all prime ideal $\fp$ for which $\fp \supseteq I$. Translated to geometry, to $I$ corresponds a $k$-variety $V(I)$ (or $k$-scheme, but we will keep on working mostly with the terminology of varieties in this example) which is defined by the equations 
\begin{equation}
\{ P(x_1,\ldots,x_m) = 0\ \vert\ P \in I \},
\end{equation} 
and so the prime ideals in $C(I)$ correspond to irreducible varieties which are subvarieties of $V(I)$. (In the old Zariski topology, $k$ was taken to be algebraically closed, and $V(I)$ consisted of the $k$-rational points defined by the equations as above. The modern variation is much richer.) Having no underlying ring nor equations at our disposal, we define the closed subsets through 
the geometric translation: we first define ``ideal subgeometries $\Gamma$'' (this time, ideal in the algebraic sense),  and then define the corresponding closed sets $C(\Gamma)$ as 
\begin{equation}
C(\Gamma)\ \ :=\ \ \{ \mbox{prime geometry}\ \fp \ \vert\ \fp \leq \Gamma \}.
\end{equation}

In the present paper, it is not (yet) essential that one uses the classical Zariski approach, or the modern one: in the former case, we can define  the set of closed sets as being 
topologically generated by the base elements $\mS$, with $\mS$ ideal geometries in the sense of section \ref{ideal}. From an incidence-geometrical point of view, it might be more convenient to work with this setting, but in the present text, we will make it clear which approach we are using at what point. \\

We also need to say what our ``prime geometries'' are.

\subsection{The category $\mQ_\ell$}
\label{ql}

Let $\widetilde{\mQ_\ell}$ be the category of generalized quadrangles with $\ell + 1$ points per line, and \ul{thick}. 
Let $\widehat{\mQ_\ell}$ be the category of which the basic objects are the elements of $\widetilde{\mQ_\ell}$, and such that 
$\widehat{\mQ_\ell}$ is closed under taking of finite products, and finite disjoint unions. 

A typical object in $\widehat{\mQ_\ell}$ looks like:
\begin{equation}
\coprod_{i = 1}^m\Big(\bigotimes_{j = 1}^{N_i}\mQ^{(ij)}\Big), 
\end{equation}
where each $\mQ^{(ij)}$ is an element of $\widetilde{\mQ_\ell}$ (and $n, m$ positive integers). \\

The objects in $\mQ_\ell$ are now defined as \ul{the (not necessarily induced) point-line subgeometries in $\widehat{\mQ_\ell}$}.

\medskip
\subsection{Remark: generalized $n$-gons embedded in generalized $m$-gons}
\label{embnm}

Any generalized $n$-gon $\Gamma$ ($n \geq 3$) can be embedded in a generalized $m$-gon $\Gamma'$ for any $3 \leq m < n$; see $\Gamma$ 
as the starting configuration for a free $m$-gon $F(\Gamma) =: \Gamma'$, and add in the next step, for every pair of elements at distance $m + 1$, a path of length $m - 1$ only using new elements. After repeating this process step-by-step, in the limit one obtains a generalized $m$-gon (see e. g. \cite[\S 4.1]{HVMsurvey}, or \cite{POL}). 
If $\Gamma$ is finite, then $\Gamma'$ has countably infinite parameters, and 
in general, $\Gamma$ neither is full nor ideally embedded in $\Gamma'$.

\medskip
\subsection{Ideal subgeometries}
\label{ideal}

Fix an element $X$ in $\mQ_\ell$. 
The set of {\em ideal subgeometries} or {\em $I$-subgeometr\-ies} of $X$ is topologically generated by the full subgeometr\-ies (in which we allow subgeometries without lines). Again, such geometries do not necessarily enjoy axiom (3). Again, we do not ask that if $u$ and $v$ are distinct collinear points in a generalized quadrangle, and $\beta$ is an $I$-subgeometry containing $u$ and $v$, then $u$ and $v$ are also collinear in $\beta$. Note that sets of mutually noncollinear points are full (in a given element 
of $\mQ_\ell$) by definition.  


\begin{proposition}
Let $\mS_1$ and $\mS_2$ be any two $I$-subgeometries of an element $X \in \mQ_\ell$. Then
\begin{itemize}
\item[{\rm (1)}]
$\mS_1 \cap \mS_2$ is also an $I$-subgeometry in $X$;
\item[{\rm (1)$'$}]
$C(\mS_1) \cap C(\mS_2) = C(\mS_1 \cap \mS_2)$;
\item[{\rm (2)}]
$\mS_1 \cup \mS_2$ is an $I$-subgeometry in $X$;
\item[{\rm (2)$'$}]
$C(\mS_1) \cup C(\mS_2) \subseteq C(\mS_1 \cup \mS_2)$, but
we do not necessarily have that $C(\mS_1) \cup C(\mS_2) = C(\mS_1 \cup \mS_2)$. 
\end{itemize} 
\end{proposition}

{\em Proof}.\quad
Parts (1), (1)$'$ and  (2) are obvious. (For part (1), note that if $U$ is a line of $X$ in $\mS_1 \cap \mS_2$, then it is also a line in both 
$\mS_1$ and $\mS_2$, so it is a full line in these $I$-subgeometries, and hence their sets of incident points coincide.)\\

The first part of (2)$'$ also is. For the last part, suppose for instance that $X$ contains full $(\ell + 1)\times(\ell + 1)$-subgrids. 
If $\Gamma$ is such a grid, take $\mS_1$ to be one regulus, and $\mS_2$ the other. Then the prime quadrangle $\Gamma$ is not contained 
in $C(\mS_1) \cup C(\mS_2)$.  
\eop \\

Clearly, the counter example at the end of the proof is just an example of a very rich class of counter examples, and this behavior deviates from algebro-geometric Zariski. 
Also, one notes that in part (1) and (1)$'$, one can consider arbitrary intersections. 

\medskip
\begin{remark}[Closed sets and primes]{\rm 
The most interesting closed sets / $I$-subgeometries are those arising from subquadrangles | i.e., those that are essentially prime.}
\end{remark}

\medskip
\subsection{Example}
\label{exunionline}

Let $X$ be some element in $\mQ_\ell$, and suppose $U$ is a line in $X$; denote the set of points incident with $U$ as $U^*$. Suppose $X$ contains a set of 
full subGQs $\{ \mQ(u)\ \vert\ u \in U^* \}$, such that for each element $\mQ(u)$, we have that $\mQ(u)$ meets $U$ precisely in the point $u$. (Note that such examples 
can easily be constructed, both in the finite and infinite case, by considering classical quadrangles.) Then the $\mQ(u)$s are ideal subgeometries of $X$, 
and their union $\Omega$ is as well if $X$ is finite. Now $\Omega \cap U = U^*$, but in $\Omega$ the points in $U^*$ are not collinear. As a point-line 
geometry, $\Omega \cap U$ is full in $U$, so it is an ideal subgeometry of the point-line geometry $U$. 




\medskip
\subsection{Prime geometries/quadrangles, and spectrum}
\label{varieties}

Since we want to focus on (motivic) incidence-geometrical behavior, our choice of prime ideals probably is the one which deviates the most from the algebro-geometric approach. 
For us, the leading prime geometries which will be the points of spectra, are the prime generalized quadrangles in the sense of section \ref{points} (Theorem \ref{decthm}). This approach makes many ideal subgeometries $\mS$ in a fixed element $X$ of $\mQ_\ell$,  ``vanish'' under the map
\begin{equation}
\alpha:  \mS \ \mapsto\ C(\mS) = \alpha(\mS)
 \end{equation}
in the modern setting. If, e.g., $\mS$ would be a generalized $m$-gon with $m > 4$, then $C(\mS)$ does not contain geometries with finite girth, and in that sense, generalized polygons with  greater gonality than $4$ do not contribute to the  story. Exactly the same remark holds for ideal subgeometries with no ordinary quadrangles as subgeometries, and 
this principle restricts the initial generality of ideal subgeometries to those at least ``related'' to GQs in some sense.\\

We define the (old/modern) {\em Zariski topology} $\USpec(X)$ of $X$ to be the topology generated by the closed sets $C(\mS)$ with $\mS$ an $I$-subgeometry of $X$ and $X \in \mQ_\ell$. Note that this is the projective version. We also have two other types|here is the \ul{list}:
\begin{itemize}
\item[PROJ] 
{\em Projective version}:  
$X$ is an element of $\mQ_\ell$ and $\USpec(X)$ is defined as above; the most interesting examples are irreducible (such examples would correspond to projective varieties) and the model examples are generalized quadrangles. 
\item[AFF]
{\em Affine version}: 
$X$ is constructed by taking away a geometric hyperplane from an element $Y$ of projective type (and geometric hyperplanes are similarly defined as in the case of quadrangles); 
the most interesting examples are irreducible (such examples correspond to affine varieties) and the model examples are affine quadrangles.  The topology of $\USpec(X)$ is the subspace topology induced by $\USpec(Y)$.
\item[QUAS]
{\em Quasi-projective version}: 
$X$ is constructed by taking away a closed set (a sub-element of projective type) from an element $Y$ of projective type; 
the most interesting examples are irreducible (such examples correspond to quasi-projective varieties).  The topology of $\USpec(X)$ is the subspace topology induced by $\USpec(Y)$.

\end{itemize}

(Note that an irreducible closed set can become reducible after taking away a closed subset.)

\subsection{Generic point(s)}

Recall that a {\em generic point} of a topological space $\mT$ is a point whose closure is the entire space; note that if $\mT$ is irreducible, then it can have at most one such point. Also, no such point necessarily exists. If $X$ is a thick generalized quadrangle which is not in the list of Theorem \ref{decthm}, then $\USpec(X)$ is irreducible, and $X$ obviously defines a (unique) generic point in $\USpec(X)$. The same is true for every primal geometry $Y$ (that is, every primal geometry $Y$ defines a unique generic point of $\USpec(Y)$). 

{\rm 
\begin{remark}[Enriching]
As we will see in section \ref{baseex}, in some specific situations, it is possible to define a richer Zariski topology on GQs $X$, if such $X$ are defined ``in a more general class of quandrangles.'' More specifically, if $X$ is defined from a combinatorial object $\mO$ which lives in some projective space over a field $\F_{p^i}$, and the same construction works over field extensions $\F_{p^j}$ with $j$ contained in some infinite set $T \subseteq \mathbb{N}$, then we consider the quadrangle over the infinite field $\ell = \cup_{j \in T}\F_{p^j}$; the topology of $X$ is then defined through the action of $\Gal(\ell/\F_{p^i})$ on certain projective subspaces which are elements (points/lines) of the quadrangle over $\ell$ (mimicking the fact that for projective varieties, prime ideals (on ``level $\F_q$'') correspond to Galois-orbits of prime ideals over an algebraic closure of $\F_q$ | see section \ref{BASE} and section \ref{Galois}). Details can be found in section \ref{baseex}.
\end{remark}
}

\medskip
\section{Isomorphisms}
\label{isom}

To determine isomorphism classes between the objects which will arise from the defining relations is a delicate matter. To give one example,
suppose we consider a class $[ L ]$, with $L$ the incidence geometry of a line with $\ell + 1$ points. As an incidence geometry, its automorphism 
group is isomorphic to the symmetric group on $\ell + 1$ letters, since it has no structure at all. The same can be said about sets of points of the same size. 
So let $\Gamma$ and $\Gamma'$ be  isomorphic generalized quadrangles, and let $C$ and $C'$ be point sets in the respective quadrangles, of the same size. 
If we see $C$ and $C'$ as isomorphic objects, the scissor relation gives us that

\begin{equation}
[\Gamma \setminus C] = [\Gamma' \setminus C], 
\end{equation}
which  is in general absurd, since $C$ could, for instance, be a set of points on one and the same line, and $C'$ could be a partial ovoid. \\

We also refer to section \ref{lines1} for further (non-)examples in this discussion.\\


So we obviously need a better definition for ``isomorphisms'' and  ``automorphism groups'' of certain geometries in order to control such phenomena. 

Suppose $A$ and $B$ are closed subsets in $X$, an element of $\mQ_\ell$; then if $x$ and $y$ are different collinear points in $A, B$, the line $xy$ is also 
a line in $A \cap B$, so the closed set $A \cap B$ keeps seeing collinearity. So a morphism between closed sets must {\em keep track of collinearity}. 

Note that the following diagram must exist and commute, if $\gamma$ is an isomorphism between $X$ and $X'$:

\begin{center}
\item
\begin{tikzpicture}[>=angle 90,scale=2.2,text height=1.5ex, text depth=0.25ex]
\node (a0) at (0,3) {$X$\ \ };
\node (a2) [right=of a0] {};
\node (a1) [right=of a2] {\ \ $X'$};

\node (b0) [below=of a0] {$\Top(X)$\ \ };
\node (b1) [below=of a1] {\ \ $\Top(X')$};

\draw[->,font=\scriptsize,thick]
(a0) edge node[left] {} (b0)
(a1) edge node[right] {} (b1)
(a0) edge node[auto] {$\cong$} (a1)
(b0) edge node[below] {$\cong$} (b1);
\end{tikzpicture}
\end{center}

Here, $\Top(\cdot)$ stands for the choice of Zariski topology we will be working with.

\begin{observation}
If $\alpha \in \Aut(X)$ (where we see $X$ as a point-line geometry), then $\alpha$ preserves $\Top(X)$, so that $\alpha \in \Aut(\USpec(X))$. \eop
\end{observation}


\subsection{Isomorphisms between partial ovoids}

One of the first basic questions which arises in the context of local isomorphisms, is how we should define them for 
partial ovoids. So let $\mC$ be a partial ovoid of a point-line geometry $X$ satisfying axiom $\widetilde{(3)}$. If we see $\mC$ 
as a ``geometry on itself,''  then $\Aut(\mC)$ is isomorphic to the symmetric group on $\vert \mC \vert$ letters, which is in general 
different from $\Aut(X)_\mC$. Even more, if $\mC$ is a more general (induced) subgeometry of $X$, the same remark is true, although 
in general we have that 
\begin{equation}
\Aut(X)_\mC\Big/\mbox{kernel} \leq \Aut(\mC).
\end{equation}

For ``good embeddings,'' often equality holds. Example: with $k$ any field, we have the following equality for orthogonal quadrangles:  
\begin{equation}
\Aut(\mQ(5,k))_{\mQ(4,k)}\Big/\mbox{kernel} = \Aut(\mQ(4,k)).
\end{equation}

The group $\Aut(X)_\mC$ contains some information about the embedding $\mC \hookrightarrow X$. In Grothendieck groups of varieties, this is not 
the case. Of course, if we work with induced subgeometries, $\Aut(\mC)$ {\em does see} the collinearity of points and the concurrency of lines, 
but we will see that we want to require more. 

\subsubsection{Size does matter?}

Let $\mC$ and $\mC'$ be partial ovoids of the same size in $X$. In the naive approach of above, it follows that $\Aut(\mC) \cong \Aut(\mC')$, so 
we could see them as isomorphic partial ovoids. But obviously $\mC$ and $\mC'$ could have very different forms and positions, so that 
$\Aut(X)_\mC$ and $\Aut(X)_{\mC'}$ could highly differ. In the next section, we consider an example in which $\mC$ and $\mC'$ are 
subGQs, rather than partial ovoids.

\subsubsection{Important example} 

Let $X, X' \cong X$ be nonclassical Kantor-Knuth GQs of the same type and with the same parameters; and let $\mC$ be a $\mQ(4,q)$-subGQ in 
orbit $\Omega_1(X)$, and $\mC' \in \Omega_2(X')$. Then $\mC \cong \mC'$ as abstract GQs, but there is no isomorphism 
\begin{equation}
X \mapsto X'
\end{equation}
mapping $\mC$ to $\mC'$. So with this definition of isomorphism, we would have 
\begin{equation}
[X \setminus \mC] \ =\ [X' \setminus \mC']
\end{equation}
while the geometries $X \setminus \mC$ and $X' \setminus \mC'$ cannot be isomorphic as subgeometries. \\

$\Big(${\em Proof}.\quad 
$\mC$ and $\mC'$ are geometric hyperplanes, so $X \setminus \mC$ and $X' \setminus \mC'$ are affine quadrangles. By Pralle \cite{Pralle}, we 
know that $X \setminus \mC \cong X' \setminus \mC'$ if and only if there is an isomorphism $\alpha: X \mapsto X'$ such that 
$\alpha(\mC) = \mC'$. Contradiction.$\Big)$\\

We suggest to approach isomorphisms between ovoids in the following way. Let $X$ and $X'$ be isomorphic geometries with $\widetilde{(3)}$, and 
let $\mC$, $\mC'$ be partial ovoids in $X, X'$. Then an {\em isomorphism} $\gamma: \mC \mapsto \mC'$ must be induced by an isomorphism 
$\overline{\gamma}: X \mapsto X'$. The following diagram commutes (where $\iota, \iota'$ are canonical embeddings): 

\begin{center}
\item
\begin{tikzpicture}[>=angle 90,scale=2.2,text height=1.5ex, text depth=0.25ex]
\node (a0) at (0,3) {$X$\ \ };
\node (a2) [right=of a0] {};
\node (a1) [right=of a2] {\ \ $X'$};

\node (b0) [below=of a0] {$\mC$\ \ };
\node (b1) [below=of a1] {\ \ $\mC'$};

\draw[->,font=\scriptsize,thick]
(b0) edge node[left] {$\iota$} (a0)
(b1) edge node[right] {$\iota'$} (a1)
(a0) edge node[auto] {$\overline{\gamma}$} (a1)
(b0) edge node[below] {$\gamma$} (b1);
\end{tikzpicture}
\end{center}

One extra motivation is the following: in incidence geometry, an ovoid of a generalized quadrangle is not intrinsically defined, but defined {\em relative to a 
generalized quadrangle}; so actually, one should see \ul{an ovoid $\mO$ of a generalized quadrangle $\Gamma$ as a pair $(\mO,\Gamma)$}. The upshot is 
that in some sense, $\Aut(\Gamma)_\mO$ sees more structure than $\Aut(\mO)$: it keeps track of the embedding 
\begin{equation}
\mO \hookrightarrow \Gamma.
\end{equation} 

\begin{remark}{\rm
If $X \not\cong X'$, partial ovoids can never be isomorphic.}
\end{remark}

As soon as enough linearity comes into play in subgeometries, we pass to the classical definition of isomorphism, though. (More about this later.)

\subsubsection{Remark: Lines, every which way}
\label{lines1}

In scheme theory, a projective line over an arbitrary field $k$ is isomorphic to an irreducible conic section:
\begin{equation}
\mathbb{P}^1(k) \cong \mC(k). 
\end{equation}

Consider $X = \mQ(4,q)$ and noncollinear points $x$ and $y$. Then $\{ x,y\}^{\perp}$ is the 
rational point set of a conic section; now consider a line $U$ in $X$ as well. Then if we would adopt 
the viewpoint of Algebraic Geometry, we would have that 
\begin{equation}
[X \setminus \{x,y\}^{\perp}]\ =\ [X \setminus U]. 
\end{equation}

As subgeometries, $\{x,y\}^{\perp}$ and $U$ are very different in nature, so in our $K_0$, this identity should {\em not} hold. Note that 
incidentally, 
\begin{equation}
\Aut(\mQ(4,q))_U\Big/ \mbox{kernel}\ \cong\ \Aut(\mQ(4,q))_{\{x,y\}^{\perp}}\Big/ \mbox{kernel}.
\end{equation}

The kernels are different though: for $U$, it is the subgroup generated by all elations and all homologies with axis $U$; for $\{ x,y\}^{\perp}$, 
it depends on the characteristic. If $q$ is odd, the group is generated by one involution and the homologies with centers $x$ and $y$. 
If $q$ is even, it is a central extension of size $2$ of $\PGL_2(q) \rtimes \texttt{Gal}(\F_q/\F_2)$.  The kernel for $U$ is in general much 
larger than the kernel for $\{x,y\}^{\perp}$.\\

\subsection{Isomorphisms between lines}

Suppose $X = \mQ(4,q)$, and $X' = \hT_2(\mO)$, with $\mO$ an oval in $\bP^2(\F_q)$, and $\mO$ not a conic (see section \ref{T2O} for the details on the $\hT_2(\mO)$-construction). Let $U$ be a line 
in $X$ and $U'$ a line in $X'$. Then $U$ carries the structure of a projective line over $\F_q$ | both on the level of induced automorphism group, 
as purely geometrical (it is a line of the ambient projective space $\bP^4(\F_q)$). On the other hand, $\Aut(X')_{U'}$ fixes a 
point of $X'$, although geometrically one could still see $U'$ as a projective line over $\F_q$ (embedding in the ambient projective space $\bP^3(\F_q)$).  

For now, we want to define $\Aut(\mbox{line})$ as $\Aut(X)_{\mbox{line}}\Big/ \mbox{kernel}$, to keep track of the embedding. If $X \cong X'$, this 
seems like a natural solution, but what if $X \not\cong X'$?

\subsubsection{Idea 1 (explained for GQs)}
\label{idea1}

There must be thick subGQs $Y \leq X$ and $Y' \leq X'$ such that there is an isomorphism $\gamma: Y \mapsto Y'$ which maps 
$U$ to $U'$. 

For instance, let $X = \mQ(5,q)$ and $X'$ a nonclassical Kantor-Knuth generalized quadrangle of order $(q,q^2)$, $U$ any line of $X$ and 
$U'$ the special line of $X'$. Then we can find $\mQ(4,q)$-subGQs $Y \leq X$ and $Y' \leq X'$ which respectively contain $U$ and $U'$, so 
in the setting of IDEA 1, $U$ and $U'$ would be isomorphic.

\subsubsection{Idea 2 (explained for GQs)}

In IDEA 2, we see a line $U$ as a couple (permutation group) $(U,\Aut(X)_U\Big/\mbox{kernel})$. So two lines $U$ are isomorphic if the corresponding permutation groups are 
equivalent. 

Consider again the GQs $X$ and $X'$ of section \ref{idea1}, let $U$ be as before, and let $U'$ be a line which meets the special line of $X'$, but is different from it. 
In the setting of IDEA 1, the lines $U$ and $U'$ are isomorphic, but in the setting of IDEA 2, they are not (since $\Aut(\mQ(5,q))_U$ induces a transitive group on $U$, and $\Aut(X')_{U'}$ fixes the point $U \cap \mbox{special line}$). 

We prefer IDEA 2, since we think that the line $U'$ should be seen as a different geometry than the special line. (``The larger geometry defines the isomorphism.'') \\

Note that if $(X,Y) = (X',Y')$, the isomorphism in IDEA 1 leads to an isomorphism in IDEA 2.  Note also that the example in IDEA 1 also works in the setting of IDEA 2 (as it should); $U'$ 
locally looks like a projective line.

\subsubsection{An extra subtlety}

Suppose that $U$ (in $X$) and $U'$ (in $X'$) are isomorphic, but {\em become} nonisomorphic, or the other way around, if we \ul{enlarge} $X$ or $X'$? 
From the viewpoint of Algebraic Geometry this should not be a problem: for $k$ an appropriate field, there exist $k$-varieties which are not isomorphic as $k$-varieties, 
but become isomorphic after base extension. Simple example: let $k = \mathbb{R}$, and put $X = V(X^2 + Y^2 - Z^2)$, and $X' = \bP^1(\mathbb{R})$. Over $\mathbb{R}$ they are not isomorphic, but they are over $\mathbb{C}$.  So the relative context {\em is} important, and that is also a take we want to consider.

\subsubsection{Example}

Let $X$ be a GQ, and suppose that $L$ is a regular line, while $L'$ is a line which is not regular. If $\Aut(X)_L \Big/ \mbox{kernel} \cong \Aut(X)_{L'} \Big/ \mbox{kernel}$ (as 
permutation groups), then once we introduce scissor relations, in IDEA 2 we will obtain a relation
\begin{equation}
[X \setminus L ]\ =\ [X \setminus L'].
\end{equation}

Of course the geometries $X \setminus L$ and $X \setminus L'$ are far from isomorphic. If $X$ would be a ``rigid GQ'' | that is, a GQ with a 
trivial automorphism group | then the aforementioned relation even turns up for any two different lines $L, L'$. 

The example of this section show the wide difference between the notion of isomorphism in IDEA 1 and IDEA 2. As stated before: changing the basic 
notions of isomorphisms is a matter of taste (or  situation), but obviously the Grothendieck ring (cf. section \ref{grogro}) will capture very different 
properties. 

In the next subsection, we explain a third idea, which is a geometric counterpart of IDEA 2, while generalizing IDEA 1 in an intrinsic way.

\subsubsection{Idea 3 (explained for GQs)}
\label{Idea3perp}

Let $X$ be a GQ, and let $U$ be a line. For the sake of convenience, we dualize to obtain a GQ $X^D$ and a point $U^T =: u$. We introduce 
the local geometry of $X^D$ at $u$ as follows: 
\begin{itemize}
\item
its points are the points of $u^\perp$ ($u$ plays a special role);
\item
its lines are of two types: lines of $X$ incident with $u$, and perps $\{ u,v\}^{\perp}$, with $v \not\sim u$; 
\item
incidence is natural. 
\end{itemize}

Denote this geometry by $\mG(u)$. It captures the geometry of $\Gamma$ ``at distance at most $1$ from $u$.'' We use the dual notions for lines, and use 
the same notation. 

Now let $V$ be another line of $X$. In IDEA 3, we say that the lines $(U,X)$ and $(V,X)$ are {\em isomorphic} if the geometries $\mG(U)$ and $\mG(V)$ are isomorphic. For now, 
we require that an isomorphism between $\mG(U)$ and $\mG(V)$ preserves the types of the points (that is, we wish that an isomorphism $\alpha: \mG(U) \mapsto \mG(V)$  
sends $U$ to $V$).  

It is obvious that if  $(U,X)$ and $(V,X)$ are isomorphic in the sense of IDEA 1, then $\mG(U)$ and $\mG(V)$ are isomorphic in the above sense, so that 
$(U,X)$ and $(V,X)$ are isomorphic in the sense of IDEA 3. I conjecture that the converse is not true, that is: an isomorphism 
\begin{equation}
\alpha: \mG(U)\ \mapsto\ \mG(V)
\end{equation}
which maps $U$ to $V$ does not necessarily give rise (in one way or the other) to an automorphism $\widehat{\alpha}$ of $X$ which 
maps $U$ to $V$. Note that I do not ask that $\widehat{\alpha}$ induces $\alpha$. 

In some cases, it is true that such an $\widehat{\alpha}$ {\em does} exist --- even inducing a given $\alpha$ --- but then much more structural information is available 
about the cover 
\begin{equation}
\pi: X \ \mapsto\ \mG(U). 
\end{equation}


\subsection{Isomorphisms between grids}

Let $\Gamma$ be a grid in $X$, and $\Gamma'$ be a grid in $X'$. Then we say that $(\Gamma,X)$ is {\em isomorphic} to $(\Gamma',X')$ 
if:
\begin{itemize}
\item[$X \cong X'$] (first approach) and there is an isomorphism $\alpha: X \mapsto X'$ which maps $\Gamma$ to $\Gamma'$.
\item[$X \cong X'$] (second approach) and $\Gamma$ and $\Gamma'$ have the same parameters, and if moreover
\begin{equation}
\Aut(X)_\Gamma\Big/\mbox{kernel}\ \cong \ \Aut(X')_{\Gamma'}\Big/\mbox{kernel}, 
\end{equation}
where we mean isomorphic {\em as permutation groups}. 
Note that if $(\Gamma,X)$ and $(\Gamma',X')$ are isomorphic in the sense of the first approach, they are also isomorphic in the second sense. 
\item[$X \not\cong X'$] Same definition as in the previous item.
\end{itemize}

In what follows, we will prefer the second definition to the first one. 

\subsubsection{Example}

Let $X = (X,\omega,T)$ be a general TGQ. Suppose that $U$, $V$ and $W$ are distinct lines incident with $\omega$ which are 
fixed by $\Aut(X)$. Note that these lines are regular lines. 
One can think of a $\hT_2(\mO)$ with $\mO$ a rigid oval (i. .e., with trivial automorphism group) as an example. 

Let $A$, $B$ and $C$ the symmetry groups of $U$, $V$ and $W$. Let $\Omega_y =: \Omega$ be the subgroup of $\Aut(X)$ which fixes some point $y$ which is opposite $\omega$. Define 
$\Gamma$ as the full grid defined by $U$, $W$ and $y$; define $\Gamma'$ as the full grid defined by $V$, $W$ and $y$. Let $Y := \proj_Wy$.

Then in the setting of IDEA 2, $\Gamma$ and $\Gamma'$ are isomorphic. The proof is easy. For the sake of convenience, we will suppose that $\Omega$ acts 
faithfully on both $\Gamma$ and $\Gamma'$. Now note that $\Aut(\Gamma)$ is given by $AC\Omega$ (and that $AB \cap \Omega = \{ \id\}$), and that $\Aut(\Gamma')$ 
is given by $BC\Omega$ (and $BC \cap \Omega = \{ \id\}$). 

Now define a map
\begin{equation}
\phi: A \cap B:\ a \mapsto b_a,
\end{equation}
such that $ab_a^{-1}$ acts trivially on $W$. This map is well-defined, since $T$ is an abelian group, and since its subgroups $A$, $B$ and $AB$ act transitively on $W \setminus \{\omega\}$. 

Now define a map
\begin{equation}
\mu: AC\Omega \mapsto BC\Omega:\ ac\omega \mapsto b_ac\omega. 
\end{equation}

Then for elements $ac\omega$ and $a'c'\omega'$ in $AC\Omega$ we have that $\mu(ac\omega)\mu(a'c'\omega') = (b_ac\omega)(b_{a'}c'\omega')$. 
On $W$, $c$ and $c'$ act as the identity, so on $W$, we have that  $(b_ac\omega)(b_{a'}c'\omega')$ equals $(b_a\omega)(b_{a'}\omega') = b_ab_{a'}[\omega,b_{a'}^{-1}]\omega\omega'$; 
as $[\omega,b_{a'}^{-1}] \in B$, and as $a'$ and $b_{a'}$ act in the same way on $W$, we have that 
\begin{equation}
b_ab_{a'}[\omega,b_{a'}^{-1}]\omega\omega' = \mu(aa'[\omega,{a'}^{-1}]\omega\omega') = \mu(a\omega a'\omega') = \mu\Big((ac\omega)(a'c'\omega')\Big)
\end{equation}
on $W$. 

Similar reasoning leads to the fact that $\mu(ac\omega)\mu(a'c'\omega') = \mu\Big((ac\omega)(a'c'\omega')\Big)$ on $V$, and since the action of an automorphism of a 
grid is completely determined by the induced action on the reguli, it follows that $\mu$ is an isomorphism of groups. 

The isomorphism of grids which yields the required equivalence of permutation actions, is the following (and is handed by the proof of the group 
isomorphism).  

\begin{equation}
\left\{
                \begin{array}{ccc}
                \Gamma &\mapsto &\Gamma' \\
                w \I W &\mapsto &w \I W \\
                u \I U &\mapsto &\proj_V(\proj_Yu)\\
                y' \I Y &\mapsto &y' \I Y \\ 
      
              \end{array}      \right.
              \end{equation}

The reader is invited to see what happens if one lets $U, V$ and $W$ be part of other than trivial line orbits (of lines incident with $\omega$).

\subsubsection{Example}

Let $k$ be a field and consider the orthogonal quadrangles $\mQ(4,k)$ and $\mQ(5,k)$. All the lines of these examples are regular lines. 
Let $\Gamma$ be a full subgrid in $\mQ(4,k)$, and $\Gamma'$ be a full subgrid of $\mQ(5,k)$. Then $\Gamma$ and $\Gamma'$ are 
isomorphic in the approach of IDEA 2; in each of the cases, the automorphism group is isomorphic to $\Big(\mathbf{P\Gamma L}_2(k) \times \mathbf{P\Gamma L}_2(k)\Big) \rtimes 
C_2$ in its natural action: both grids are products of two projective lines $\mathbb{P}^1(k)$. (The kernels of $\Aut(\mQ(4,k))_\Gamma$ and $\Aut(\mQ(5,k))_{\Gamma'}$ in their action on $\Gamma$ and $\Gamma'$ 
are different in general.)


\subsubsection{Example --- CAUTION}

Let $X$ be a GQ, and suppose $\Gamma$ and $\Gamma'$ are full grids in $X$ (so that they have the same parameters). Suppose furthermore that 
\begin{equation}
(\Gamma,X)\ \cong\ (\Gamma',X)
\end{equation}
in our second approach. Then by the scissor relations that we will impose on the Grothendieck ring to be defined (cf. section \ref{grogro}), we obtain that 
\begin{equation}
[X \setminus (\Gamma,X)]\ =\ [X \setminus (\Gamma',X)].
\end{equation}

If the geometries $[X \setminus \Gamma]$ and $[X \setminus \Gamma']$, which are generalizations of affine quadrangles, would happen to be isomorphic by an isomorphism $\beta$, 
then if the work of Pralle would also apply to this more general case, then there would exist an automorphism $\widehat{\beta}$ of $X$ which sends 
$\Gamma$ to $\Gamma'$, and which induces  $\beta$. So under these assumptions, we would have that $\Gamma$ and $\Gamma'$ are isomorphic in the sense of 
our first approach. 

\subsubsection*{Question}

Suppose $X, \Gamma$ and $\Gamma'$ are as in the previous paragraph, and suppose $\beta: X \setminus \Gamma\ \mapsto\ X \setminus \Gamma'$ is an 
isomorphism. Does there exist and automorphism $\widehat{\beta}$ of $X$ which induces $\beta$?

\subsubsection{Important Remark}

In the Grothendieck ring of $k$-varieties $K_0(\mV_k)$ ($k$ a field), the equality
\begin{equation}
\label{xisy}
[X]\ =\ [Y]
\end{equation}
with $X$ and $Y$ $k$-varieties, does not necessarily imply that $X$ and $Y$ are isomorphic. In fact, it is an open question that $X$ and $Y$ be birationally 
equivalent under the equality $[X] = [Y]$. 

It is extremely important to note that likewise, (\ref{xisy}) does not imply that $X$ and $Y$ have isomorphic automorphism groups. As one of our guiding 
examples, this is important to keep in mind.

\subsection{Isomorphisms between spans/perps and lines}

If we allow dualities as isomorphisms on the level of simple examples, then in some circumstances stars 
can be isomorphic to lines, and perp geometries to dual perp geometries, etc. In this paper, we do not consider 
dualities as isomorphisms, but that depends entirely on the taste of the reader (or of the problem which has to be handled). 

On the other hand, allowing dualities as isomorphisms, one would have to deal with rather unnatural decompositions in the Grothendieck 
ring which is to be constructed! And we would end up with a bad notion of Krull dimension (see section \ref{Krull} for examples of dual generalized quadrangles with 
different dimensions).


\subsection{Intrinsic automorphisms}

In each of the cases considered above, we have a quick look at the {\em automorphism groups} of the objects in the various settings. We use the notation 
used in the previous paragraphs.

\subsubsection{Partial ovoids}
\label{sspo}

If $(\mO,X)$ is a partial ovoid, then $\Aut(X)_{\mO}$ is the automorphism group of $(\mO,X)$. Note that we do not 
mod out a possible kernel here | the kernel of the action of $\Aut(X)_{\mO}$ on $\mO$ contains information about the embedding 
$\mO \hookrightarrow X$. 

\subsubsection{Lines, IDEA 2}

If $(U,X)$ is a line, then $\Aut(X)_U\Big/\mbox{kernel}$ is the automorphism group of $(U,X)$ in the setting of IDEA 2. The remark in \S \ref{sspo} about the kernel 
does not hold here, due to the fact that we have required less information about the embedding in the defining data of a line.

If $X = U$, then $\Aut(U)_U$ is isomorphic to the symmetric group on the points incident with $U$.

\subsubsection{Lines, IDEA 1}

For all $Y \leq X$, $U$ in $Y$ as a full subgeometry, an automorphism of $Y$ fixing $U$ would define an automorphism of $U$. One must decide when 
two such automorphisms $\alpha \in \Aut(Y)_U$ and $\beta \in \Aut(Y')_U$, where $Y' \leq X$ and $U$ is fully contained in $Y'$, define the same 
automorphism of $U$. One possible natural approach would be to say that they define the same such automorphism if $\alpha$ and $\beta$ 
agree on (some subset of) $Y \cap Y'$. 

So let $\mY$ be the set of all geometries $Y$ as above (supposed sufficiently linear), and for $Y, Y' \in \mY$, write that $Y \preceq Y'$ if 
$Y'$ is a subgeometry of $Y$, and if we can restrict each element of $\Aut(Y)_U$ to $Y'$. The latter property is needed to make 
$\Big(\mY,\preceq\Big)$ into a directed system through the restriction of automorphisms. Obviously, in general this system does not behave well: if we consider elements $\widetilde{Y}$ and 
$\widetilde{Y'}$ in $\mY$, there is no obvious upper bound for both geometries in this setting. So we relax the definition to the following: we say that $\alpha$ and $\beta$ (notation of the previous paragraph) 
are equal if they agree on {\em some subgeometry} $Y'' \in \mY$, $Y, Y' \preceq Y''$. Note that this is equivalent with saying that $\alpha = \beta$ if and only if 
they induce the same action on $U$ (as $U \in \mY$).  It follows that 
\begin{equation}
\varinjlim_{\Big(\mY,\preceq \Big)}\Aut(Y)_U \ \cong\ \Aut(U),
\end{equation}
where the latter is the naive automorphism group of $U$ (without seeing the embedding $U \hookrightarrow X$). 

In the case that $X = U$, we obtain the same automorphism group (and the same action, of course) as in the previous subsection.

\subsection{Grids}

If $(\Gamma,X)$ is a grid, then by our very definition, $\Aut(X)_\Gamma/\mbox{kernel}$ is the automorphism group of $\Gamma$. If $X = \Gamma$, we 
obtain a group isomorphic to $\texttt{S}_u \times \texttt{S}_v$ if $\Gamma$ is an $(u \times v)$-grid with $u \ne v$; if $u = v$, we obtain  $\Big(\texttt{S}_u \times \texttt{S}_v\Big)\rtimes C_2$.

\medskip

\section{Trace geometries} 
\label{tracegeom}

Motivated by the idea in subsection \ref{Idea3perp}, we introduce yet another approach to the isomorphism problem of thin subgeometries of 
(say) generalized quadrangles. 

So let $X$ be a thick generalized quadrangle with $\ell + 1$ points per line. Let the full subgeometry $\mT$ be of one of the following types: (a) partial ovoid; (b) perp geometry with base point $u$; (c)  grid. Define a point-line geometry $\hTr(\mT,X)$ as follows.

\begin{itemize}
\item
Points are the points of $\mT$.
\item
Lines are of (at most) two types: (i) the lines of $\mT$, and (ii) the point sets $v^{\perp} \cap \mT$, where in the latter expression $\mT$ is considered as a point set, 
and $v$ is a point in $X \setminus \mT$.
\item
Incidence is natural.
\end{itemize}

Note that if $\mT$ is of type (b), and if the geometry is not ideal, then we might pick a point $v \sim u$ in $X \setminus \mT$: the point set $v^{\perp} \cap \mT$ then 
equals $\{ u \}$, but we do not consider such a set as a line.

The idea of this construction is obvious: we \ul{enrich} the geometry $\mT$ with data coming from the embedding $\mT \hookrightarrow X$. 

\subsection{Trace isomorphism}

Say that $\mT$ and $\mT'$ (both subgeometries of the aforementioned type (a), (b) or (c)) are {\em trace isomorphic} if their {\em trace geometries} $\hTr(\mT,X)$ and $\hTr(\mT',X)$ 
are isomorphic as point-line incidence geometries. We distinguish two sub-types of such isomorphisms: on the one hand, {\em $p$-trace isomorphisms} which 
keep track of the lines of type (i) (so which map lines of type (i) to lines of type (i)), and general trace isomorphisms, which do not necessarily have this property. \\

Note that $p$-trace isomorphisms preserve the type of the geometry. For general trace isomorphisms, the question needs a more detailed analysis. In any case, 
a geometry of type (b) contains a number of concurrent lines which meet all other lines, and if it is not a projective plane, they are they only lines with this property. 
A geometry of type (c) contains two parallel classes of lines which meet all other lines, so such a geometry can never be isomorphic to a geometry of type (b). 
It is less clear to me how these geometries compare to geometries of type (a), certainly in the infinite case. In the finite case, we could give a simple number theoretic argument to exclude possible trace isomorphisms. If a geometry of type (a) (living in a GQ of finite order $(s,t)$) would be isomorphic to a geometry of type (b) or (c) (living in a GQ of finite order $(a,b)$), then obviously $t + 1$ (the line size in type (a)), must equal $a + 1$ (since in type (b) or (c), there are lines of size $a + 1$). In case of geometries of type (b), we also need that $a = b$ for such an isomorphism to exist. Comparing the number of points of the geometry in case of type (b), we see that $s^2 + 1 = (s + 1)s + 1$, which is false. In case of type (c), we have that $s^2 + 1 = (s + 1)^2$, which is also false.\\

The following question is obviously important.

\begin{question}
Are all trace isomorphisms $p$-trace isomorphisms?
\end{question}

I claim that the answer is negative. I will describe a general class of counter examples for geometries of type (b); suppose the base $u$ is a projective point in $X$; as we know, this means by definition that $\hTr(\mT,X)$ is a projective plane, and such examples exist in abundance in both the finite and infinite case. {\em Example}: let $k$ be a commutative field, and consider the orthogonal quadrangle $\mQ(4,k)$. Now suppose that the automorphism group of the plane $\hTr(\mT,X)$ does not fix $u$ (and let $\gamma$ be an automorphism not fixing $u$); again, the class of orthogonal quadrangles $\mQ(4,k)$ gives such examples since the planes are isomorphic to $\mathbb{P}^2(k)$ (group acts transitively on the points). Now, as soon as we consider two such geometries $\hTr(\mT,X)$ and $\hTr(\mT',X)$ and one isomorphism 
\begin{equation}
\beta:\  \hTr(\mT',X)\ \mapsto\ \hTr(\mT,X)
\end{equation}
which preserves line types, the isomorphism $\gamma \circ \beta$ does not. \\

In the finite case, these might be the only examples. Let $X$ have finite parameters $(\ell,m)$. Suppose $\gamma$ is as above. Then the lines incident with $\gamma(u)$ are: the line $u\gamma(u)$ and $m$ lines of type (ii). Since the number of points per line of type (ii) equals $m  + 1$, and since the lines of type (i) have $\ell + 1$ points, we know that $\ell = m$. 
The lines on $u$ have a particular property: they mutually intersect precisely in $u$; so $\gamma(u)$ also has this property. Now let $u''$ be any point in $\hTr(\mT,X)$ which is not collinear with $u''$ in $X$. Then there is at least one, and hence precisely one, line of type (ii) containing both $\gamma(u)$ and $u''$. Now in principle, we could have constructed $m$ lines containing $\gamma(u)$ and $u''$ | each one per point in $\{ \gamma(u), u'' \}^{\perp} \setminus \{u\}$, and we conclude that they all coincide. 
It follows that the pair $\{ \gamma(u),u''\}$ is regular in $X$. This might be an indication that $u$ is regular, and $\hTr(\mT,X)$ a projective plane.\\

I suspect that in the infinite case, there are (many) other types of examples. \\

Now we turn to grids; so suppose $\mT$ is of type (c). Suppose that $\gamma$ is an automorphism of $\hTr(\mT,X)$ which does not fix the set of lines of type (i). Note that in $\hTr(\mT,X)$, every line of type (ii) meets {\em all} lines of type (i). So if $U$ is a line of type (i) and $\gamma(U)$ is not, then $\gamma(V)$ is also not of type (i) for each line $V$ in the same regulus $\mR(U)$ of $\mT$ as $U$. 

Now let $R$ be a line of type (ii) which is not contained in $\gamma(\mR(U))$. Then if $R$ would meet some line $W \in \gamma(\mR(U))$ in at least two points, it cannot meet all the lines of $\gamma(\mR(U))$, contradiction. So $R$ must coincide with $W$. In other words, if $u, v$ are different points of $W$, and $w \in \{ u, v\}^{\perp} \setminus \mT$, then $w$ is collinear with all points of $W$.
This fact has strong combinatorial implications. Suppose for instance that 
$X$ is a finite GQ of order $(s,s)$. Then the points of $R$, together with the points in $R^{\perp}$ form a complete dual $(s - 1) \times (s + 1)$-grid, which is a very uncommon object in such GQs. (Almost no such examples are known | see the final remark in \cite{KTcomparcs}.)
It can be shown (see the final remark in the paper \cite{KTcomparcs}) that the existence of complete dual $(s - 1) \times (s + 1)$-grids gives rise to the existence of complete $(s^2 - s)$-arcs, and by \cite{KTcomparcs} we have the following result. 

\begin{theorem}[Theorem 4.1, \cite{KTcomparcs}]
\label{4.1}
Let $\mS$ be a known thick finite GQ of order $(s, t)$, and suppose $\mS$ has a complete $(st - t/s)$-arc. Then we necessarily have one of the following possibilities.
\begin{itemize}
\item
$\mS \cong \mQ(4, 2)$ and up to isomorphism there is a unique example.
\item
$\mS \cong \mQ(5, 2)$ and up to isomorphism the arc is unique.
\item
$\mS \cong \mQ(4, q)$ with $q$ odd. 
\end{itemize}
\end{theorem}

The construction of complete $(s^2 - s)$-arcs is as follows: start from a complete dual $(s - 1) \times (s + 1)$-grid in $X$, a GQ of order $(s,s)$. Let $\mR$ be the dual regulus of size $s + 1$; then the $2(s + 1)$ lines which are not lines of the dual grid form a full grid with parameters $(s,1)$. Now consider all the points which are not contained in that grid, and which are not collinear with a point of the dual grid: there are $(s^2 - s)$ such points, and it is easy to see that they form a complete arc. If we now inspect the known examples of finite GQs (Theorem \ref{4.1}, first and last case), we see that only the last case is possible with $q - 1 = 2$. \\

The infinite case is different, and not much is known. Still, the construction of extremal complete arcs from extremal dual grids can be adapted. We sketch the idea, and leave the details to the reader. So let $\Gamma$ be a full grid in a generalized quadrangle $\mS$, with the following property:
\begin{quote}
each line of $\mS \Gamma$ meets $\Gamma$ (necessarily in a unique point).
\end{quote}

This property guarantees that $\Gamma$ is ``as big as possible'' in $\mS$. 

Now suppose $\widetilde{R}$ is a set of mutually noncollinear points in $\mS \setminus \Gamma$ such that if $u \ne v$ in $R$, then $u^{\perp} \cap \Gamma = v^{\perp} \cap \Gamma = \widetilde{R}'$, and if $r$ is a point in $\widetilde{R}'$, then each line incident with $r$ and not contained in $\Gamma$, contains a point of $\widetilde{R}$. 
Then $\widetilde{R}$ and $\widetilde{R}'$ form the point reguli of the analogon of a complete dual $(s - 1) \times (s + 1)$-grid in $\mS$ (with the additional property that one of them lies in a full grid). 

Now perform the exact same construction as in the finite case to obtain the analogon of a complete $(s^2 - 1)$-arc.

\subsection{Embedding properties}

For our purposes, only $p$-trace isomorphisms are important, since they keep track of embeddings. 

The following observation is straightforward. 

\begin{observation}
If $\mT$ and $\mT'$ are isomorphic in the setting of IDEA 1 (so as induced subgeometries of $X$), then they are also $p$-trace isomorphic. \eop
\end{observation}

\begin{observation}[Embedding property]
Let $X \leq \widehat{X}$, with $\widehat{X}$ a generalized quadrangle in which $X$ is full. Then $\hTr(\mT,X) \ne \hTr(\mT,\widehat{X})$. \eop
\end{observation}

\begin{question}
Can $\hTr(\mT,X) \cong \hTr(\mT,\widehat{X})$?
\end{question}

This question is interesting and most relevant for our discussion. Quite strangely, the answer is ``yes.'' Consider a GQ $X'$ with full grid $\mT$, and let 
\begin{equation}
\alpha:\ X'\ \hookrightarrow\ X' 
\end{equation}
be a full embedding of $X'$ in itself. Then obviously $\hTr(\mT,X) \cong \hTr(\mT,\widehat{X})$. We now describe such an example.
Consider a quadric in an infinite-dimensional real vector space $V$
consisting of vectors of the form $(X_0, X_1, X_2, X_3, Y_0, Y_1, Y_2, \ldots)$, with quadratic form
\begin{equation}
q:\ X_0 X_1\ +\ X_2 X_3\ +\ Y_0^2\ +\ Y_1^2 +Y_2^2\ + \cdots 
\end{equation}

This yields a generalized quadrangle $X'$, whose lines are projective lines over $\mathbb{R}$. Consider the morphism
\begin{equation}
\phi:\ V\ \to\ V:\ (X_0, X_1, X_2, X_3, Y_0, Y_1, Y_2, \ldots)\  \mapsto\  (X_0, X_1, X_2, X_3, 0,  Y_0, Y_1, \ldots). 
\end{equation}

This morphism leaves the quadratic form invariant, and hence induces an injective morphism
$X' \mapsto X'$, which is not a bijection, and fixing certain lines pointwise (since the
$(X_0, X_1, X_2, X_3)$ piece is left invariant). Also, $X'$ is orthogonal, so all its lines are regular. \\

More generally, say that $\mT$ and $\mT'$ (both subgeometries of the aforementioned type) are {\em trace isomorphic} if their {\em trace geometries} $\hTr(\mT,X)$ and $\hTr(\mT',X')$ 
are isomorphic as point-line incidence geometries.

\subsection{Parallel points (and lines of type (i))}
\label{PPs}

Note that it is not essential to identify the lines of (i) if one keeps track of the embedding $\mT \hookrightarrow X$: for, let $x, y$ be different points of $\mT$. If $x$ and $y$ 
are not collinear in $\mT$, then consider a point $z \in \{ x, y\}^{\perp} \setminus \mT$. Such points always exist. Then $z^{\perp} \cap \mT$ is a line of 
$\hTr(\mT,X)$ which connects $x$ and $y$, making these points collinear. If $x$ and $y$ are collinear in $\mT$, they are only collinear through the unique line of $X$ connecting them. So if we would only consider the lines of type (ii) in $\hTr(\mT,X)$, two different points only are 
not collinear precisely if they are collinear in $\mT$. Moreover, maximal sets of mutually noncollinear points in $\hTr(\mT,X)$ \ul{in this setting} (``parallel classes of points''), precisely are the point sets of the lines in $\mT$.

\subsection{The thick case}  

When $\mT$ is a full thick subquadrangle, there is an interesting application as well. Define $\hTr(Y,X)$, with $\mT = Y$  a thick full subquadrangle of $X$, in the same way as above. The geometry $\hTr(Y,X)$ has two line types: lines of $Y$ and subtended ovoids of $Y$ from points in $X \setminus Y$. 

First of all, remark that the consideration in subsection \ref{PPs} remains true. If we would leave out the lines of type (i), the parallel classes of points recover the former GQ-lines. 
For the sake of convenience, define $\underline{\hTr(Y,X)}$ as $\hTr(Y,X)$ without the lines of type (i). (Clearly, these geometries have the same automorphism groups, and are essentially the same.)

The geometries $\underline{\hTr(Y,X)}$ have been extensively studied in literature, and have been denoted ``$\mA$'' in the paper \cite{covers}. Many (difficult) open problems occur in this theory. In \cite{covers}, it is furthermore suppose that each subtended ovoid/line in $\hTr(Y,X)$ is subtended by a constant number of points.

\subsection{Isomorphism/embedding problems}

\begin{question}
(When) can $X$ be reconstructed from $\hTr(\mT,X)$?
\end{question}

This answers is most probably too difficult without further assumptions. Also important in the tick case. 

\begin{question}
(When) can automorphisms of $\hTr(\mT,X)$ be extended to automorphisms of $X$?
\end{question}

Same remark. 

\begin{question}
How do $p$-trace isomorphisms between $\hTr(\mT,X)$-geometries related to isomorphisms in the setting of IDEA 2?
\end{question}

\medskip
\section{Theory of lines, and birationality}
\label{theline}

\subsection{}
Let $A$ be a line in a generalized quadrangle $\Gamma$. The induced topology is the cofinite one: the open sets are 
the empty set, $A$, and the complements of finite point sets. 

\subsection{}
Let $B$ be an subset of $A$; then $B$ naturally defines a line, which we call {\em subline} of $(U,\Gamma)$. It is defined as $(B,{\Big(\Aut(\Gamma)_A\Big)}_B\Big/\mbox{kernel})$. 
For the sake of convenience, we assume that $B$ has at least two points. Now suppose the element $\beta \in \Aut(\Gamma)$ stabilizes $B$. As $B$ has at least 
two points, it follows that $U$ also is stabilized by $\beta$. So 
\begin{equation}
(B,{\Big(\Aut(\Gamma)_A\Big)}_B\Big/\mbox{kernel})\ \cong\ (B,{\Aut(\Gamma)}_B\Big/\mbox{kernel}). 
\end{equation}

\subsection{}
Call two elements $X, Y$ in $\mQ_{\ell}$ {\em birationally equivalent} if there are dense open sets $U \subset X$ and $V \subset Y$ such that 
$U$ and $V$ are isomorphic (as incidence geometries). (Recall that a set $D \subset Z$ in a topological space $X$ is {\em dense} if it meets 
every nonempty open set of $Z$.) 

\subsection{}
Now let $X$ and $Y$ be GQs in $\mQ_\ell$, and consider lines $(U,X)$ and $(V,Y)$. By the above, these lines are birationally equivalent if 
there are isomorphic open dense sublines. \\

{\bf CASE $\ell$ finite.}\quad 
the topologies are discrete, so the only dense open sets are the entire sets. So $(U,X)$ and $(V,Y)$ are birationally equivalent {\em if and only if 
they are isomorphic.}\\

{\bf CASE $\ell$ infinite.}\quad
All nontrivial open sets have infinite size (with finite complements), so any neighborhood of a point contains a point of any nonempty open set. 

For the rest of this subsection, we suppose that $\ell$ is infinite. 

\begin{observation}
\label{obsbirlines}
Let $X= Y$. Then $(U,X)$ is birationally equivalent to $(V,Y)$ if and only if $(V,Y)$ has an open subline of type $(W,W)$, with $W \subseteq V$ (the induced 
group being $\mathrm{Sym}(W)$). \eop \\  
\end{observation}  

\begin{corollary}
\label{corbirlines}
In the situation of Observation \ref{obsbirlines}, we have that $Y$ cannot be a TGQ.
\end{corollary}

{\em Proof}.\quad
Suppose $Y = (Y^x,T)$ is a TGQ by way of contradiction, and let $(V,Y)$ have an open subline of type $(W,W)$. We first suppose that $V$ is not 
incident with $x$. Let $A$ be the group of symmetries with axis $\proj_xV =: [V]$. It is an abelian subgroup of $T$, and it fixes every line of $[V]^{\perp}$. 

Suppose that $x$ is the only translation point of $Y$. 

Let $S$ be a subgroup of $\Aut(Y)$ which induces $\mathrm{Sym}(W)$ on $W$; as it fixes 
$x$ and $V$, it also fixes $[V]$. Whence 
\begin{equation}
[A,S] \leq A. 
\end{equation}

Consider different points $w, w' \in W$, and an element $\omega$ in $S$ that fixes both $w$ and $w'$; if $a \in A$ is the element which maps $w$ 
to $w'$, then $[\omega,a]$ fixes both $w$ and $w'$. As $[\omega,a]$ is a symmetry with axis $[V]$, it follows that $[\omega,a] = \id$. On the other hand, for any 
point $w'' \in W \setminus \{ w,w'\}$, we can find such an $\omega$ for which ${w''}^{\omega} \ne {w''}^{a\omega a^{-1}}$. For such an $\omega$, we have that 
$[\omega,a] \ne \id$, and this is a contradiction. \\

If $V \I x$, let $A$ be the group of symmetries with axis $C \ne V$, $C \I x$, and repeat the argument of above to obtain a contradiction. \\

Now suppose that $Y$ has more than one translation point. First suppose that there is a line $L$ of translation points, and that all translation points are incident with $L$. Then $L$ is fixed by $S$. If $V = L$, we can take $x \in W$ and consider the action of $S_x$ on the set $W \setminus \{x\}$; then we just proceed as before. Now let $V \sim L \ne V$. Then $S$ fixes the translation point $L \cap V$, and we can also proceed as before. Finally, let $V \not\sim L$. Project $W$ on $L$ to obtain the point set $W'$. Then $(S,W)$ is permutation equivalent to $(S,W')$, and we can yet again proceed as before. \\

Next, suppose that all points of $X$ are translation points. We provide two different methods, the second one of which is elementary.\\

{\bf METHOD 1}.\quad 
Since all points are translation points, we know by Tits \cite{TiQuad} (or by elementary combinatorial methods) that $Y$ is a Moufang quadrangle. By Tits and Weiss \cite{TiWe}, $Y$ is orthogonal, and hence its lines are projective lines. This means that for all lines $D$ in $Y$, we have that 
\[    \Aut(Y)_D/N \ \cong\ \PGL_2(k) \rtimes \Aut(k),           \]
with $N$ the kernel of the action of $\Aut(Y)_D$ on $D$. This leads to a contradiction, since now obviously $V$ cannot obtain a cofinite set $W$ such that $\mathrm{Sym}(W)$ is induced by a subgroup of $\PGL_2(k) \rtimes \Aut(k)$. \\

For obvious reasons, this method is not desirable, so we also describe a much more basic approach.\\ 

{\bf METHOD 2}.\quad 
Let $x$ be a point in $W$, and consider the permutation group $(S_x, W \setminus \{x\})$. Let $C \I x$ be a line different from $V$, let $y \I C$ be different from $x$, and let $V' \I y$ be different from $C$. Also, let $A$ be the group of symmetries with axis $C$. 
Consider the group $\Aut(Y)_{[V]}$ of all automorphisms of $Y$ which fix $V$ pointwise; since all points of $Y$ are translation points, one easily observes that it acts transitively on the lines opposite $V$. It follows readily that we can find an automorphism group $\widehat{S}$ 
in $\Aut(Y)_{(V,x,C,y,V')}$ which induces $S_x$ on $W \setminus \{x \}$. 

Then 
\[   [\widehat{S}, A]\  \leq\ A,     \]
and we can repeat the argument as in the case of one translation point. 
\eop \\

{\bf Second proof of Corollary \ref{corbirlines} for linear TGQs}.\quad
There is an entirely different ``projective'' proof of the previous corollary, if one supposes the TGQ to have a 
projective representation (which is always the case in characteristic $\ne 0$). So suppose $(Y^x,T)$ is linear, and let $k$ be the kernel of 
the TGQ.  We use the terminology from the paper \cite{TGQkernel}. \\

Represent $(Y^x,T)$ in the projective spave $\mathbb{P}$ over $k$ (the corresponding egg lying in the hyperplane $\mathbb{P}'$). 
To $[V]$ corresponds a subspace of $\mathbb{P}$ which is not contained in $\mathbb{P}'$. We will freely use the next result from \cite{TGQkernel}: 

\begin{observation}(Transfer of automorphisms \cite[Observation 6.5]{TGQkernel})
\label{transfer}
Let $(\Gamma,u,C)$ and $(\Gamma',u',C')$ be TGQs, and suppose
\[ \gamma:  (\Gamma,u,C)\ \mapsto\ (\Gamma,u',C') \]
is an isomorphism (which maps $u$ to $u'$). Suppose moreover that the kernel of $(\Gamma,u,C)$ contains a prime field $\wp$. Then $\gamma$ induces a semilinear isomorphism between the projective spaces $\widetilde{\bP}$ and $\widetilde{\bP}'$ which arise from interpreting $C$, respectively $C'$, as a vector space over $\wp$, respectively $\wp^\gamma$. Moreover, this isomorphism maps the egg of $\Gamma$ in $\widetilde{\bP}$ to the egg of $\Gamma'$ in $\widetilde{\bP}'$.
\end{observation}

We only consider the generic situation in which $V$ is not incident with $x$, and in which $S$ fixes $x$ (the other cases follows easily from this case). We suppose without loss of generality that $V \cap [V]$ is not a point of $W$; otherwise, we replace $W$ by $W \setminus \{ V \cap [V] \}$. By Observation \ref{transfer}, $S$ induces a collineation group of the projective space $\bP$  on the infinite point set $W$, which is cofinite in some affine subspace of $\bP \setminus \bP'$ of dimension $\geq 1$. This is clearly not possible. \eop \\

\medskip
\subsection{An interesting example}
\label{intexline}

Let $X$ be a TGQ with translation point $u$ and translation group $T$ (of sufficiently large order to avoid degeneracies). Let $\Gamma$ be a full grid inside $X$, and let $U$ be a line in $\Gamma$. 
We have that
\begin{equation}
[X]\ =\ [X \setminus (\Gamma,X)]\ +\ [(\Gamma,X)]\ =\ [X \setminus (\Gamma,X)]\ +\ [(\Gamma,X) \setminus (L,(\Gamma,X))]\ +\ [(L,(\Gamma,X))]. 
\end{equation}

On the other hand, we also have that
\begin{equation}
[X]\ =\ [X \setminus (L,X)]\ +\ [(L,X)]\ =\ \Big[\Big(X \setminus (L,X)\Big) \setminus \Big((\Gamma,X) \setminus (L,X)\Big)\Big] \ +\ [(\Gamma,X) \setminus (L,X)]\ +\    \ [(L,X)].
\end{equation}

As obviously 

\begin{equation}
\left\{
                \begin{array}{ccc}
[X \setminus (\Gamma,X)] &= &\Big[\Big(X \setminus (L,X)\Big) \setminus \Big((\Gamma,X) \setminus (L,X)\Big)\Big] \\

[(\Gamma,X) \setminus (L,(\Gamma,X))] &= &[(\Gamma,X) \setminus (L,X)]\\ 
\end{array}
\right.
\end{equation}
\medskip
(in the first equality, the representants are evaluated on a purely geometric basis, and in the second equality, the representants come with the same induced permutation group), we have that 
\begin{equation}
\label{eqrememb}
[(L,(\Gamma,X))]\ =\ [(L,X)].
\end{equation}

We interpret equation (\ref{eqrememb}) in this specific case, as the fact that $K_0(\mQ_\ell)$ \ul{``remembers the embedding.''} \\ 

Note that we (oviously) have an embedding of lines 
\begin{equation}
(L,(\Gamma,X))\ \hookrightarrow\ (L,X). 
\end{equation}

Of course, it remains to find examples of triples $(X,\Gamma,L)$ which satisfies the assumptions made initially. We provide one example. 

Suppose $X$ is a nonclassical Kantor-Knuth generalized quadrangle, and let $[\infty]$ be the special line (this is the only line which is fixed by $\Aut(X)$). 
All the lines in $[\infty]^{\perp}$ are regular. We consider a full grid $\Gamma$ which contains one point $v$ of $[\infty]$ (such a grid obviously exists). 
Let $U$ be any line of $\Gamma$ which is not incident with $v$. Put $A := \Aut(X)$. We want to compare ${(A_\Gamma)}_U$ and $A_U$. 

First of all, not that $A_U$ acts $2$-transitively on the points incident with $U$. This is easy to see, since all lines in $\{[\infty],U\}^{\perp}$ are 
axes of symmetry, so the group generated by these symmetries already acts doubly transitively on the lines in $\{[\infty],U\}$, whence on the 
points of $U$. 

Now consider ${(A_\Gamma)}_U = {(A_U)}_\Gamma$. Every element of ${(A_U)}_\Gamma$ fixes the point $v$, so also the point $w := 
\proj_Uv$. It follows that
\begin{equation}
{\Big( A_\Gamma\Big)}_U \ \ne\ A_U.
\end{equation}

So the lines $(L,X)$ and $(L,(\Gamma,X))$ are indeed nonisomorphic, but define the same class in the $K_0(\mQ_\ell)$ which contains $[X]$.

\medskip
\subsection{Rigid lines}

Suppose $X$ is a {\em rigid GQ} --- that is, $\Aut(X) = \{ \id\}$. Then in IDEA 2 all lines in $X$ are isomorphic, while in IDEA 1, no two lines are isomorphic. 
In this paragraph, we assume to work in IDEA 2. 

Let $U$ be a line in $X$.
Now consider the grid $\Gamma := (U,X) \times (U,X)$; as $\Aut(U) := \Aut\Big((U,X)\Big) = \{ \id \}$, it follows that $\Aut(\Gamma)$ cannot contains isomorphisms 
fixing any line in $\Gamma$ with a nontrivial action on that line. It follows easily that this is equivalent to the property that non nontrivial element of $\Aut(\Gamma)$ 
fixes a line of $\Gamma$.

\newpage
\section{Product}
\label{prod}

In the Grothendieck ring $K_0(\mV_k)$, we have the important identity
\begin{equation}
\bL^{\otimes n}\ =\ \bL^n\ =\ [\bA^n(k)]. 
\end{equation}

In the category of incidence geometries, we have several products at our disposal, usually based on a product on some underlying graph. For now, we opt for the {\em Cartesian product}
of the underlying collinearity graph $C(\Gamma)$ of a point-line geometry $\Gamma$. The {\em collinearity graph} $C(\Gamma)$ of $\Gamma$ has as vertices the points of 
$\Gamma$, and $u \sim v \ne$ if $u$ and $v$ are collinear. By definition, a vertex is adjacent with itself.
If $C_1$ and $C_2$ are graphs, the {\em Cartesian product} $C_1 \otimes C_2$ has as vertices the elements of 
the Cartesian product of the underlying vertex sets, and $(u,v) \sim (x,y)$ if 
\begin{equation}
u = x\ \mbox{and}\ v \sim y \ne v,\ \mbox{or}\ u \sim x \ne u\ \mbox{and}\ v = y. 
\end{equation}

\subsection{Product of lines}

Consider the class of an abstract line ($\ell + 1$ points incident with one line): $[L]$. Then
\begin{equation}
\Big[L\Big]\ \times\ \Big[L\Big]\ =\ \Big[ L \otimes L\Big],
\end{equation}
and the latter is the class of a thin quadrangle of order $(\ell,1)$. This is precisely what we want | in some sense it corresponds to the identity 
$[\mathbb{A}^1(k)] \times [\mathbb{A}^1(k)] = [\mathbb{A}^2(k)].$ 
(The Cartesian product of two complete graphs $K_n$ and $K_m$ is 
the rook graph | the line graph of the complete bipartite graph $K_{n,m}$, which is isomorphic to a dual $(n \times m)$-grid.) 

In general, $[L ]^n := [L]^{\otimes n}$ looks like an affine space of dimension $n$ with only $n$ parallel classes of lines. It has 
no triangles as subgeometry.

\subsection{On $A \setminus B$}

Suppose $A$ is a point-line geometry, and let $B$ be a subgeometry. Consider the geometry $A \setminus B$ (and $\USpec(A \setminus B)$). Let $U$ be a line of $A$, 
and suppose $B$ contains at least two points of $U$, while $U$ is incident with at least two points of $A$ which are not in $B$. Then since $B$ and $A \setminus B$ are subgeometries of $A$, the line $U$ is a line of both $B$ and $A \setminus B$. This seems to contradict the expression

\begin{equation}
\label{AB}
A\ =\ B\ \coprod\ \Big( A \setminus B \Big) 
\end{equation}
on the level of lines. It is not, though: as we have seen in much detail, the line $(U,A)$ is different than the line $(U,B)$.
(It also works on the level of points; so we do not have to agree that expressions such as (\ref{AB}) should be interpreted solely on the level of the underlying point sets.)

If $U$ is a line of $A$ with at least two distinct points, and all points of $U$ are also points of $B$ ($U$ is full in $B$), $A \setminus B$ does not contain any point of $A$ incident with $U$, so 
$U$ is not a ``proper  line'' in $A \setminus B$. There are two cases to consider: 
\begin{itemize}
\item[(a)]
if $U^* \subseteq B$ and $U$ is not a line in $B$, the points of $U^*$ are not collinear in $B$, and $A \setminus B$ contains 
the empty line $U^e$;
\item[(b)]
if $U^* \subseteq B$ and $U$ is in $B$ as well, then $U$ is a full line in $B$, and $A \setminus B$ does not contain $U$, nor 
any point of $U^*$. 
\end{itemize}
(Compare this situation  to the setting of $k$-varieties for some field $k$!) Still, we do consider ``empty lines'' | see section \ref{calc} for more details. (Note that if $U^e$ is an empty line, then $C(U^e) = \emptyset$.)

Now let $U$ be a line of $A$ and $B$ (so with at least two distinct points in $B$), such that all points but one point $u$ incident with $U$ are  points of $B$; then $A \setminus B$ contains the point $u$, but the line $U$ only contains one point in $A \setminus B$. Usually, we do not see $U$ as a line in $A \setminus B$. On the other hand, in the context of $\Fun$-geometry, for instance in the theory of M\'{e}rida-Thas of \cite{MMKT} it has been proven quite useful to see the geometry of an incident point-line pair as nontrivial; it is a combinatorial depictation of the absolute affine line $\Spec(\Fun[X])$ (also called the ``absolute flag''). We imagine that the geometry $A \setminus B$ then remembers the direction of the line $U$ which has been removed from $A \setminus B$  through deletion of its points in $B$, by recording it in a ``hairy point'' picture as below. So in some situations, we will allow this possibility (and therefore extend the definition of product of $\mQ_\ell$-varieties below to make sense of this situation). 

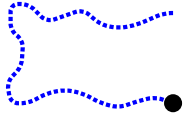
\begin{center}
\begin{tikzpicture}[scale=1.2]
\draw[ultra thick,densely dotted, blue,rounded corners=1mm] (0,1)--(-0.1,1)--(-0.4,0.87)--(-0.6,0.83)--(-0.8,0.87)--(-1,1.04)--(-1.4,0.9)--(-1.6,1.1)--(-1.8,1.1)--(-1.8,0.8)--(-1.66,0.7)--(-1.68,0.4)--(-1.84,0.24)--(-1.8,0)--(-1.6,0)--(-1.4,0.11)--(-1.2,0.15)--(-1,0.11)--(-0.8,0)--(-0.6,-0.05)--(-0.4,0.02)--(-0.2,0.07)--(0,0);
\fill [black] (0,0) circle (0.1);
\end{tikzpicture} 
\end{center}

Depicting hairy points as above (one might add the name of the former line to remember its origin), one easily extends the product definition of above to point-line geometries with absolute flags as follows; add a hairy point picture in the collinearity graphs for each line $U$ which was removed as above from the point $u$, and then add the picture each time 
$u$ occurs in some $(u,v)$ in  the product of the graphs.


\subsection{Properties of the product}

If $\Gamma$ is a point-line geometry, let the {\em line spectrum} of $\Gamma$ be the set of possible line sizes in $\Gamma$. 
In general, we have the following:

\begin{proposition}
\label{prodwell}
\begin{itemize}
\item[{\rm (1)}]
Let $\Gamma_1$ and $\Gamma_2$ be point-line geometries with $\widetilde{(3)}$ and respective line spectra $S_1$ and $S_2$. (We suppose that each line has at least two different points for the sake of convenience.)  Then $\Gamma_1 \otimes \Gamma_2$ has line spectrum $S_1 \cup S_2$, and it satisfies $\widetilde{(3)}$ as well.   
\item[{\rm (2)}]
Let $X$ and $Y$ be elements of $\mQ_\ell$, and let $C$ be a closed subset of $X$. Then $Y \otimes C$ is closed in $Y \otimes X$, and 
\begin{equation}
Y \otimes \Big(X \setminus C\Big) = (Y \otimes X) \setminus (Y \otimes C).
\end{equation}
\end{itemize}
\end{proposition}

{\em Proof}.\quad
(1)\quad
Let $\mG_1$ and $\mG_2$ be the collinearity graphs of $\Gamma_1$ and $\Gamma_2$. Then the lines of $\Gamma_1$ and $\Gamma_2$ correspond to the cliques of $\mG_1$ and $\mG_2$ (and they all have size at least $2$). We can see $\mG_1 \otimes \mG_2$ as a number of copies of $\mG_1$, each for each vertex of $\mG_2$, and obviously any clique 
of $\mG_1 \otimes \mG_2$ either lies in (a copy of) $\mG_1$ or $\mG_2$. So the line spectrum of $\Gamma_1 \otimes \Gamma_2$ is $S_1 \cup S_2$. Also, if $\Gamma_1 \otimes \Gamma_2$ would have a triangle, then the vertices corresponding to its points are contained in a clique of $\mG_1 \otimes \mG_2$, which would be a clique in either $\mG_1$ or $\mG_2$, and this is an obvious contradiction. \\

(2)\quad
First note that $Y \otimes C$ is full in $Y \otimes X$. 
Then, as $X \setminus C$ is the induced subgeometry  (in $X$) on the point set of $X \setminus C$, and as a similar remark holds for $(Y \otimes X) \setminus (Y \otimes C)$, it follows that 
$Y \otimes (X \setminus C)$ and $(Y \otimes X) \setminus (Y \otimes C)$ have the same vertex set (as graphs); obviously, they also have the same adjacencies, so (2) holds. (For expressions which involve empty lines, we refer to section \ref{calc}.)  \\

The proof is finished. \eop \\

By (2), it follows that under the agreement of our nomenclature, 

\begin{equation}
Y \otimes X = \Big((Y \otimes X) \setminus (Y \otimes C)\Big)\ \coprod\ (Y \otimes C),
\end{equation}
so that we have:

\begin{proposition}
\label{propring}
The subgroup $N$ of $K_0(\mQ_\ell)$ freely generated by all relations of the form 
\[ \Big[ X \Big] - \Big[ X \setminus C \Big] - \Big[ C \Big] \]  
is an ideal under our choice of multiplication.  \eop  
\end{proposition}

\begin{proposition}
Let $\mC$ and $\widetilde{\mC}$ be geometries that satisfy $\widetilde{(3)}$, with each line having at least $2$ different points. 
If $\mC \otimes \widetilde{\mC}$ satisfies $(3)$, then either $\mC$ and $\widetilde{\mC}$ are point sets, or $\mC$ and $\widetilde{\mC}$ are lines. 
\end{proposition}

{\em Proof}.\quad
First note that if $\mC$ or $\widetilde{\mC}$ do not satisfy ${(3)}$, $\mC \otimes \widetilde{\mC}$ also does not, so we may assume that 
both $\mC$ and $\widetilde{\mC}$ have ${(3)}$.

If $\mC$ is a point set, then $\widetilde{\mC}$ cannot contain lines, since $\mC \otimes \widetilde{\mC}$ is a disconnected set of copies of $\widetilde{\mC}$. 
Now suppose both $\mC$ and $\widetilde{\mC}$ have lines, $L$ and $\widetilde{L}$. Suppose $\mC$ contains a point $c$ which is not incident with $L$. Since 
(3) is satisfied for $\mC$, there is a unique line $M$ such that $M \I c$ and $M \sim L$. Let $\mC'$ be the subgeometry defined by $c, M, M \cap L, L$; then locally 
at $\mC' \otimes \widetilde{\mC}$, $\mC \otimes \widetilde{\mC}$ looks like a union of two grids $M \otimes \widetilde{L}$ and $L \otimes \widetilde{L}$ meeting in a 
copy of $\widetilde{L}$. Obviously (3) cannot be satisfied. It follows that all points of $\mC$ and $\widetilde{\mC}$ are incident with $L$ and $\widetilde{L}$. \eop \\

\begin{corollary}
If $\mC$ is prime in $X$ and $\widetilde{\mC}$ is prime in $\widetilde{X}$, then $\mC \otimes \widetilde{\mC}$ is only prime in $X \otimes \widetilde{X}$ if 
both $\mC$ and $\widetilde{\mC}$ are point sets, or lines.  \eop 
\end{corollary}

Suppose $X$ and $\widetilde{X}$ are point-line geometries with $\widetilde{(3)}$. Then, in general, we have 
\begin{equation}
\cl(X \otimes \widetilde{X})\ \ne\ \cl(X) \otimes \cl(\widetilde{X}), 
\end{equation}
while obviously $\cl(X \otimes \widetilde{X}) \supseteq \cl(X) \otimes \cl(\widetilde{X})$. 
Here, $\cl(Y)$ denotes the set of closed sets which are contained in $Y$, and 
\[   \cl(X) \otimes \cl(\widetilde{X}) \ =\ \{ C_1 \otimes C_2 \ \vert\ C_1\ \mbox{closed in }\ X,\ C_2\ \mbox{closed in}\ \widetilde{X}  \}.    \] 

In $X \otimes \widetilde{X}$, a perp geometry with at least one line, never is isomorphic to a product (this compares to the product of topological spaces, which usually
also is richer than the product set of the open sets).

\subsubsection{Fiber products}

Let $(\mV,\mE)$ be a graph. 
Define for $i = 1,2$ a map
\begin{equation}
\psi_i:  \mV \times \mV \mapsto \{0,1,2\}:\ (u,u') \mapsto \psi_i(u,u'),
\end{equation}
such that $\psi_1(u,u') = 1$ if $u = u'$, $\psi_1(u,u') = 2$ if $u \sim u' \ne u$, and $\psi_1(u,u') = 0$ otherwise; 
define $\psi_2(u,u') = 1$ if $u \sim u' \ne u$, $\psi_2(u,u') = 2$ if $u = u$, and $\psi_2(u,u') = 0$ otherwise. \\

\tikzset{->-/.style={decoration={
  markings,
  mark=at position 1 with {\arrow{>}}},postaction={decorate}}}

Below, we see graph morphisms between graphs $\Gamma = (V,E)$ and $\Gamma' = (V',E')$ as maps $\beta: V \times V \mapsto V' \times V'$ which preserve adjacency; 
we assume that every vertex is adjacent with itself (since in an incidence geometry we suppose that every point is collinear with itself). In fact, we denote 
the corresponding graph morphism also as $\beta$ (which initially is a map between vertex sets), so that the map between $V \times V$ and $V' \times V'$ is in fact $\beta \times \beta$.  

Fix the graphs $\Gamma$ and $\Gamma'$, and consider the category $\mC$ in which the objects are commutative diagrams



\begin{center}
\begin{tikzpicture}
  \node (s) {$V'' \times V''$};
  \node (xy) [below=2 of s] {$\{0,1,2\}$};
  \node (x) [left=of xy] {$V \times V$};
  \node (y) [right=of xy] {$V' \times V'$};
  \draw[->-] [ultra thick] (s) to node [sloped, above] {$\varphi_2$} (y);
  \draw[<-] [ultra thick] (x) to node [sloped, above] {$\varphi_1$} (s);
  \draw[<-] [ultra thick] (xy) to node [below] {$\psi_1$} (x);
  \draw[<-] [ultra thick] (xy) to node [below] {$\psi_2$} (y);
\end{tikzpicture}
\end{center}
where $\Gamma'' = (V'',E'')$ is a graph, and the $\varphi_i$, $i = 1,2$, are graph morphisms. The morphisms in $\mC$ are natural. 

Then with $\Gamma \otimes \Gamma'$ the fiber product represented by $\Big(V \times V'\Big) \times \Big( V \times V' \Big)$ as above, with $\varphi_1 = \pi_1: 
\Big(V \times V'\Big) \times \Big( V \times V' \Big) \mapsto V \times V: ((a,b),(c,d)) \mapsto (a,c)$ and $\varphi_2 = \pi_2: 
\Big(V \times V'\Big) \times \Big( V \times V' \Big) \mapsto V \times V: ((a,b),(c,d)) \mapsto (b,d)$, 
it is easy to see 
that the fiber product is a terminal object in $\mC$ (and so is a categorical fiber-like product). First note that the $\pi_i$ are indeed graph morphisms | consider, e.g., 
$\pi_1$. If $(a,b) \sim (c,d)$ in $\Gamma \otimes \Gamma'$, we either have $a = c$ and $b \sim d \ne b$, or $a \sim c \ne a$ and $b = d$, or $a = c$ and $b = d$. In the first and last case, $((a,b),(c,d))$ is mapped on $(a,c) = (a,a)$; in the second case, the image is $(a,c)$ with $a \sim c \ne a$. The case $\pi_2$ is of course similar. \\

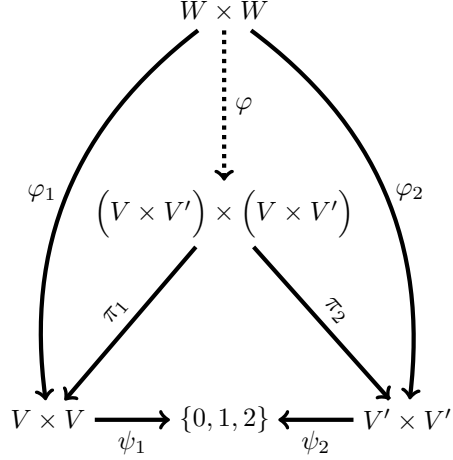
\begin{center}
\begin{tikzpicture}
  \node (s) {$\Big( V \times V' \Big) \times \Big( V \times V' \Big)$};
  \node (z) [above=2 of s] {$W \times W$};
  \node (xy) [below=2 of s] {$\{0,1,2\}$};
  \node (x) [left=of xy] {$V \times V$};
  \node (y) [right=of xy] {$V' \times V'$};
  \draw[->-] [ultra thick] (s) to node [sloped, above] {$\pi_2$} (y);
  \draw[<-] [ultra thick] (x) to node [sloped, above] {$\pi_1$} (s);
  \draw[<-] [ultra thick] (xy) to node [below] {$\psi_1$} (x);
  \draw[<-] [ultra thick] (xy) to node [below] {$\psi_2$} (y);
  \draw[->] [ultra thick, dotted] (z) to node [right] {$\varphi$} (s);
 \draw[->,ultra thick] (z) edge [bend right] node [left]  {$\varphi_1$} (x);
\draw[->,ultra thick] (z) edge [bend left] node [right] {$\varphi_2$} (y);
\end{tikzpicture}
\end{center}

\medskip
Now let $o$ be any object in $\mC$ defined by the graph $\Gamma^* = (V^*,E^*)$, and graph morphisms $\varphi_1: V^* \times V^* \mapsto V \times V$, $\varphi_2: V^* \times V^* \mapsto V' \times V'$. Let $(w,w') \in V^* \times V^*$ such that $w \sim w'$, and such that
\begin{equation}
\psi_1 \circ \varphi_1(w) = \psi_2 \circ \varphi_2(w').
\end{equation}

Then $\varphi(w,w')$ must be equal to $((\varphi_1(w),\varphi_2(w)),(\varphi_1(w'),\varphi_2(w'))$. If $w = w'$, $\psi_1(\varphi_1(w),\varphi_1(w')) = 0 = \psi_2(\varphi_2(w),\varphi_2(w'))$, so that 
$\varphi_2(w) \sim \varphi_2(w') \ne \varphi_2(w)$. So $\varphi(w) \sim \varphi(w')$ in $\Gamma \otimes \Gamma'$. Now let $w \ne w'$. 
If $\varphi_1(w) = \varphi_1(w')$, a similar argument leads to $\varphi(w) \sim \varphi(w')$ in $\Gamma \otimes \Gamma'$. 
Suppose this is not the case. As $\varphi_1(w) \sim \varphi_1(w') \ne \varphi_1(w)$, we have that $\psi_1(\varphi_1(w),\varphi_1(w')) = 1 = \psi_2(\varphi_2(w),\varphi_2(w'))$, so that 
$\varphi_2(w) = \varphi_2(w')$. It follows again that $\varphi(w) \sim \varphi(w')$ in $\Gamma \otimes \Gamma'$.

\medskip
\subsection{Neutral element}

The class of a point, $[ e ]$, obviously is a neutral element for our product, unique up to unique isomorphism.

\medskip
\subsection{Automorphisms of products}

Let 
\begin{equation}
X = Y \otimes Z, 
\end{equation}
with $X, Y, Z$ elements in $\mQ_\ell$. Then obviously $\Aut(Y) \times \Aut(Z)$ is a subgroup of $\Aut(X)$. We say that an element $U \in \mQ_\ell$ 
is {\em primal} (as opposed to being {\em prime} in the Zariski context) if it is not the product of elements which are different than points. We use the same nomenclature for the underlying collinearity graphs. 

\begin{theorem}[Sabidussi \cite{sabi}]
\label{Sabu}
If $\Gamma_1$ and $\Gamma_2$ are graphs, and both are prime, then $\Aut(\Gamma_1 \otimes \Gamma_2) \cong \Aut(\Gamma_1) \times \Aut(\Gamma_2)$ if $\Gamma_1$ is 
not isomorphic to $\Gamma_2$. If they are isomorphic, we have that  $\Aut(\Gamma_1 \otimes \Gamma_2) \cong \Big(\Aut(\Gamma_1) \times \Aut(\Gamma_2)\Big)\rtimes \mathrm{Sym}(2)$.
\end{theorem}

{\em Problem}. Let $X$ be a nonclassical $\hT_2(\mO)$ of order $(\ell,\ell)$, and let $W = (W,X)$ be a line which is not incident with the translation point. What is 
$\Aut\Big((W,X) \otimes (W,X)\Big)$? According to Theorem \ref{Sabu}, we would have that it is isomorphic to  $\Aut\Big((W,X)\Big) \times \Aut\Big((W,X)\Big)$ ``up to a twist,'' 
and that implies that the outcome is in general different than $\Aut(W \otimes W)$, where $W = (W,W)$ is seen as an abstract line (without the embedding in $X$).  (If $\Aut\Big((W,X) \otimes (W,X)\Big)$ would be 
bigger, then a stabilizer of a copy of $(W,X)$ would be too big according to our formalism.)\\

We say that $(W,X)$ {\em embeds} in $(W',Y)$ if there is an injective morphism 
\begin{equation}
\epsilon: (W,\Aut(X)_W/\mbox{kernel})\ \hookrightarrow\ (W',\Aut(Y)_{W'}/\mbox{kernel}) 
\end{equation}
in the category of permutation representations of degree $\vert W \vert = \vert W' \vert$ (= the category of lines of size $\vert W \vert = \vert W'\vert$). \\

The following diagram commutes:

\begin{center}
\item
\begin{tikzpicture}[>=angle 90,scale=2.2,text height=1.5ex, text depth=0.25ex]
\node (a0) at (0,3) {$(W,X)$\ \ };
\node (a2) [right=of a0] {};
\node (a1) [right=of a2] {\ \ $(W,X) \otimes (W,X)$};

\node (b0) [below=of a0] {$(W,W)$\ \ };
\node (b1) [below=of a1] {\ \ $(W,W) \otimes (W,W)$};

\draw[->,font=\scriptsize,thick]
(a0) edge node[left] {$\epsilon$} (b0)
(a1) edge node[right] {$\epsilon \times \epsilon$} (b1)
(a0) edge node[auto] {$\iota_1$} (a1)
(b0) edge node[below] {$\iota_1'$} (b1);
\end{tikzpicture}
\end{center}

All the maps are embeddings; $\epsilon$ exists because $\Aut(W,W) \cong \mathrm{Sym}(W)$, and $\iota_1, \iota_1'$ are 
``embeddings into the first component." On the level of automorphism groups, we obtain, with the same notation for the induced morphisms on 
the automorphism groups: 

 \begin{center}
\item
\begin{tikzpicture}[>=angle 90,scale=2.2,text height=1.5ex, text depth=0.25ex]
\node (a0) at (0,3) {$\Aut(X)_W/\mbox{kernel}$\ \ };
\node (a2) [right=of a0] {};
\node (a1) [right=of a2] {\ \ $\Big(\Aut(X)_W/\mbox{kernel}\times\Aut(X)_W/\mbox{kernel} \Big)\rtimes \mathrm{Sym}(2)$};

\node (b0) [below=of a0] {$\mathrm{Sym}(W)$\ \ };
\node (b1) [below=of a1] {\ \ $\Big(\mathrm{Sym}(W) \times \mathrm{Sym}(W)\Big)\rtimes \mathrm{Sym}(2)$};

\draw[->,font=\scriptsize,thick]
(a0) edge node[left] {$\epsilon$} (b0)
(a1) edge node[right] {$\epsilon \times \epsilon$} (b1)
(a0) edge node[auto] {$\iota_1$} (a1)
(b0) edge node[below] {$\iota_1'$} (b1);
\end{tikzpicture}
\end{center}

The notion of embeddings can be easily and naturally generalized to all the geometries that we consider, and we will use the general notion 
throughout.

\newpage
\medskip
\section{The ring $K_0(\mQ_\ell)$}
\label{grogro}

We are finally ready to define the Grothendieck ring.

Suppose $\mQ_\ell$ is defined as in section \ref{ql}. The various notions of isomorphisms were described in section \ref{isom} (and the choice of isomorphisms heavily impacts the structure of the ring as we will see below).
Then $K_0(\mQ_\ell)$ is defined as the free abelian group generated by the isomorphism classes $[Q]$ of objects in $\mQ_\ell$, 
modulo the scissor relations
\begin{equation}
\Big[Q \Big]\ \ =\ \ \Big[Q \setminus C\Big]\ \ +\ \ \Big[C\Big], 
\end{equation}
with $C$ a closed set in $Q$ (in the old or new Zariski setting). Clearly, the class of the  ``empty generalized quadrangle'' {\Large $\varnothing$} is the unique neutral element for taking sums.\\ 

{\em Example}: we have, in $K_0(\mQ_q)$ with $q$ a prime power: 
\begin{equation}
\Big [\mQ(4,q) \Big]\ \ =\ \ \Big[\mQ(4,q) \setminus \mQ(3,q)\Big]\ \ +\ \ \Big[\mQ(3,q)\Big]. \\
\end{equation}

\medskip
Products are defined as in section \ref{prod} (and the neutral element for this product is the class $[e]$ of a point).

Then by section \ref{prod}, we can conclude the following result.

\begin{proposition}[The ring $K_0(\mQ_\ell)$]
We have that $K_0(\mQ_\ell)$ is a commutative unital ring.  
\end{proposition}

{\em Proof}.\quad
Immediately by Proposition \ref{propring}. \eop \\



\medskip
\subsection{Impact of isomorphisms}

By the scissor relations, classes in the Grothendieck group become much larger than initially defined: representants of the same 
class need not be isomorphic ``at the end.''  The choice of the notion of isomorphism is crucial in this context | below, we will give some 
examples which clearly show the pivotal role of these notions. (Changing the notion of isomorphism changes the ring entirely.)

Let $X$ be a generalized quadrangle with a regular line $U$, and suppose $V$ is not a regular line. We suppose furthermore that 
$\Aut(X)_U/\mbox{kernel} \cong \Aut(X)_V/\mbox{kernel}$. Then by the scissor relations, we have that 
\begin{equation}
[X \setminus U]\ =\ [X \setminus V]. 
\end{equation}

Clearly the geometries $X \setminus U$ and $X \setminus V$ cannot be isomorphic. If we would use another notion of isomorphism, such as 
the one in IDEA 1, they {\em would} be isomorphic. 

\begin{remark}{\rm 
Exactly the same phenomenon happens in the Grothendieck of varieties $K_0(\mV_k)$: the equality 
\begin{equation}
\label{eqxy}
[X]\ =\ [Y]
\end{equation}
does not necessarily imply that $X$ and $Y$ are isomorphic; one important open question asks whether  (\ref{eqxy}) implies 
that $X$ and $Y$ are birationally equivalent.} 
\end{remark}

Now let $X$ be a generalized quadrangle of order $(u,u)$, and let $\Gamma$ and $\Gamma'$ be full grids in $X$. We suppose that 
$(\Gamma,X) \cong (\Gamma',X')$. The scissor relations imply that $[X \setminus \Gamma] = [X \setminus \Gamma']$. As $X \setminus \Gamma$ and 
$X \setminus \Gamma'$ are affine quadrangles, the property that they would be isomorphic (by an isomorphism $\alpha$), implies that 
there is an isomorphism
\begin{equation}
\overline{\alpha}: X \mapsto X
\end{equation}
which maps $\Gamma$ to $\Gamma'$, and which induces $\alpha$. So this would mean that $\Gamma$ and $\Gamma'$ are isomorphic in a stronger sense 
than the one we want to use.

\medskip
\section{Krull dimension}
\label{Krull}

For a generalized quadrangle $\Gamma$, or any other of the connected variety-like elements in $\mQ_\ell$, we define the {\em Krull dimension} as 
the supremum of the lengths $d$ of maximal chains 
\begin{equation}
\fp_0 \subset \fp_1 \subset \cdots \subset \fp_d = \Gamma,
\end{equation}
where each $\fp_i$ is a (nonempty) point in the Zariski topology, and inclusions are with respect to the inclusion of incidence geometries. \\

\subsection{Example}

Let $\Gamma \cong \mQ(5,q)$, $q$ a prime power. Then a maximal chain is of the form
\begin{align*}
\fp_0 = \mbox{point} &\subset \fp_1 = \mbox{full line} \subset \fp_2 = \mbox{two distinct intersecting full lines}\\ 
&\subset \fp_3 =   \mbox{full grid} \subset \fp_4 = \mbox{$\mQ(4,q)$-subGQ} \subset \fp_5 = \Gamma. \\
\end{align*}

So $\mQ(5,q)$ and $\mQ(4,q)$ have dimension $5$ and $4$. It is also clear that any grid has dimension $3$. 

\subsection{Example}

Let $\Gamma \cong \mH(3,q^2)$, $q$ a prime power. Then a maximal chain is of the form
\begin{align*}
\fp_0 = \mbox{point} &\subset \fp_1 = \mbox{full line} \subset \fp_2 = \mbox{full ideal perp geometry} \subset \fp_3 = \Gamma. \\
\end{align*}

So the dimension of $\mH(3,q^2)$ is $3$, while on the other hand $\mH(3,q^2)$ is isomorphic to the point-line dual of $\mQ(5,q)$. \\

\subsection{Example}

Let $q$ be a prime power and consider $\Gamma \cong \hW(q)$. If $q$ is even, we know that $\hW(q) \cong \mQ(4,q)$, so its dimension is $4$. 
On the other hand, if $q$ is odd, a maximal chain is of the form 
\begin{align*}
\fp_0 = \mbox{point} &\subset \fp_1 = \mbox{full line} \subset \fp_2 = \mbox{full ideal perp geometry} \subset \fp_3 = \Gamma, 
\end{align*}
and so its dimension is $3$!\\

\begin{theorem}
Let $\Gamma$ be a thick generalized quadrangle. Then its dimension is at least $3$. 
\end{theorem}

{\em Proof}.\quad
We always have chains of the form 
\begin{align*}
\fp_0 = \mbox{point} &\subset \fp_1 = \mbox{full line} \subset \fp_2 = \mbox{full ideal perp geometry} \subset \fp_3 = \Gamma.
\end{align*}
\eop \\

Needless to say, there are many degenerate examples in lower dimensions. \\

\subsection{Affine discussion}

In the case of affine examples, we can also use the affine list in section \ref{affine} to introduce a notion of Krull dimension. 

\subsection*{Example}

We construct an affine quadrangle $\underline{\Gamma}$ be taking away a $\mQ(4,q)$-subGQ $\Delta$ in $\Gamma \cong \mQ(5,q)$. (Note that $\underline{\Gamma}$ defines an open set in $\Gamma$.) One example of a chain of maximal length is the following:
\begin{align*}
\fp_0 = \mbox{point} &\subset \fp_1 = \mbox{full line in $\underline{\Gamma}$} \subset \fp_2 = \mbox{two distinct intersecting full lines in $\underline{\Gamma}$}\\ 
&\subset \fp_3 =   \mbox{full grid in $\underline{\Gamma}$} \subset \fp_4 = \underline{\Gamma}' \subset \fp_5 = \underline{\Gamma}. \\
\end{align*}

Here, $\underline{\Gamma}'$ is the intersection of a $\mQ(4,q)$-subGQ in $\Gamma$ different from $\Delta$ with $\underline{\Gamma}$; note that $\underline{\Gamma}$ can intersect in different ways with $\Delta$ so that ``nonisomorphic chains'' can be obtained (of the same maximal length). Here, we assume that $\Delta$ and the $\mQ(4,q)$-subGQ intersect in a full ideal perp geometry of $\Delta$. 
Note also that every full line in $\Gamma \cong \mQ(5,q)$ contains precisely one point of $\Delta$. \\

We conclude that the dimension of $\underline{\Gamma}$ is $5$. \\

\subsection{Larger dimensions}

We now describe an example of a thick generalized quadrangle of non-finite dimension. The details are contained in \cite{KTAb}. In that latter paper, a generalized quadrangle is described of order $(\vert k \vert, \vert k \vert)$, with $k$ any infinite field, which contains an infinite chain of full subGQs
\[  \fp_0\ \subset\ \fp_1\ \subset\ \fp_2\ \subset\ \cdots                     \]  

The construction starts from an infinite GQ $\Omega$ with a linear representation in $\bP^n(k)$, with $n = \vert \mathbb{N} \vert$ and $k$ an infinite field. In $\bP^{n - 1}(k)$, the point set which defines $\Omega$ is a set which meets every line of $\bP^{n - 1}(k)$ in precisely two different points. This construction could probably also be used (with finite $n$) to construct generalized quadrangles of finite dimension $n \geq 6$. 

The interest of such constructions lies in the following theorem:

\begin{theorem}[Dimension Theorem]
If $\Gamma$ is a thick generalized quadrangle of dimension at least $6$, then $\Gamma$ is not finite. 
\end{theorem}

{\em Proof}.\quad
By \cite[section 2.2.2]{PTsec}, the longest possible chains we can theoretically have in a thick finite GQ $\Gamma$ of order $(s,t)$ contains full subgrids, and are of the form
\begin{align*}
\fp_0 = \mbox{point} &\subset \fp_1 = \mbox{full line in ${\Gamma}$} \subset \fp_2 = \mbox{two distinct intersecting full lines in ${\Gamma}$}\\ 
&\subset \fp_3 =   \mbox{full grid in ${\Gamma}$} \subset \fp_4 = {\Gamma}' \subset \fp_5 = {\Gamma},
\end{align*}
where $\Gamma'$ is a subGQ of order $(s,s)$. 
\eop \\

Note that all $\mQ(5,q)$-GQs and all Kantor-Knuth GQs have dimension $5$.

\medskip
\section{Infinite $\ell$, and Lefschetz motives}
\label{Lef}

Suppose $\Gamma$ is an abstract grid with $\ell + 1$ points per line. Let $\ell$ be finite. Define a geometry $P(N)$ to have two intersecting lines, with $N$ points per lines ($N$ a positive integer). 
Then we have the following identity in $K_0(\mQ_\ell)$:

\begin{equation}
\Big[ \Gamma = (\ell + 1) \times (\ell + 1)-\mbox{grid} \Big]\ =\    \Big[ P(\ell + 1) \Big]\ + \       \Big[ (\ell \times \ell)-\mbox{grid} \Big]
\end{equation}
so that inductively we get that 
\begin{equation}
\Big[ \Gamma = (\ell + 1) \times (\ell + 1)-\mbox{grid} \Big]\ =\    \Big[ P(\ell + 1) \Big]\ + \       \Big[ P(\ell) \Big]\ + \ \ldots\ +\ \Big[ P(1) \Big]. 
\end{equation}

Now let $\ell$ be not finite, and let $m$ be a positive integer. Then 
\begin{equation}
\Big[ \Gamma = (\ell + 1) \times (\ell + 1)-\mbox{grid} \Big]\ =\    \Big[ P(\ell + 1) \Big]\ + \       \Big[ P(\ell) \Big]\ + \ \ldots\ +\ \Big[ P(\ell + 2 - m) \Big]\ +\ \Big[ (\ell + 1 - m)\times(\ell + 1 - m)-\mbox{grid} \Big]. 
\end{equation}

If we argue that $P(\ell + 1) \cong P(\ell) \cong \cdots \cong P(\ell + 1 - m)$, and likewise that $(\ell + a)\times(\ell + a)$-grids and $(\ell - b)\times(\ell - b)$-grids  are isomorphic for finite $a$ and $b$, then we 
obtain 
\begin{equation}
0 = m\Big[ P(\ell + 1) \Big],  
\end{equation}
so that either we have zero divisors, or 
$\Big[ P(\ell + 1) \Big] = 0$, or $m = 0$ for any positive integer $m$. In the latter case, we would wind up in a set-up for working in characteristic $1$, but that is not what we want, 
obviously. 

The problem lies deeper: for suppose that $\mU$ is a hyperbolic quadric in $\mathbb{P}^3(k)$, with $k$ an infinite field. If we remove the closed set $P(\vert k\vert + 1)$ from $\mU$ in 
$K_0(\mV_k)$ (with ``fully embedded lines''), then $\mU \setminus P(\vert k\vert + 1)$ is not isomorphic to $\mU$ as a projective variety, and the reason essentially is that if we consider $\mathbb{P}^1(k)$ and 
remove (example given) a point, we obtain $\Spec(k[X])$ and not again a projective line over $k$. In general, if we remove a finite number of points from a projective line $\mathbb{P}^1(k)$, the obtained object is not isomorphic to $\mathbb{P}^1(k)$. If we were to remove an {\em infinite} number of points, the situation is very different of course | see the examples at the end of section \ref{geomgr}. \\

If we go back to $K_0(\mQ_\ell)$, we see that we need to mend this problem in one way or the other (for finite $\ell$ no problem occurs). Essentially four ``questions'' arise:
\begin{itemize}
\item[(F1)]
we want to see a difference between lines which come from generalized quadrangles in $\mQ_\ell$, and lines which are pieces of such lines;
\item[(F2)]
the dual problem of (F1);
\item[(F3)]
we want to see a difference between grids which are full grids in generalized quadrangles in $\mQ_\ell$, and grids which are pieces of such grids, and 
abstract grids which are products of lines (be it of full lines, or pieces of lines);
\item[(F4)]
formulation of (F3) for higher products. 
\end{itemize}

Since $K_0(\mQ_\ell)$ is generated by isomorphism classes of elements of $\mQ_\ell$, each line is ``of geometric origin": it is a piece of a full line 
in a generalized quadrangle. This is not the case for grids: either they are of ``geometric origin,'' or they are a product of lines, but then the latter are of geometric origin. For that matter we introduce {\em parents}: if, for example, a grid is of geometric origin, then we know it comes with some  embedding in a thick generalized quadrangle $\Gamma$; we call $\Gamma$ a {\em parent} of the grid. Objects in $K_0(\mQ_\ell)$ can obviously come with many parents. 

In any case, these simple facts allow us to distinguish between the various types of lines and grids. And that is exactly how we do it. For example, we agree that 
a fully embedded grid $\Gamma$ in a thick generalized quadrangle $X$  with $\ell + 1$ points per lines {\em can never be isomorphic to an abstract grid} $\Gamma_{ab}$ with $\ell + 1$ points per line, even 
if they would happen to be isomorphic in the theory of IDEA 2:
\begin{equation}
\Big(\Aut(X)_\Gamma,\Gamma\Big)\ \cong\ \Big((\mathrm{Sym}(\ell + 1) \times \mathrm{Sym}(\ell + 1))\rtimes C_2,\Gamma_{ab}\Big).
\end{equation}
    
More details can be found in section \ref{laws}.\\

\subsection{Geometric grids in geometric grids}
\label{geomgr}

The following result  considers a grid which is fully embedded in a generalized quadrangle, and a subgrid which is fully embedded in a subquadrangle. In some circumstances, 
such subgrids cannot exist, and such results are of obvious interest for the discussion in this section. 

\begin{theorem}
\label{fullgr}
Let $\Gamma$ be a full grid in a generalized quadrangle $X$ of order $(\ell,\omega)$; then for finite $m$, the grid $\widehat{\Gamma}$ which arises by removing $m$ lines (and their points) from each 
of the reguli, can not be a full grid in an ideal subquadrangle of $X$.  
\end{theorem}

{\em Proof}.\quad
Suppose by way of contradiction that $\widehat{X}$ is a subGQ of $X$ which contains $\widehat{\Gamma}$ as a full grid. Let $U$ be a line which was removed to 
construct $\widehat{\Gamma}$. By projecting each point of $\widehat{X}$ on $U$, we obtain a set of lines $S_U$ in $\widehat{X}$ which forms a spread in the latter quadrangle. 
Obviously, the lines of $S_U$ are precisely the 
lines of $\widehat{\Gamma}$ which are contained in one of its reguli. But that means that the line of the opposite regulus meet all the lines of $S_U$, which is an obvious 
contradiction. \eop \\  

Note that the statement remains true if we allow $m$ to be not finite. \\

Obviously the ideal property is crucial in the proof of the latter theorem. In the next theorem, we will not ask idealness. We will also generalize the numeric assumption. 

\begin{theorem}
\label{subgr}
Let $\Gamma$ be a full grid in a generalized quadrangle $X$ of order $(\ell,\omega)$, and consider a cardinal number $m$ for which $\ell + 1 - m > m$; then the grid $\widehat{\Gamma}$ which arises by removing $m$ lines (and their points) from each of the reguli, can not be a full grid in a subquadrangle of $X$.  
\end{theorem}

{\em Proof}.\quad
Let $\widehat{X}$ be a hypothetical subquadrangle of $X$ which contains $\widehat{\Gamma}$ as a full grid. 
Suppose that $U$ is a line which is not contained in $\Gamma$, and let $S_U$ be the set of points which are incident with $U$ but not contained in $\widehat{X}$. Let $V$ be a line of $\Gamma \setminus \widehat{\Gamma}$, and suppose that $U \ne V$. Let $\mU$ be the set of lines of $\Gamma'$ which meet $V$. By projecting each line of $\mU$ on $U$, we immediately see that the intersection points of such lines with $V$, must be collinear with points in $S_U$. As the set of these intersection points has size $\vert \mU \vert = \ell + 1 - m$, and as $\ell + 1 - m > m = \vert S_U \vert$, we conclude that triangles appear. This ends the proof of the statement.  
\eop \\

If we remove the condition $\ell + 1 - m > m$, the statement of the last theorem is no longer true. For example, consider a thick orthogonal quadrangle $\mQ$  defined in a projective space 
$\mathbb{P}^n(k)$, with $k$ a field, by a quadratic form. By field reduction (to a field $f$ of the same size as $k$), we can obtain a thick subquadrangle which is also orthogonal, say $\widehat{\mQ}$. As 
each line in an orthogonal quadrangle is regular, any full grid $\widehat{\Gamma}$ in $\widehat{\mQ}$ lies in a full grid $\Gamma$ of $\mQ$, and with the notation of Theorem \ref{subgr}, we have 
that $\ell = \vert k \vert$, and $\ell + 1 - m = \vert f \vert + 1$, so that $m = \vert k \vert = \vert f \vert$. There are many such examples.  \\

Note that we can take $f$ isomorphic to $k$ (for example, let $k := \mathbb{Q}(x_1,x_2,\ldots)$ with the variables $x_i$ indexed by positive integers, and let $f := \mathbb{Q}(x_2,x_3,\ldots)$). In this case, 
we have that 
\begin{equation}
\Gamma \cong \Gamma',
\end{equation}
so even if $\Gamma'$ is a piece of $\Gamma$, it can still be isomorphic to $\Gamma$ (and the same can be said for the projective $k$-lines which are the lines of $\Gamma$).

\medskip
\subsection{Further additive laws}
\label{laws}

Motivated by the discussion in the previous section, we impose the following extra rules on isomorphisms between degenerate elements of $\mQ_\ell$, in case $\ell$ is not finite. In the finite case, 
the rules are trivial. \\

\begin{itemize}
\item[\framebox{L1}]
A line $U$ minus a finite nonzero number of points is not isomorphic to $U$ itself. \\
\item[\framebox{L2}]
A grid $\Gamma$ with a finite nonzero number of full lines removed is not isomorphic to $\Gamma$ itself. \\
\item[\framebox{L3}]
{\em General form}: an element $X$ of $\mQ_\ell$ with a finite nonzero number of points and/or full lines removed is not isomorphic to $X$ itself. 
\end{itemize}





\medskip
\subsection{Calculus of the empty line}
\label{calc}

Suppose $L$ is a full (``proper'') line in some element $X$ in $\mQ_\ell$, and define $L^e := L \setminus L^*$. As $L^*$ is a full subgeometry of $L$ (which is considered as a point set without further structure), we have 

\begin{equation}
\Big[L \Big]\ := \ \Big[ L^e \Big]\ + \ \Big[ L^* \Big].  \\
\end{equation}

Let $\mS$ be any other object in $\mQ_\ell$ which contains points. Then the collinearity graph of $L^e$ is empty, so the product $L^e \otimes S$ seems to be ill defined. 
On the other hand, we have that 

\begin{equation}
\Big[ L^e \otimes S \Big]\ =\ \Big[(L \setminus L^*) \otimes S \Big]\ = \ \Big[ L \otimes S \Big]\ -\ \Big[ L^* \otimes S \Big]. \\
\end{equation}

If $S^*$ is the point set of $S$, and if we assume for the sake of convenience that $S^*$ is finite, then obviously

\begin{equation}
\label{law1}
 \Big[ L \otimes S \Big]\ -\ \Big[ L^* \otimes S \Big]\ =\    \Big[ (L \otimes S) \setminus (L^* \otimes S)\Big]     \ =\ \Big[\coprod_{s \in S^*}L^e\Big]\ =\ \vert S^* \vert \cdot \Big[ L^e \Big]. \\
\end{equation}

\medskip
By Remark \ref{rem2.2} below, $L^e \cong N^e$ for all lines $L, N$ in the setting of IDEA 2. In fact, even through a unique isomorphism. So 
\begin{equation}
\Big [N^e \otimes S \Big]\ =\ \vert S^* \vert \cdot \Big[ N^e \Big]\ =\ \vert S^* \vert \cdot \Big[L^e \Big]\ =\ \Big[L^e \otimes S \Big], 
\end{equation}
hence this class is independent of the empty line chosen. In any case, when constructing $L^e \otimes S$, $L^e$ takes away (``resolves'') the points of $S$, while 
duplicating itself $\vert S^* \vert$ times.  

Now suppose that $S = M \setminus M^*$, another empty line (and coming from the full proper line $M$ in some element $Y \in \mQ_\ell$). Then

\begin{equation}
\label{law2}
\Big[ L^e \otimes S \Big]\ =\ \Big[(L \setminus L^*) \otimes (M \setminus M^*) \Big]\ = \ \Big[ L \otimes M \Big]\ -\ \Big[ L^* \otimes M \Big]  \  -\ \Big[ L \otimes M^* \Big]\ +\ \Big[ L^* \otimes M^* \Big]\ =\ \Big[ \mbox{\Large $\varnothing$} \Big]. \\
\end{equation}

Combining the laws (\ref{law1}) and (\ref{law2}) in more general expressions, one indeed also has a well-defined product in $K_0(\mQ_\ell)$ if empty lines are involved. \\

It is very important to remark that in some expressions of type $L \setminus A$ in $\mQ_\ell$, when base extension (from $\ell$ to $\nabla$) is defined for $L$ and $A$ such as in the examples which we will later 
describe in section \ref{BASE}, it is possible that while in $\mQ_\ell$ we have an empty line, still $(L \otimes_\ell \nabla) \setminus (A \otimes_\ell \nabla)$ is not empty | as if an empty line is defined over ``the empty field.'' (One can for instance 
easily adapt the examples in section \ref{exunionline} for this purpose.) This is yet another motivation for defining empty lines.\\ 

\begin{remark}
\label{rem2.2}
{\rm 
Note that $\Aut(L^e)$ is isomorphic to $\{ \id \} \cong \texttt{S}_0$ (so $L^e = (L^e,\{\id\})$).}
\end{remark}


\medskip
\section{Base extension and closed points}
\label{BASE}

Suppose $\chi$ is a scheme of finite type over $\mathbb{Z}$. Let $\F_q$ be a finite field, with $q$ a power of the prime $p$. 

\subsection{Rational points}

Just for the sake of reminder, we provide the little diagram below, which depicts some relations between the points of a Zariski-topological space of an affine scheme ($\Spec(A)$, with $A$ a commutative, unital ring). \\ 

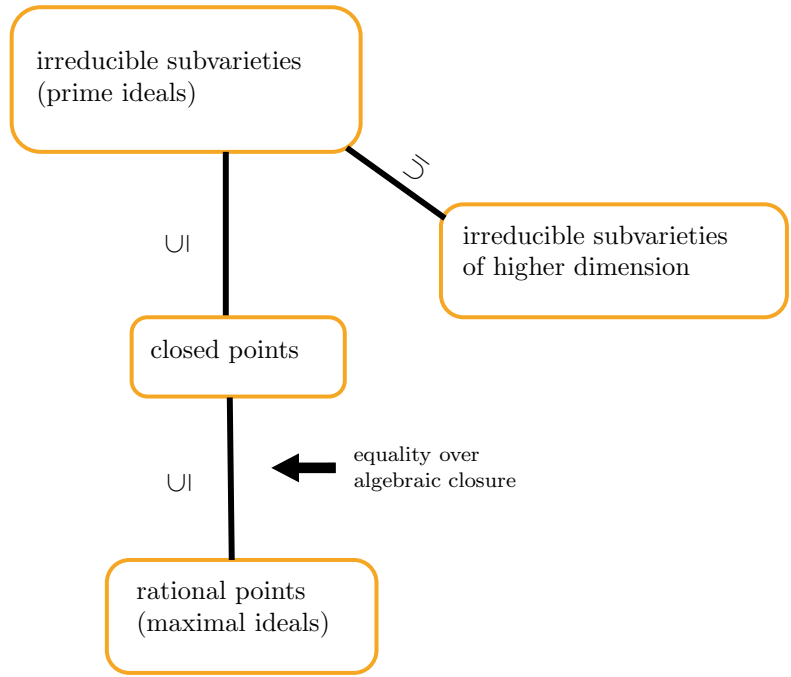
\tikzset{every picture/.style={line width=0.75pt}} 

\begin{center}
\item
\begin{tikzpicture}[x=0.75pt,y=0.75pt,yscale=-1,xscale=1]

\draw  [color={rgb, 255:red, 245; green, 166; blue, 35 }  ,draw opacity=1 ][line width=1.5]  (58,112.1) .. controls (58,104.04) and (64.54,97.5) .. (72.6,97.5) -- (219.4,97.5) .. controls (227.46,97.5) and (234,104.04) .. (234,112.1) -- (234,155.9) .. controls (234,163.96) and (227.46,170.5) .. (219.4,170.5) -- (72.6,170.5) .. controls (64.54,170.5) and (58,163.96) .. (58,155.9) -- cycle ;
\draw [line width=2.25]    (166,170.5) -- (166,253.5) ;
\draw  [color={rgb, 255:red, 245; green, 166; blue, 35 }  ,draw opacity=1 ][line width=1.5]  (118.5,262) .. controls (118.5,257.58) and (122.08,254) .. (126.5,254) -- (218,254) .. controls (222.42,254) and (226,257.58) .. (226,262) -- (226,286) .. controls (226,290.42) and (222.42,294) .. (218,294) -- (126.5,294) .. controls (122.08,294) and (118.5,290.42) .. (118.5,286) -- cycle ;
\draw  [color={rgb, 255:red, 245; green, 166; blue, 35 }  ,draw opacity=1 ][line width=1.5]  (274.5,208.4) .. controls (274.5,202.1) and (279.6,197) .. (285.9,197) -- (437.6,197) .. controls (443.9,197) and (449,202.1) .. (449,208.4) -- (449,242.6) .. controls (449,248.9) and (443.9,254) .. (437.6,254) -- (285.9,254) .. controls (279.6,254) and (274.5,248.9) .. (274.5,242.6) -- cycle ;
\draw [line width=2.25]    (227,168.5) -- (276.5,204) ;
\draw  [color={rgb, 255:red, 245; green, 166; blue, 35 }  ,draw opacity=1 ][line width=1.5]  (106,387.9) .. controls (106,381.6) and (111.1,376.5) .. (117.4,376.5) -- (230.6,376.5) .. controls (236.9,376.5) and (242,381.6) .. (242,387.9) -- (242,422.1) .. controls (242,428.4) and (236.9,433.5) .. (230.6,433.5) -- (117.4,433.5) .. controls (111.1,433.5) and (106,428.4) .. (106,422.1) -- cycle ;
\draw [line width=2.25]    (168,294.5) -- (169,376.5) ;
\draw  [fill={rgb, 255:red, 0; green, 0; blue, 0 }  ,fill opacity=1 ] (217.75,332.17) -- (202.4,332.17) -- (202.4,336.5) -- (190,330) -- (202.4,323.5) -- (202.4,327.83) -- (217.75,327.83) -- cycle ;\draw  [fill={rgb, 255:red, 0; green, 0; blue, 0 }  ,fill opacity=1 ] (221,332.17) -- (220.35,332.17) -- (220.35,327.83) -- (221,327.83) -- cycle ;\draw  [fill={rgb, 255:red, 0; green, 0; blue, 0 }  ,fill opacity=1 ] (219.7,332.17) -- (218.4,332.17) -- (218.4,327.83) -- (219.7,327.83) -- cycle ;

\draw (69,118) node [anchor=north west][inner sep=0.75pt]   [align=left] {irreducible subvarieties\\ (prime ideals)};
\draw (126.5,264) node [anchor=north west][inner sep=0.75pt]   [align=left] {closed points};
\draw (284,206) node [anchor=north west][inner sep=0.75pt]   [align=left] {irreducible subvarieties \\of higher dimension};
\draw (134.4,224) node [anchor=north west][inner sep=0.75pt]  [font=\large,rotate=-270]  {$\subseteq $};
\draw (263.7,188.99) node [anchor=north west][inner sep=0.75pt]  [font=\large,rotate=-215.93]  {$\subseteq $};
\draw (119.4,385.5) node [anchor=north west][inner sep=0.75pt]   [align=left] {rational points \\ (maximal ideals)};
\draw (135.4,344.5) node [anchor=north west][inner sep=0.75pt]  [font=\large,rotate=-270]  {$\subseteq $};
\draw (229,317.5) node [anchor=north west][inner sep=0.75pt]  [font=\footnotesize] [align=left] {equality over \\algebraic closure};
\end{tikzpicture}
\item
\item
(We see $A$ as a $k$-algebra with $k$ a field. The ``algebraic closure'' above refers to an algebraic closure of $k$.)\\

\end{center}

\bigskip
The {\em $\F_q$-rational points} of $\chi$ are given by the homomorphisms
\begin{equation}
\alpha: \Spec(\F_q) \ \mapsto\ \chi.
\end{equation}

With $\chi_q := \chi \otimes_{\Spec(\Z)}\Spec(\F_q)$, they are also given by the $\F_q$-homomorphisms $\beta: \Spec(\F_q) \mapsto \chi_q$. 
(If $\chi_q$ is a projective variety given by homogeneous equations $f_i = 0$ with the $f_i$ in $\F_q[X_0,\ldots,X_d]$, the $\F_q$-rational 
points are points which have homogeneous coordinates with entries solely in $\F_q$.) Let $\cl(\chi)$ be the set of closed points of $\chi$. Write 
$\overline{\chi} := \chi \otimes_{\Spec(\Z)}\Spec(\overline{\F_q})$, with $\overline{\F_q}$ an algebraic closure of $\F_q$. Also, let 
\begin{equation}
F: \overline{\F_q}\ \mapsto\ \overline{\F_q}: \ x\ \mapsto\ x^p 
\end{equation}
be the Frobenius automorphism. The Frobenius automorphism acts on $\chi(\overline{\F_q})$, the set of $\overline{\F_q}$-rational
points of $\chi$. For $x \in \chi(\overline{\F_q})$, we have that $x^{F^m} = x$ if and only if $x$ is an $\F_{p^m}$-rational point of $\chi$. Each $\Aut(\overline{\F_q}/\F_p)$ 
orbit is finite, and its size is $m$, if its points are $\F_{p^m}$-rational points, but not $\F_{p^n}$-rational for $n < m$.  The $\Aut(\overline{\F_q}/\F_p)$-orbits in $\chi(\overline{\F_q})$ correspond
bijectively to the closed points of $\chi$. The $\F_q$-rational points of $\chi$ (or $\chi_q$) are those closed points that are fixed by $F^h$, with $q = p^h$.

\subsection{Absolute Galois groups}
\label{absgal}

Fix a prime $p$.
Let $S \subseteq \mathbb{N}^\times$, and suppose $(S,\preceq)$ is a directed set, where $\preceq$ is defined as follows:
$a \preceq  b$ if $a$ divides $b$ (which happens if and only if $\F_{p^a}$ is a subfield of $\F_{p^b}$). Let $u \preceq v$; then we 
have a natural group morphism
\begin{equation}
\pi_{uv}:\ \Aut(\F_{p^v}/\F_p)\ \mapsto\ \Aut(\F_{p^u}/\F_p)
\end{equation}
since $\F_{p^u}$ is the unique subfield of $\F_{p^v}$ of order $p^u$;
note that in general, with $n \in \mathbb{N}^\times$, $\Aut(\F_{p^n}/\F_p)$ is a cyclic group of order $n$ which is generated by the Frobenius automorphism $F$. Also, 
note that $\pi_{vu}$ has a kernel which is isomorphic to $\Aut(\F_{p^v}/\F_{p^u})$.

We obtain a projective system $(\Aut(\F_{p^m}/\F_p),\pi_{mn})$ over $(S,\preceq)$, and the projective limit\\ 
$\varprojlim \Aut(\F_{p^m}/\F_p)$ is well defined. 
If $S = \mathbb{N}^\times$, then 
\begin{equation}
\varprojlim \Aut(\F_{p^m}/\F_p) = \Aut(\overline{\F_p}/\F_p) = \mathrm{Gal}(\overline{\F_p}/\F_p) \cong \widehat{\Z},
\end{equation}
the profinite completion of the integers.

\subsection{Galois action on maximal and prime ideals}
\label{Galois}

We first state and prove the following result. 

Fix a finite prime field $\F_p$. 
Let $i \in \Z$ be a positive integer different from $0$, and let $S := \{ n\cdot i \ \vert\ n \in \Z, n > 0 \}$, and suppose $\emptyset \ne T \subseteq S$ is such that $(T,\preceq)$ is 
a directed set, with $\preceq$ as in section \ref{absgal}.  Define the group homomorphisms 
\begin{equation}
\iota_{uv}^{\F}:\ \F_{p^u} \hookrightarrow \F_{p^u}
\end{equation}
as before, taken that $u \preceq v$ in $T$.  Define $\ell := \varinjlim_T \F_{p^j}$. 

The following theorem is well known in case of, e.g., Galois groups of algebraic closures ($\overline{k}$ of a field $k$) acting on maximal ideals (closed points) over a given maximal 
ideal in a $k$-algebra $B$.\\ 

\Big({\em Example}: action of $\mathrm{Gal}(\C/\mathbb{R})$ on the prime ideals $(x + i)$ and $(x - i)$ in the polynomial ring $\C[x]$; the orbit 
$\Big\{ (x + i), (x - i) \Big\}$ corresponds to the closed point $(x^ 2 + 1)$ in $\Spec(\mathbb{R}[x])$.\Big)\\

The fact that we are working with projective limits of Galois groups gives us enough flexibility to go a bit further. 

\begin{theorem}[Galois action of $\Aut(\ell/\F_{p^i})$]
\label{galactthm}
Suppose $\mX \mapsto \F_{p^i}$ is a scheme defined over $\F_{p^i}$.  Let $\cup X_j$ be an open cover of $\mX$ such that 
the $X_j$s are spectra $\Spec(A_j)$ of $\F_{p^i}$-algebras. Let $A$ be any such $\F_{p^i}$-algebra. Let $\mX_\ell := \mX \otimes_{\Spec(\F_{p^i})}\Spec(\ell)$ and 
$A_\ell := A \otimes_{\F_{p^i}} \ell$. Then for each prime ideal $\frak{p}$ of $A$, 
we have that $\Aut(\ell/\F_{p^i})$ acts transitively on the set of prime ideals of $\mX_\ell$ over $\frak{p}$. 
\end{theorem}

{\em Proof.}\quad
First of all, note that for any $j \in T$, we have that $\Aut(\F_{p^j}/\F_{p^i})$ acts transitively on the set of prime ideals of $A \otimes_{\F_{p^i}}\F_{p^j}$ 
over $\frak{p}$, since $\F_{p^j}/\F_{p^i}$ is a finite Galois extension. So, if $T$ is a finite set, $\ell$ is also finite, and the statement is true. So we may suppose that 
$T$ is not finite. 
Also, let $\varphi \in \Aut(\F_{p^m}/\F_{p^i})$; then $\varphi$ acts on the 
prime ideals of $A \otimes_{\F_{p^i}}\F_{p^m}$. If $n \preceq m$ in $T$, and $\varphi$ sends $\frak{p'}$ to $\varphi(\frak{p}')$, where $\frak{p}'$ is a prime 
ideal of $A \otimes_{\F_{p^i}}\F_{p^m}$, then it sends $\frak{p} \cap (A \otimes_{\F_{p^i}} \F_{p^n})$ to $\varphi(\frak{p}) \cap (A \otimes_{\F_{p^i}} \F_{p^n})$ (action of $\pi_{nm}(\varphi)$).

Now let $\frak{p}_1$ and $\frak{p}_2$ be prime ideals in $A_\ell$ over $\frak{p}$. For each $r \in T$, we have that $\frak{p}_1 \cap (A \otimes_{\F_{p^i}} \F_{p^r})$ and  
$\frak{p}_2 \cap (A \otimes_{\F_{p^i}} \F_{p^r})$ are primes over $\frak{p}$, so there is an element $\varphi_r$ sending the former prime ideal to the latter. 
Now consider an element 
\begin{equation}
\overline{\varphi} := \Big(\widetilde{\varphi_r}\Big)_{r \in T} \in \Aut(\ell/\F_{p^i}) = \varprojlim_T\Aut(\F_{p^u}/\F_{p^i}) 
\end{equation}
such that each $\widetilde{\varphi_r}$ sends $\frak{p}_1 \cap (A \otimes_{\F_{p^i}} \F_{p^r})$ to  
$\frak{p}_2 \cap (A \otimes_{\F_{p^i}} \F_{p^r})$. To show that such $\overline{\varphi}$ exist, proceed as follows. For each $u \in T$, let $P_u$ be the set of all elements in $\Aut(\F_{p^u}/\F_{p^i})$ mapping 
$\frak{p}_1 \cap (A \otimes_{\F_{p^i}} \F_{p^u})$ to  
$\frak{p}_2 \cap (A \otimes_{\F_{p^i}} \F_{p^u})$. The transition maps $\pi_{vw}$ induce maps $\underline{\pi_{vw}}: P_w \mapsto P_v$ for $v \preceq w$, making 
$\Big(P_u,\underline{\pi_{uv}}\Big)$ a projective system over $T$. Every $P_u$ is nonempty and finite, so $\varprojlim P_u$ is nonempty as well (see e.g. Serre \cite[Chapter II, \S 2.1]{cours} or \cite[Lemma 4.21.7, \texttt{https://stacks.math.columbia.edu/tag/002z}]{Stack}). Now take $\overline{\varphi} \in \varprojlim P_u$. \\

As 
\begin{equation}
A \otimes_{\F_{p^i}}\ell\ =\ \bigcup_{s \in T}\Big(A \otimes_{\F_{p^i}}\F_{p^s}\Big),
\end{equation}
it follows that $\overline{\varphi}$ sends $\frak{p}_1$ to $\frak{p}_2$. This concludes the proof of the statement. 
\eop \\



\medskip
\subsection{Base extension: an example}
\label{baseex}

In this section we explain an example of base extension of nonclassical GQs over fields containing some given field $\F_{2^r}$, which allows us to pass to an infinite subfield $\ell$ of an algebraic closure of $\F_{2^i}$. We then use the action of a suitable Galois group on the points of the topology of the quadrangle over $\ell$, so as to perform a 
kind of Galois descent, to define new points and closed sets of the initial finite GQs.

We first need some terminology. \\

\subsection{Ovals and quadrangles}
\label{T2O}

An {\em oval} of an axiomatic projective plane $\bP$ (of order at least $2$), is a set of points $\mO$ such that each line meets $\mO$ in at most two distinct points, and such that each 
point $x$ of $\mO$ is incident with precisely one line which does not meet $\mO$ in any other point (such a line is called a {\em tangent line} of $\mO$ at that point). If $\bP$ is finite of order $n$, 
the second defining property can be re-phrased as: $\vert \mO \vert = n + 1$. Still in the case where $n$ is finite, and assuming that $n$ is odd furthermore,
the latter quantity is the maximal number of points 
for a point set of $\bP$ in which no three distinct points are incident with the same line. If $n$ is even, any set $\mO$ of size $n + 1$ with this property can be extended to a set 
of size $n + 2$ still having this property; this is because one can show that the $n + 1$ tangents of $\mO$ meet in one and the same point of the plane, which we call the {\em nucleus} of $\mO$. 
A set of $n + 2$ distinct points $\mH$ in $\bP$, still assumed of finite even order $n$, such that any line meets $\mH$ in $0$ or $2$ distinct points, is called a {\em hyperoval}. 
In the general case, a {\em hyperoval} $\mH$ is a set of points in an axiomatic projective plane of order $n \geq 2$, such that any line which is incident with one point of $\mH$, is 
incident with precisely one other point of $\mH$.  \\

Now let $\mH$ be a hyperoval of a projective plane $\bP^2(k)$ over a field $k$ (in general, this does not force $k$ to have characteristic $2$). Embed 
$\bP^2(k)$ as a hyperplane in $\bP^3(k)$. Let $\eta$ be one distinguished point of $\mH$ (so that it is the nucleus of the oval $\mH \setminus \{ \eta\}$). 
Now we define the following point-line incidence structure, $\Gamma(\mH \setminus \{\eta\},\eta) = (\mP,\mL,\I)$.

The elements of $\mP$ are:
\begin{itemize}
\item[(1)]
the points of $\bP^3(k) \setminus \bP^2(k)$;
\item[(2)]
the planes in $\bP^3(k)$ which meet $\bP^2(k)$ in a line containing $\eta$ (and hence an other point of $\mH$);
\item[(3)]
a symbol $(\infty)$.
\end{itemize}

Lines come in two guises:
\begin{itemize}
\item[(a)]
lines of $\bP^3(k)$ that meet $\bP^2(k)$ in a point of $\mH \setminus \{\eta\}$;
\item[(b)]
the points of $\mH \setminus \{\eta\}$.
\end{itemize}

Incidence works as follows: the symbol $(\infty)$ is incident with all lines of type (b); points of type (1) are incident with the lines of type (a) in which they 
are contained as subspaces; finally, points of type (2) are incident with the one line of type (b) they contain as subspaces, and with all the lines of type (b) 
which they contain as subspaces.

The one proves that $\Gamma(\mathcal{H} \setminus \{ \eta\},\eta)$ is a GQ of order $(\vert k \vert,\vert k \vert)$.

\subsubsection{Direct limit of ovals}

The starting point for our construction is the following result of Segre. 

\begin{theorem}[Segre]
\label{segre}
Let $i$ and $h$ be positive integers such that $\mathrm{gcd}(i,h) = 1$, and let $q = 2^h$. Then the set of points in $\bP^2(q)$ given 
by homogeneous coordinates
\begin{equation}
\{ (1 : t : t^{2^i})\ \vert \ t \in \F_q \}\ \cup\ \{ (0 : 0 : 1) \}
\end{equation}
is an oval $\mO(i,h)$ with nucleus $\eta = (0 : 1 : 0)$. 
\end{theorem}

If $\mathrm{gcd}(i,h) \ne 1$, then the set $\{ (1 : t : t^{2^i})\ \vert \ t \in \F_q \} \cup \{ (0 : 0 : 1), (0 : 1 : 0)\}$ does not form a hyperoval; see Hirschfeld \cite[\S 8.4, Corollary 3]{Hirsch}. 

Each of the ovals constructed from Theorem \ref{segre} has as nucleus $\eta = (0 : 1 : 0)$ (in planes represented over different fields, that is to say). 

Now define a direct system of fields, as follows. First of all, fix one positive integer $i > 1$ (for the sake of convenience, this could be some prime). Now let $\mS(i)$ 
be the set of all positive integers $r$ such that $\mathrm{gcd}(r,i) = 1$. We define a directed set $(\mS(i),\preceq)$ as follows:  $a \preceq  b$ if $a$ divides $b$
(which happens if and only if $\F_{2^a}$ is a subfield of $\F_{2^b}$). We now consider the family $\{ \mO(i,h)\ \vert \ h \in \mS(i) \}$. If $u \preceq v$, 
then define the natural embeddings
\begin{equation}
\begin{cases}
&\iota^\mO_{uv} := \mO(i,u) \ \hookrightarrow\ \mO(i,v), \\
&\iota^\F_{uv}  := \F_{2^u} \ \hookrightarrow\ \F_{2^v}, \\
&\iota^\bP_{uv} := \bP^2(2^u) \ \hookrightarrow\ \bP^2(2^v).
\end{cases}
\end{equation}

Obviously $\Big( \mO(i,u),\iota_{uv}^{\mO} \Big)$ is a direct system over $\mS(i)$,   as well as $\Big( \F_{2^u},\iota_{uv}^{\F} \Big)$ and  $\Big( \bP^2(2^u),\iota_{uv}^{\bP} \Big)$. 

\begin{theorem}
We have that $\varinjlim\mO(i,u)$ is an oval of $\bP^2(\ell)$, with nucleus $(0: 1 : 0)$; here, $\ell = \varinjlim\F_{2^u}$ and $\bP^2(\ell) = \varinjlim\bP^2(2^u)$. 
\end{theorem}

{\em Proof.}\quad
We only need to verify that $\widetilde{\mO} := \varinjlim\mO(i,u)$ is an oval in $\bP^2(\ell)$ with nucleus $(0 : 1 : 0)$; the rest is obvious. 
Note that in $\bP^2(\ell)$
\begin{equation}
\label{equnion}
\widetilde{\mO}\ =\ \bigcup_{r \in \mS(i)}\mO(i,r),
\end{equation}
where the homogeneous coordinates are chosen with respect to the same base in each of the planes $\bP^2(2^u)$. 

Now consider any point $x$ of $\widetilde{\mO}$, and suppose $L \I x$ is a line of $\bP^2(\ell)$ which does not contain $(0 : 1 : 0)$, nor any other point of $\widetilde{\mO}$. Note that 
$(0 : 1 : 0)$ is contained in any of the planes $\bP^2(2^h)$. 
Let $y \ne x$ be any other point incident with $L$. By eq. (\ref{equnion}), there is some $m \in \mS(i)$ such that $(0 : 1 : 0)$, $x$ and $y$ are points of $\bP^2(2^m)$. 
So $L$ is also a line of $\bP^2(2^m)$. Since $\mO(i,m)$ is an oval with nucleus $(0 : 1 : 0)$, we conclude that $L$ contains at least one point of $\mO(i,m) \cup \{(0 : 1 : 0)\}$ besides 
$x$, contradiction. We conclude that any line which meets $\widetilde{\mO}$ in at least one point, contains at least one other point of $\widetilde{\mO} \cup \{ (0 : 1 : 0) \}$. 
By a similar argument, it follows that any line which meets $\widetilde{\mO}$ in at least one point, contains precisely two distinct points of $\widetilde{\mO} \cup \{ (0 : 1 : 0) \}$. 
The statement is proved. 
\eop \\

\medskip
\subsubsection{The quadrangles}

Passing to the corresponding GQs, we obtain a direct system of GQs $\Gamma(\mO(i,u),\eta)$ over the directed set $\mS(i)$, and 
\begin{equation}
\varinjlim \Gamma(\mO(i,u),\eta) = \Gamma(\widetilde{\mO},\eta). 
\end{equation}

Note that $\Gamma(\widetilde{\mO},\eta) = \cup_{u \in \mS(i)}\Gamma(\mO(i,u),\eta)$.

Although $\ell$ is not algebraically closed, it comes as close to being to algebraically closed as the property that the sets $\mO(i,r)$ should all be ovals | and 
the sets of $\Gamma(\mO(i,h),\eta)$ be elements of $\mQ_{\F_{2^h}}$ | allows. So we view $\Gamma(\widetilde{\mO},\eta)$ as a synthetic analogon 
of a base extension to an algebraically closed field of a given scheme over a finite field. \\

\medskip
\subsubsection{Galois action and enriching}
\label{galact}

Now define $\USpec(\Gamma(\widetilde{\mO},\eta))$ as before. In the context of the example under review, we can ``enrich'' the Zariski topology of the finite quadrangles 
$\Gamma(\mO(i,h),\eta)$ via the one of $\USpec(\Gamma(\widetilde{\mO},\eta))$, as follows.  

\bigskip
\fbox{\begin{minipage}{5cm}
{\bf CLOSED POINTS}
\end{minipage}}

First of all, consider the closed points in $\USpec(\Gamma(\widetilde{\mO},\eta))$. Following our definition of Zariski topology, we have three types of closed points 
in $\USpec(\Gamma(\widetilde{\mO},\eta)) =: \overline{\Gamma}$: the symbol $(\infty)$, the planes which meet $\widetilde{\mO}$ in a tangent line (``planar points''), and 
the points of $\bP^3(\ell) \setminus \bP^3(\ell)$ (``affine points''). Recall that in our context $\ell$ plays the \ul{role of algebraic closure} of the finite fields in consideration. 

Consider $\Gamma(\mO(i,h),\eta)$. Let $T := \{ n\cdot h\ \vert\ n \in \Z, n > 0 \} \cap \mS(i)$, and note that 
\begin{equation}
\varinjlim_T\F_{2^u} = \varinjlim_{\mS(i)}\F_{2^u} = \ell.
\end{equation}

\ul{Descent on affine points.}\quad 
The closed points of affine type of $\Gamma(\mO(i,h),\eta)$ are given by $\Aut(\ell/\F_{2^h})$-orbits in the affine point set of $\Gamma(\widetilde{\mO},\eta)$ (so in the set 
$\bP^3(\ell) \setminus \bP^2(\ell)$). 

Note that $\bP^3(2^h) \setminus \bP^2(2^h)$ is a subset (if its elements are seen as singletons): they are the fixed points of the Frobenius automorphism
\begin{equation}
F:\ \mathbf{x} \mapsto \mathbf{x}^{2^h},
\end{equation}
where $\mathbf{x}$ is a point of $\bP^3(\ell) \setminus \bP^2(\ell)$ expressed with homogeneous coordinates, and $F$ acts on these coordinates. \\

\ul{Descent on planar points}.\quad
Points of this type are $\Aut(\ell/\F_{2^h})$-orbits of planar points of $\Gamma(\widetilde{\mO},\eta)$; after introducing homogeneous coordinates, elements 
of $\Aut(\ell/\F_{2^h})$ act on the points of the corresponding planes as in the case of affine points, and the planar points of $\Gamma(\mO(i,h),\eta)$ are singled 
out by the Frobenius automorphism in a similar way. \\

\ul{The symbol $(\infty)$}.\quad
We do not define a nontrivial descent on this point. The reason is obvious: it should be considered as a rational point of $\Gamma(\mO(i,h),\eta)$, and 
so constitutes a trivial orbit under $\Aut(\ell/\F_{2^h})$. \\

\bigskip
\fbox{\begin{minipage}{5cm}
{\bf OTHER PRIME IDEALS}
\end{minipage}}

\medskip
Other prime ideals in the new setting are defined in a similar way. Consider for instance a full subgrid in $\Gamma(\widetilde{\mO},\eta)$; theoretically, such a 
grid could be of several different types, but the generic type is the following: consider a plane $\delta$ in $\bP^3(\ell)$ which meets $\widetilde{\mO}$ (precisely) in two different points $u$ and $v$;
then $\delta$ defines such a grid. The lines of the grid are the lines $u$ and $v$ of the quadrangle, and the lines of $\delta$ which meet $\widetilde{\mO}$ precisely in $u$ or 
in $v$ (so they are not lines of $\bP^2(\ell)$); the points are the affine points in $\delta$, and the planes of $\bP^3(\ell)$ which are different from $\bP^2(\ell)$ which 
contain either $u\eta$ or $v\eta$. 

Consider $\mO(i,h)$ for some fixed $i$ and $h$. Then note that with the notation of above, we could have that two elements, resp. one element, resp. no elements in $\{ u, v \}$ are contained in $\mO(i,h)$ (for instance, $\delta$ could have an empty intersection with $\mO(i,h)$). So different types of grids arise. On each such grid we can then consider the action of $\mathrm{Gal}(\ell/\F_{2^h})$ (keeping Theorem \ref{galactthm} in the back of our mind). 
Similar remarks can be made about, e. g., lines arising from enriching (see also the next section).

\bigskip
\fbox{\begin{minipage}{5cm}
{\bf NEW CLOSED SETS}
\end{minipage}}

\medskip
\subsubsection{Lines}

If we look at generic lines (those not incident with $(\infty)$), then we look at all such lines of $\Gamma(\widetilde{\mO},\eta)$ | call the set of all such lines $\widetilde{\mL}$ | and consider the action of $\Aut(\ell/\F_{p^h})$ on $\widetilde{\mL}$. (Note that $\Aut(\ell/\F_{p^h})$ also acts on the points of $\widetilde{\mO}$ | see the next section.) 
Then the $\Aut(\ell/\F_{p^h})$-orbits in $\widetilde{\mL}$ define the new lines at level $\F_{p^h}$. We can single out the lines of $\Gamma(\mO(i,h),\eta)$ which are not incident with $(\infty)$ as those which are fixed by each element of $\Aut(\ell/\F_{p^h})$; similarly we can single out the lines over field extensions of $\F_{p^h}$. We can work in the same way for other prototypes of closed sets, making distinction between ``rational sets'' and ``orbit sets.''

One question which remains to be considered is the following variation:
\begin{question}
Let $\I$ be an ideal geometry $($over $\F_q$$)$; how does $C(\I)$ look like in the extended framework?
\end{question}

In case of schemes, ideals (varieties) are not defined by their rational points: an ideal corresponds to a set of equations which remains the same after base extension. In our case, an analogue to an ideal such as $(X^2 + 1)$ cannot be given (if we imagine to be working over $\mathbb{R}$), as we only work with ``points coming from Incidence Geometry'' (here: rational points).  Instead, we should be {\em given} the ideals $(X + i)$ and $(X - i)$ plus the action of $\mathrm{Gal}(\C/\mathbb{R})$ to define $(X^2 + 1)$ as an orbit. And once an ideal geometry is given (over $\F_q$), matters are more complicated. We propose the following variation:\\ 

{\bf PROCEDURE}.\quad
The geometry $\I$ is given over $\F_q$ (as a point-line subgeometry of $\Gamma(\mO(i,h),\eta)$). Over $\ell$, we need to consider all point-line subgeometries $\widetilde{\I}$ of $\Gamma(\widetilde{\mO},\eta)$ which only contain elements  of $\I$ over $\F_q$, and which are on the other hand in some sense {\em contained in $\I$ once $\I$ is considered over $\ell$.} 
(noting that $\I$ is also a point-line subgeometry of $\Gamma(\widetilde{\mO},\eta)$). We first define $\overline{\I} = \I \otimes_{\F_q} \ell$, the base extension of $\I$ to 
$\ell$.  We impose the following rules to obtain $\overline{\I}$. \\

\begin{itemize}
\item[\framebox{C1}]
{\em Lines}. If $U$ is a line in $\I$, it is full over $\F_q$ by definition. Over $\ell$, $U \otimes_\ell \ell$ also becomes a full line in $\overline{\I}$.  \\
\item[\framebox{C2}]
{\em Ideal perps}. As in the case of full lines, ideal perps in $\I$ remain ideal perps over $\ell$ in $\overline{\I}$. \\
\item[\framebox{C2.B}]
{\em Non-ideal perps}. Suppose $\mG$ is a non-ideal perp in $\I$ with base point $x$. 
Let $V$ be the set of all full (thin or thick) subGQs in $\I$ (with an order) containing $x$. Then $\mG \otimes \ell$ is the union of all perp extensions 
performed in all $\mS' \in V$ in which $\mG \cap \mS'$ is ideal (following (C2) and (C3) below), so that we obtain 
\[ \bigcup_{\mS' \in V, \mS' \cap \mG\ \mbox{ideal in}\ \mS'}\Big(\mG \cap \mS'\Big)\otimes \ell.\] 
The set of lines $W \I x$, where $W$ is a line in $\mG$ which is not contained in an element $\mS' \in V$ for which $\mG \cap \mS'$ is ideal in $\mS'$, remains unchanged, and is extended following (C1).   
\\
\item[\framebox{C3}]
{\em Subquadrangles}. Generalizing (C1) and (C2), each full subquadrangle $Y$ (both thick and thin) which is contained in $\I$ (over $\F_q$), becomes a full subquadrangle 
of $\overline{\I}$ over $\ell$, under the following agreement: $Y \otimes_\ell \ell$ is the generalized quadrangle generated by the set of full lines over $\ell$ of $Y$. (So by definition, $\overline{\I}$ contains $Y \otimes_\ell \ell$.) \\
\end{itemize}

By (C3), full grids become full grids.

Now consider again any ideal geometry $\I$ over $\F_q$. 
Then $\overline{\I}$ defines a closed set $C(\overline{\I})$ (over $\ell$), and using descent we obtain a closed set $\underline{C(\overline{\I})}$ over $\F_q$. 
Note that the action of $\texttt{Gal}(\ell/\F_q)$ is well defined, since it acts on $\widetilde{\Gamma}$ (so images of ideal geometries under elements of $\texttt{Gal}(\ell/\F_q)$ are again subgeometries of $\widetilde{\Gamma}$). In short, we have the following diagram:

\begin{equation}
\I \ \ \xrightarrow[]{\otimes_\ell}\ \ \overline{\I} \ \ \xrightarrow[]{\text{in}\ \USpec(\Gamma(\widetilde{\mO},\eta))} \ \ C(\overline{\I}) \ \ \xrightarrow[]{\text{descent by}\ \texttt{Gal}(\ell/\F_q)}\ \ C(\I) = \underline{C(\overline{\I})}.
\end{equation}

\medskip
\subsubsection{Comparison to schemes}

Fix a field $k = \F_{2^h}$, and let $i$ be a positive integer such that $\mathrm{gcd}(i,h) = 1$. Then observe that $\mO(i,h) \cup \{ (0 : 0 : 1) \}$ is the set of $k$-rational points of 
the algebraic curve $\mC$ defined by the equation
\begin{equation}
Y^{2^i} = ZX^{2^i - 1}. 
\end{equation}

Let $\ell$ be, as before, the limit $\varinjlim\F_{2^u}$. Then the $\ell$-rational points of the scheme 
\begin{equation}
\Proj\Big(\ell[X,Y,Z]/(Y^{2^i} - ZX^{2^j - 1})\Big) 
\end{equation}
are precisely 
the points of $\widetilde{\mO} \cup \{ (0 : 0 : 1)\}$.

\medskip
\subsubsection{General picture for GQs coming from ovals}

The example given in this section easily generalizes to other GQs coming from ovals. We leave the details to the reader. 


\medskip
\subsubsection{General recipe for enriching the Zariski topology}

Of course one has to be cautious with this enriching procedure, since it \ul{depends on the class in which the GQ is considered}. \\

Let  $h$ be a positive odd integer, and let $q = 2^h$. Then the set of points in $\bP^2(q)$ given 
by homogeneous coordinates
\begin{equation}
\mO(h)\ := \ \{ (1 : t : t^6)\ \vert \ t \in \F_q \}\ \cup\ \{ (0 : 0 : 1) \}
\end{equation}
is an oval in $\bP^2(q)$, with nucleus $\eta = (0 : 1 : 0)$ (see Brown \cite{brown} for all the details, and more).\\

Similarly as above, we put $\ell = \cup_{h\ \mbox{odd}}\F_{2^h}$; it is an infinite field contained in an algebraic closure $\overline{\F_2}$ of $\F_2$ (union of finite extensions of $\F_2$). Denote the oval which arises similarly as above in $\bP^2(\ell)$, by $\mO(\ell)$. Now let $h = 1$: then $\Gamma(\mO(1),\eta)$ is a GQ of order $(2,2)$ which is isomorphic to $\mQ(4,2)$. Now endow $\mQ(4,2)$ with the enriched topology $\tau_1$ coming from $\Gamma(\mO(\ell),\eta)$. On the other hand, we can endow $\mQ(4,2)$ also with the ``full Zariski topology'' coming from $\mQ(4,\overline{\F_2})$ (seen as projective variety); this is nothing else than the ``full Zariski topology'' $\tau_2$ of $\mQ(4,2)$ (seen as projective variety). The topologies $\tau_1$ and $\tau_2$ are different, and both contain the synthetic Zariski topology which we defined without base extension. It thus makes sense to mention the construction class $\mC$ in which we see an example, so as to ``$\mC$-enrich'' its topology as above. \\

\bigskip
\tikzset{every picture/.style={line width=0.75pt}} 
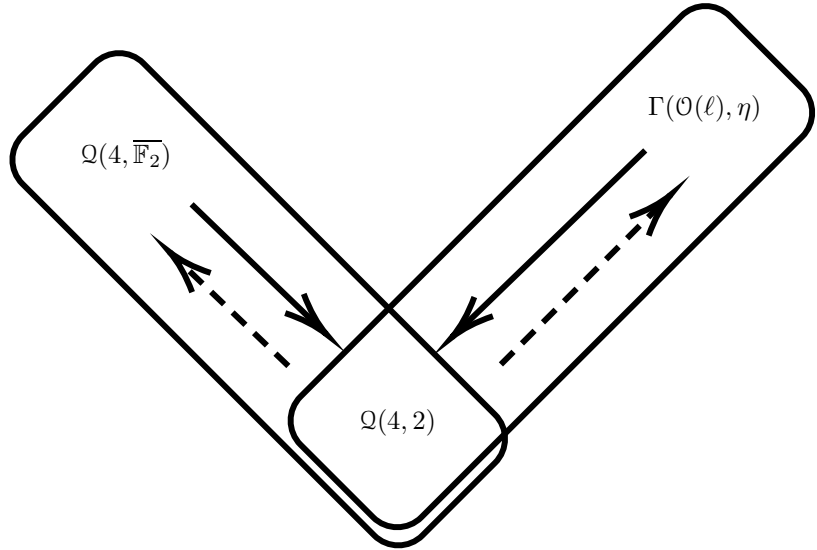
\begin{center}
\item
\begin{tikzpicture}[x=0.75pt,y=0.75pt,yscale=-0.8,xscale=0.8]

\draw  [line width=2.25]  (134.06,72.94) .. controls (142.48,64.52) and (156.12,64.52) .. (164.53,72.94) -- (386.27,294.68) .. controls (394.69,303.09) and (394.69,316.73) .. (386.27,325.15) -- (340.57,370.85) .. controls (332.15,379.27) and (318.51,379.27) .. (310.1,370.85) -- (88.36,149.11) .. controls (79.94,140.7) and (79.94,127.06) .. (88.36,118.64) -- cycle ;
\draw  [line width=2.25]  (579.91,88.91) .. controls (588.32,97.32) and (588.32,110.96) .. (579.91,119.37) -- (339.87,359.41) .. controls (331.46,367.82) and (317.82,367.82) .. (309.41,359.41) -- (263.73,313.73) .. controls (255.32,305.32) and (255.32,291.68) .. (263.73,283.27) -- (503.77,43.23) .. controls (512.18,34.82) and (525.82,34.82) .. (534.23,43.23) -- cycle ;
\draw [line width=2.25]  [dash pattern={on 6.75pt off 4.5pt}]  (256,264) -- (183.89,194.77) ;
\draw [shift={(181,192)}, rotate = 403.83000000000004] [color={rgb, 255:red, 0; green, 0; blue, 0 }  ][line width=2.25]    (27.98,-8.42) .. controls (17.79,-3.57) and (8.47,-0.77) .. (0,0) .. controls (8.47,0.77) and (17.79,3.57) .. (27.98,8.42)   ;
\draw [line width=2.25]  [dash pattern={on 6.75pt off 4.5pt}]  (390,263) -- (494.17,158.83) ;
\draw [shift={(497,156)}, rotate = 495] [color={rgb, 255:red, 0; green, 0; blue, 0 }  ][line width=2.25]    (27.98,-8.42) .. controls (17.79,-3.57) and (8.47,-0.77) .. (0,0) .. controls (8.47,0.77) and (17.79,3.57) .. (27.98,8.42)   ;
\draw [line width=2.25]    (195,162) -- (277.12,241.22) ;
\draw [shift={(280,244)}, rotate = 223.97] [color={rgb, 255:red, 0; green, 0; blue, 0 }  ][line width=2.25]    (27.98,-8.42) .. controls (17.79,-3.57) and (8.47,-0.77) .. (0,0) .. controls (8.47,0.77) and (17.79,3.57) .. (27.98,8.42)   ;
\draw [line width=2.25]    (481,128) -- (361.88,243.22) ;
\draw [shift={(359,246)}, rotate = 315.95] [color={rgb, 255:red, 0; green, 0; blue, 0 }  ][line width=2.25]    (27.98,-8.42) .. controls (17.79,-3.57) and (8.47,-0.77) .. (0,0) .. controls (8.47,0.77) and (17.79,3.57) .. (27.98,8.42)   ;

\draw (123,119.4) node [anchor=north west][inner sep=0.75pt]    {$\mathcal{Q}(4,\overline{\mathbb{F}_{2}})$};
\draw (299,289.4) node [anchor=north west][inner sep=0.75pt]    {$\mathcal{Q}(4,2)$};
\draw (481,91.4) node [anchor=north west][inner sep=0.75pt]    {$\Gamma(\mathcal{O}(\ell),\eta )$};

\end{tikzpicture}
\item
\item
Different base extensions can give rise to different topologies. In this picture, the dotted arrows stand for base extension; the other arrows for Galois descent. 
\end{center}

\subsection{Variation}

Consider a nonclassical Kantor-Knuth generalized quadrangle $\mS = \mK(q,n,\gamma)$ ($\gamma \ne \id$). As we have seen, it is the union of $\mQ(4,q)$-subGQs (in that each point and line is in at least one such subGQ). 

\begin{question}
Can we endow each $\mQ(4,q)$ with the classical Zariski topology $($where we see $\mQ(4,q)$ as a projective variety$)$, such that the subspace topology induced on each intersection of $\mQ(4,q)$-subGQs is the same $($and hence yielding a topology on $\mS$$)$? 
\end{question}

Note that for arbitrary $\mQ(4,q)$-subGQs $\mQ_1$ and $\mQ_2$ contained in the same subGQ-orbit $\Omega_i$, and which intersect in a grid $\Gamma$, there is no automorphism of $\mS$ which maps $\mQ_1$ to $\mQ_2$ and which leaves each point and line of $\Gamma$ invariant (we explore this property further in this paragraph). (We do not claim that the existence of such automorphisms is necessary for obtaining the exact same topology on $\Gamma$ whether we start from $\mQ_1$ or $\mQ_2$, but if it would exist, the property would come for free in this specific case.) In fact, since we have the different subGQ-orbits $\Omega_1$ and $\Omega_2$, there is even no $\alpha \in \Aut(\mS)_\Gamma$ sending $\mQ_1$ to $\mQ_2$ when 
$\mQ_1$ and $\mQ_2$ are picked in respectively $\Omega_1$ and $\Omega_2$. So we could only look at $\Omega_1$ for that matter. Here is a more 
detailed analysis: in $\Omega_1$ only two different $\mQ(4,q)$-subGQs contain $\Gamma$ | say $\mS_1$ and $\mS_2$. The (nontrivial) involution $\sigma_1$ which fixes $\mS_1$ elementwise stabilizes $\mS_2$, and the same property holds for the involution $\sigma_2$ which fixes $\mS_2$ elementwise. The group $\epsilon := \langle \sigma_1, \sigma_2 \rangle$ has size $4$ and is the full group of automorphisms of $\mS$ which fixes $\Gamma$ elementwise. It also fixes each of $\mS_1, \mS_2$. 
For any $\mS' \in \Omega_2$, there is only one element in $\Aut(\mS')$ (which also lives in $\Aut(\mS)_{\mS'}$) which fixes $\Gamma$ pointwise,  so since no element in $\Omega_2$ is fixed elementwise by a nontrivial automorphism of $\mS$, it follows that $\vert {\mS'}^\epsilon \vert = 2$.  So even in the orbits $\Omega_1$ and $\Omega_2$, our idea cannot work. \\

Now suppose again that $\mQ_1 \cap \mQ_2 = \Gamma$ is a grid. Let $S_1 = x^\perp \cap \Gamma$ with $x$ a point in $\mQ_1 \setminus \Gamma$; is there a point $y \in  \mQ_2 \setminus \Gamma$ such that 
\begin{equation}
y^\perp \cap \Gamma\ =\ S_1?
\end{equation}
From this perspective, it seems even more pessimistic: there are precisely $q + 1$ $\mQ(4,q)$-subGQs $\mQ_i$ which intersect two-by-two in $\Gamma$, and one can show that inside each $\mQ_i$, we have that for each point $z \not\in \Gamma$ there is exactly one second point $z' \not\in \Gamma$ such that 
\begin{equation}
z^\perp \cap \Gamma\ =\ {z'}^\perp \cap \Gamma.
\end{equation} 
But by \cite[section 1.2.4]{PTsec} the maximum number of points of $\mS$ in $S_1^{\perp}$ is $q + 1$. So there is no way that in general, closed sets such as $x^\perp \cap \Gamma$ (in $\mQ_1$) can be obtained in the same way through $\mQ_2$. We have met this situation several times in the course of this paper: the local geometry induced by $\mQ_1$ on $\Gamma$ and the local geometry induced by $\mQ_2$ can in principle be different. In this specific classical case, we know that there is an isomorphism 
\[ \gamma: \mQ_1 \mapsto \mQ_2 \] 
for which $\Gamma^\gamma = \Gamma$ ($\gamma$ is not assumed to be in $\Aut(\mS)$, of course), so  we also know the induced geometries are isomorphic. If, more generally, $\mQ_1$ and $\mQ_2$ would merely be assumed to be isomorphic nonclassical quadrangles meeting in the common full subgrid $\Gamma$, with both $\mQ_1$ and $\mQ_2$ a subquadrangle of some larger quadrangle, then there is no need at all that the induced geometries on $\Gamma$ by projection from external points in $\mQ_i \setminus \Gamma$ as above, be isomorphic. 

Especially in the infinite formulation of this discussion, it might be easy enough to construct such examples.


\newpage
\appendix

\section{Quadrics in the Grothendieck ring}
\label{quadapp}

In the next theorem, $k$ is a finite field. Only in the case of dimension $1$, we will work with general fields. An adaptation of the proof 
yields the same result for fields with $\mu$-invariant at most $2$; see \cite{MOquad}. 

This section was inspired by \cite{MOquad}, but contains more details 
in the dimension 1 case. (Only the non-algebraic geometer will benefit from these details.)  Below, $K_0(\mV_k)$ denotes the Grothendieck ring 
of $k$-varieties. (Isomorphisms are taken in the category of $k$-varieties $\mV_k$.)  

\begin{theorem}
For every quadric hypersurface $\mQ_n$ in $\mathbb{P}^n(k)$, $n \geq 1$ and $k$ a finite field, there exists a field extension $K / k$ of degree at most $2$, such that the class 
$[\mQ_n]$ is contained in the $\mathbb{Z}[\mathbb{L}]$-module generated by $1$ and the class $[\Spec(K)]$ in the Grothendieck ring $K_0(\mV_k)$.
\end{theorem}
{\em Proof.}\quad 
The proof is by induction. First we handle dimension $1$. 

\section*{Dimension $n = 1$ case} 

Consider a quadric $\mQ_1$ in $\mathbb{P}^1(k)$, and suppose it is given by the quadratic polynomial equation $f(x,y) = f = 0$ (so that $\mQ_1 = \Proj\Big(k[x,y]/(f)\Big)$). 

Below, we use the fact that if $\ell$ is a field, and $A$ a finitely generated $k$-algebra, then $\Proj(A[x]) \cong \Spec(A)$. 

We have a number of possibilities, depending on how $f$ factors in $k[x,y]$. 

\begin{itemize}
\item
{\em $f$ is irreducible}. 
In that case, we have that 
\begin{align}
\mQ_1 &= \Proj\Big(k[x,y]/(x^2 + bxy + cy^2)\Big)\ \cong\ \Spec\Big(k[x]/(x^2 + bx + c)\Big)\ \cong\ \Spec(K), 
\end{align}
with $K$ a field extension of degree $2$ of $k$. Note that $\Proj\Big(k[x,y]/(x^2 + bxy + cy^2)\Big)$ and $\Spec\Big(k[x]/(x^2 + bx + c)\Big)$ both 
are one-point schemes with structure sheaf $$\mathrm{pt} \mapsto K,$$ 
and that  $\Spec\Big(k[x]/(x^2 + bx + c)\Big)$ is obtained as the complement of $y = 0$ in $\Proj\Big(k[x,y]/(x^2 + bxy + cy^2)\Big)$. 

\item
{\em $f = a(x + dy)(x + ey)$ with $e \ne d$.} 
We obtain two different points: by the Chinese remainder theorem and the fact that for commutative rings $A, B$ 
we have that $(A \times B)[x] \cong A[x] \times B[x]$, the following identities are true in the Grothendieck ring: 
\begin{align}
\mQ_1 &= \Proj\Big(k[x,y]/(a(x + cy)(x + ey))\Big)\ \cong\ \Proj\Big((k[x,y]/(x + cy)) \times (k[x,y]/(x + ey))\Big) \nonumber \\ 
&\cong \ \Proj\Big( k[x] \times k[x]\Big)\ \cong 
\ \Proj\Big( (k \times k)[x] \Big) 
 \cong\ \Spec(k \times k), 
\end{align}
and hence $\mQ_1 \cong \Spec(k \times k)$. (The product used is the $k$-algebra product.) So
\begin{align}
[\mQ_1]\ &=\ [\Spec(k \times k)]\ =\ [\Spec(k) \coprod \Spec(k)] \nonumber  \\
&= [\Spec(k)]\ +\ [\Spec(k)]\ =\ 2.
\end{align} 
\item
{\em $f = a(x + dy)^2$.} Then 
\begin{equation}
\mQ_1 = \Proj\Big(k[x,y]/(a(x + cy)^2)\Big) \cong \Proj\Big(   k[\epsilon,\delta]/(\epsilon^2)  \Big) \cong \Spec\Big(k[\epsilon]/(\epsilon^2)\Big), 
\end{equation}
so that $[\mQ_1] = 1$. 
\end{itemize}

\section*{Dimension $n = 2$ case}

As a general reference for this part of the proof, consult \cite[chapter 22]{HT1}. \\

Now consider a quadric $\mQ_2$ in $\mathbb{P}^2(k)$, and suppose it is given by the quadratic polynomial equation $f(x,y) = f = 0$ (so $\mQ_2 = \Proj\Big(k[x,y,z]/(f)\Big)$). Below, we will use the fact that for a $k$-scheme $\mX$, its class $[\mX]$ 
coincides with the class $[\mX_{\mathrm{red}}]$ of its reduced scheme. (The latter's topology is homeomorphic to that of 
$\mX$, and its sheaf of rings is given by reducing the rings in the structure sheaf of $\mX$ | that is, each ring in the structure sheaf of $\mX$ is replaced by the quotient of the ring by its ideal of nilpotent elements.) 

Again, we need to consider some cases. 

\begin{itemize}
\item
{\em $\mQ_2$ is a conic.} 
In that case, we have that $\mQ_2 \cong \mathbb{P}^1(k)$, so that $[\mQ_2] = \mathbb{L} + 1$.
\item
{\em $\mQ_2$ is a cone $u\mQ$, where $\mQ$ is isomorphic to $\Proj(k[u,v]/(g))$, with $g = a(x + dy)(x + ey)$ with $e \ne d$.} 
In that case, $[\mQ_2] = 1 + \mathbb{L}[\mQ] = 1 + 2\mathbb{L}$.  
\item
{\em $\mQ_2$ is isomorphic to $\Proj(k[x,y,z]/(x^2 + bxy + cy^2))$ with $x^2 + bxy + cy^2$ irreducible over $k$.} 
Then 
\begin{align*}
[\mQ_2]\ &=\ [\mQ_2 \cap (x = 0)]\ + \ [\mQ_2 \setminus (\mQ_2 \cap (x = 0))]  \\ 
&=\ [\Proj\Big(k[y,z]/(y^2)\Big)]\ +\ [\Spec\Big(k[y,z]/(cy^2 + by + 1)\Big)]\ \\
&= \ [\Proj\Big(k[y,z]/(y)\Big)]\ +\ [\Spec\Big(k[y]/(cy^2 + by + 1) \otimes k[z]\Big)]\ \\ 
&= \ [\Proj\Big(k[y]\Big)\ +\ [\Spec\Big(k[y]/(cy^2 + by + 1)\Big) \times \Spec\Big( k[z]\Big)] \ \\
&=\ 1 + \ [\Spec\Big(k[y]/(cy^2 + by + 1)\Big)] \cdot [\Spec\Big( k[z]\Big)]\ \\
&=\ 1 + [\Spec(K)]\mathbb{L}
\end{align*}
with $K / k$ a field extension of degree $2$.  Note that $\Proj(k[y,z]/(y))$ is the reduced scheme of $\Proj(k[y,z]/(y^2))$. (Note that $\mQ_2$ is a cone $u\mQ$ with $u$ a point and $\mQ \cong \Proj\Big(k[x,y]/(x^2 + bxy + cy^2)\Big)$ with $x^2 + bxy + cy^2$ irreducible.)

\item
{\em $\mQ_2 \cong \Proj(k[x,y,z]/(z^2))$.} 
Then 
\begin{align*}
[\mQ_2]\ &=\ [\mQ_2 \cap (x = 0)]\ + \ [\mQ_2 \setminus (\mQ_2 \cap (x = 0))]\ \\
&= \ [\mathrm{point}]\ +\ [\Spec\Big(k(y,z)/(z^2)\Big)]\ = \ 1\ + \ [\Spec\Big(k(y,z)/(z)\Big)], 
\end{align*}
noting that $\Spec(k(y,z)/(z))$ is the reduced scheme of $\Spec(k[y,z]/(z^2))$. So $[\mQ_2]  = \mathbb{L} + 1$. 
\end{itemize}

\section*{Dimension $n \geq 3$ case}

As a general reference for this part of the proof, consult \cite[chapter 22]{HT1}. \\ 

Let $u$ be a rational point of $\mQ_n$ (and note that such points exist since $n \geq 2$ and $k$ is finite), and project $\mQ_n$ onto a $\mathbb{P}^{n - 1}(k)  \subset \mathbb{P}^n(k)$ which 
does not contain $u$. If $\Pi_u$ is the tangent hyperplane of $\mQ_n$ at $u$, then this projection away from $u$ yields an 
isomorphism between $\mQ_n \setminus (\Pi_u \cap \mQ_n)$ and $\mathbb{P}^{n - 1}(k) \setminus \mathbb{P}^{n - 2}(k)$, with 
$\mathbb{P}^{n - 2}(k) = \mathbb{P}^{n - 1} \cap \Pi_u$ (this projection is a birational morphism).  \\

Also, we know that $\Pi_u \cap \mQ_n$ is a cone $u\mQ'_{n - 2}$, in which $\mQ'_{n - 2}$ is a smooth quadric hypersurface in $\mathbb{P}^{n - 2}(k)$.  

So we have:
\begin{align*}
[\mQ_n] \ &=\ \Big[\mQ_n \setminus (\Pi_u \cap \mQ_n) \Big]\ +\ \Big[\Pi_u \cap \mQ_n\Big] \\
&= \Big[\mathbb{P}^{n - 1}(k) \setminus \mathbb{P}^{n - 2}(k)\Big]\ +\ \Big[u\mQ'_{n - 2}\Big]\\
&= \mathbb{L}^{n - 1}\ + \ (1 + \mathbb{L} \cdot [\mQ'_{n - 2}]).  
\end{align*}

By the induction hypothesis, there exists a field extension $K / k$ of degree at most $2$ for which $[\mQ'_{n - 2}]$ is in the $\mathbb{Z}[\mathbb{L}]$-module generated by $1$ and $\Spec(K)$. So by the last identity, we have that $[\mQ_n]$ is also in this module.  

The proof is done. \eop \\

\section{GQ examples}

We have a quick look at the classes of the finite orthogonal quadrangles. We fix the finite field $\F_q$. 

\subsection*{$\mQ(4,q)$-quadrangles}

Let $u$ be a (rational) point of $\mQ(4,q)$, and project $\mQ(4,q)$ onto a $\mathbb{P}^{3}(q)  \subset \mathbb{P}^4(q)$ which 
does not contain $u$. Let $\Pi_u$ be the tangent hyperplane at $u$. The projection away from $u$ yields an 
isomorphism between $\mQ(4,q) \setminus (\Pi_u \cap \mQ(4,q))$ and $\mathbb{P}^{3}(q) \setminus \mathbb{P}^{2}(q)$, with 
$\mathbb{P}^{2}(q) = \mathbb{P}^{3} \cap \Pi_u$. 

Now $\Pi_u \cap \mQ(4,q)$ is a cone $u\mC$, with $\mC$ conic in $\mathbb{P}^2(q)$. As $\mC \cong \mathbb{P}^2(q)$, we obtain

\begin{equation}
[\mQ(4,q)]\ =\ \mathbb{L}^3\ + \ 1\ +\ \mathbb{L}(\mathbb{L} + 1)\ =\ \mathbb{L}^3\ +\ \mathbb{L}^2\ +\ \mathbb{L}\ +\ 1.  
\end{equation}

\subsection*{$\mQ(5,q)$-quadrangles}

Using the same notation as in the previous cases, and performing again a projection as above, we obtain that $[\mQ(5,q)] = \mathbb{L}^4 + 1 + \mathbb{L}[\mE]$, in which $\mE = (\Pi_u \cap \mQ(5,q)) \cap \mathbb{P}^3(q)$, with $\mathbb{P}^3(q) \subset \Pi_u$ not containing $u$. We determine the class of $\mE$ in precisely the same way (projection away  from a rational point $u' \in \mE$), to obtain that $[\mE] = \mathbb{L}^2 + 1 + \mathbb{L}[\mQ]$, in which $\mQ = (\Pi_{u'} \cap \mQ(5,q)) \cap \mathbb{P}^1(q)$, with $\mathbb{P}^1(q) \subset \Pi_{u'}$ not containing $u'$; in the analysis of the $1$-dimensional case, we concluded that $[\mQ] = [\Spec(K)]$ with $K$ a field extension of degree $2$ of $k$. 

We obtain that 

\begin{equation}
[\mQ(5,q)]\ =\ \mathbb{L}^4 + \mathbb{L}^3 + \mathbb{L}^2[\Spec(K)]\ +\ \mathbb{L}\ +\ 1.
\end{equation}

\subsection*{$\mQ(3,q)$-quadrangles}

Again applying the same method as above, or simply observing that $\mQ(4,q) \cong \mathbb{P}^1(q) \times \mathbb{P}^1(q)$, we conclude that 

\begin{equation}
[\mQ(3,q)] \ =\ (\mathbb{L} + 1)^2. 
\end{equation}


 





\end{document}